\documentclass[11pt,a4paper,maplems]{article}
\usepackage{color,amsmath,amssymb,lscape,graphicx,ifthen}
\usepackage[colorlinks=true,
linkcolor=webgreen,
filecolor=webbrown,
citecolor=webgreen]{hyperref}

\definecolor{webgreen}{rgb}{0,.5,0}\definecolor{webbrown}{rgb}{.6,0,0}

\newcommand{\seqnum}[1]{\href{http://oeis.org/#1}{\underline{#1}}}

\def\pbn{\par\bigskip\noindent}

\def\psn{\par\smallskip\noindent}
\def\pn{\par\noindent}
\def\Beq{\begin{equation}}
\def\Eeq{\end{equation}}
\def\Beqarray{\begin{eqnarray}}
\def\Eeqarray{\end{eqnarray}}

\def\notequiv{\not\equiv}
\def\sspp{\,+\,}
  
\def\sspm{\,-\,}
\def\spm{\ -\ }  
\def\sspeq{\,=\,}
\def\speq{\ =\ }
\def\sspneq{\, \neq \,}
\def\sspdef{\, := \,}
\def\sspfed{\, =: \,}
\def\sspequiv{\,\equiv\,}
\def\ssppm{\,\pm\,}

\def\sspleq{\,\leq\,}
\def\sspgeq{\,\geq\,}
\def\sspgr{\,>\, }
\def\sspkl{\,<\, }
\def\spkl{\ <\ }
\def\sspin{\,\in\,}
\def\spin{\ \in\,}
\def\sspto{\,\to\,}
\def\ssptimes{\,\times\,}

\def\sspcdot{\,\cdot\,}
\def\sspapprox{\,\approx\,}

\def\eg{{\it e.g.},\,}
\def\Eg{{\it E.g.},\,}
\def\ie{{\it i.e.},\,}

\def\viz{{\it viz}\ }

\def\sspnotequiv{\,\not\equiv  \,}

\def\rpapf{{\it rpapf}\,}
\def\rpapfs{{\it rpapf}s\,}
\def\modan#1#2{{#1}\,(\text{mod}\,{#2})\,}
\def\LTvez#1#2{\left( \begin{array}{c} #1 \\ #2\end{array}\right)} 

\def\LTmatz#1#2#3#4{\left(\begin{array}{cc}#1&#2 \\ #3&#4 \end{array}\right)} 
\def\dstyle#1{$\displaystyle #1 $}
\def\Eq#1{eq.\,(#1)\,}
\def\floor#1{\left\lfloor{#1}\right\rfloor}
\def\ceil#1{\left\lceil{#1}\right\rceil}
\def\abs#1{\vert {\,#1\,} \vert}
\def\Caseszwei#1#2#3#4{\left\{ \begin{array}{ll}#1&\mbox{#2}\\ &\\#3&\mbox{#4}\end{array}\right.}
\def\range#1#2{#1,\,\count15=#1 \advance\count15 by +1 \number\count15,\,...,\,#2\,}
\def\via{{\it via}\,}

\normalbaselines
\begin{document}
\bibliographystyle{unsrt}
\vbox {\vspace{6mm}}
\begin{center}
{\Large {\bf  On Positive Integer Descartes-Steiner Curvature Quintuplets}}\\ [9mm]
Wolfdieter L a n g \footnote{ 
\href{mailto:wolfdieter.lang@partner.kit.edu}{\tt wolfdieter.lang@partner.kit.edu},\quad 
\url{http://www.itp.kit.edu/~wl}
                                          } \\[3mm]
\end{center}
\vspace{2mm}
\begin{abstract}
In Descartes' five circle problem integer curvatures (inverse radii) are considered. The positive integer curvature triple $[c_1,\, c_2,\, c_3]$ (dimensionless), with non-decreasing entries for three given mutually touching circles, leading to integer curvatures $[c_{4,-},\, c_{4,+}]$ for the two circles touching the given ones is called a {\sl Descartes-Steiner} triple. They come in two types: $[c,\,c,\,d]$ (or $[c,\,d,\,d]$) and triples with distinct entries. The first case is related to Pythagorean triples. The distinct curvature case is more involved and needs a combined representations of certain binary quadratic forms of the indefinite and definite type. The degenerate case when a straight line touches the three given touching circles can also be characterized completely. 
\par\smallskip\noindent
\end{abstract}
\section{Introduction}
\hskip 1cm 
Motivated by the section on {\sl Descartes} in \cite{Strick} where his three circle problem is discussed we consider the case of positive integer curvatures (inverse radii or bends) of three given mutually touching circles (we do not consider vanishing curvatures, \ie the straight line degenerate cases). The touching circular disks are assumed to have as their pairwise intersections only the touching point, \ie the three touching circles are outside one another, see {Figure 1}. The question is for which triple $[c_1,\,c_2,\,c_3]$, with non-decreasing positive integer entries (in  some inverse length unit), the curvatures for the two circles which touch each of the three given ones are also integer. {\sl Descartes} \cite{Descartes} showed in a letter from 1643 to Princess Elizabeth of Bohemia how to find a quadratic equation for the radius of the inner touching circle, leaving details as an exercise for a laborious calculator, and {\sl Steiner} \cite{Steiner} (not mentioning {\sl Descartes}) found a nice form for this quadratic equation for the curvature $c_4$. For {\sl Descartes}' problem see also \cite{MathWorld} , \cite{Wikipedia} and \cite{Lagarias} with the {\sl Beecroft}, {\sl Soddy}, {\sl Pedoe} and {\sl Coxeter} references). This quadratic equation is
\Beq
\left(\sum_{j=1}^4\,c_j\right)^2 \sspm 2\,\left(\sum_{j=1}^4\, c_j^2 \right)\sspeq 0\ .
\Eeq  
The two solutions for $c_4$, called $c_{4,+}$ and $c_{4,-}$, are \psn
\Beq
c_{4,\pm} \sspeq c_1\sspp c_2\sspp c_3 \ssppm 2\,\sqrt{c_1\,c_2\sspp c_1\,c_3\sspp c_2\,c_3} .
\Eeq
Vanishing $c_{4,-}$ indicates a degenerate case if a straight line touches the three given (non-degenerate) circles (see {\sl Figure 5}). $c_{4,-}$ (which can also be negative) stands for the circle of curvature $|c_{4,-}|$ circumscribing the three given ones. $c_{4,+}$ gives the curvature of the inner circle. See {\sl Fig. 2 (b)} of \cite{Lagarias} and the examples in the {\sl Note 1} below after {\sl Proposition 4}. We shall use the abbreviation $q\sspeq\sqrt{c_1\,c_2\sspp c_1\,c_3\sspp c_2\,c_3}\sspgr 0$. 
\psn
A computer program can be written for finding the primitive {\sl Descartes-Steiner} ($DS$)-triples: Select from all in`teger triples $0\sspkl c_1\sspleq c_2,\sspleq c_3\sspeq N\sspin \mathbb N$ the ones with $\gcd(c_1,\,c_2,\,c_3)\sspeq 1$ and integer $q^2$, then compute $c_{4,\pm}$. This produces a list $L_N\sspeq [c_1,\,c_2,\,c_3,c_{4,-},\,c_{4,+}, q\sspgr 0]$ for all primitive $DS$ triples with curvatures $\sspleq N$. 
In this way {\sl Tables 2,\,3,\,4} and {\sl 5} will later be found.
\psn
{\sl Section 2} gives the coordinates of the centers of the three touching circles and of the circum- and inscribed ones. A classical problem for certain areas involving these five circles is discussed. Then the definition of primitive $DS$-triples and some of their general properties are treated. 
\psn 
In {\sl Section 3} the primitive $DS$-triples of types $[c,\,c,\,d]$ and $[c,\,c,\,d]$ are studied, and mapped to primitive Pythagorean triples.
\psn
In {\sl Section 4} the distinct curvature triples come in three cases $q\sspeq c_3$, $q\sspleq c_3$, and $q\sspgr c_3$, named  ${\bf i)}$, ${\bf ii)}$, and $\bf iii)$, respectively.\pn
The curvatures for case ${\bf i)}$ are obtained from the proper solutions of the indefinite binary form $X^2\sspm 2,\,Y^2$ of discriminant $8$ representing the numbers $-s(n)$, with $s(n)\sspeq $\seqnum{A058529}$(n)$, with $n \sspgeq 2$. Such linkable OEIS \cite{OEIS} A-numbers will be given later without quotation. A uniqueness property selects in each pair of conjugate families of infinitely many solutions one with $0\sspkl X\sspkl Y$. See {\sl Proposition 6}, following many {\sl Lemmata}. 
\psn 
The curvatures for the cases $\bf ii)$ and $\bf iii)$ are obtained similarly. There two binary quadratic forms enter, one indefinite like in case $\bf i$ but with $Y\sspto \widehat Y$, the other positive definite $t^2\sspp 2\,k^2$, of discriminant $-8$, with $q\sspeq c_3\sspm k$ and $q\sspeq c_3\sspp k$, with $k\sspgr 0$, for the tow cases $\bf ii)$ and $\bf iii)$, respectively. $\widehat Y\sspeq Y\sspp k$ and $\widehat Y\sspeq Y\sspm k$, respectively. These forms have to represent simultaneously a number $-a$ and $a\sspgr 2$, respectively. For the indefinite form only one solution for each pair of conjugate solutions will qualify with the respective positive $\widehat Y$, and the definite form also has one solution with $k\sspgr 0$. The possible values for $a$ are in both cases determined, and come in three types, ${\bf a)}, {\bf b)} and {\bf c})$, depending on the proper or improper solutions of the two forms. See {\sl Proposition 7} and {\sl Proposition 8} for the two cases.
\pn
These solutions have then to be investigated for $0\sspkl X\sspkl Y$ and primitivity of $[c_1,\,c_2\, c_3]$. Examples for these investigations are given. There is no general way to find the final curvature solutions for a given candidate $a$.
\psn    
{\bf Note in passing:} This touching circle problem bears strong resemblance to {\sl G. F. Malfatti}'s circles (\cite{Mathworld2},\cite{Wiki2}, {\sl Martin} \cite{Martin}, pp. 92-96) treated in connection with {\sl Malfatti}'s problem from (1803). There the circle problem is to inscribe in a general triangle three mutually (externally) touching circles, each of which touches two of the sides of the triangle. {\sl J. Steiner} \cite{SteinerWerke} generalized this to the present five touching circle problem. 
\pbn
\section{Touching Circles and DS curvature triples}
\psn
{\bf Exercise 1: Choice of centers of the three touching circles}\psn
Given the (dimensionless) curvature triple $[c_1,\,c_2,\,c_3]$, with non-decreasing positive entries, we can choose the center of the circle with radius \dstyle{r_1\sspeq\frac{1}{c_1}} as the origin in the $(x,\,y)-$plane: $M_1\sspeq (0,\,0)$. The center of the next circle with radius \dstyle{r_2\sspeq\frac{1}{c_2}} is chosen (if necessary after a rotation) as $M_2\sspeq (r_1\sspp r_2,0)$. Then the center of the third circle with radius \dstyle{r_3\sspeq\frac{1}{c_3}} is chosen (if necessary after a reflection on the $x$-axis) in the upper half-plane and it will be $M_3\sspeq (x_3,\,y_3)$ with $y_3^2 \sspeq (r_1\sspp r_3)^2\sspm x_3^2\sspeq (r_2\sspp r_3)^2\sspm (r_1\sspp r_2\sspm  x_3)^2$. This means that
\psn
\Beqarray
x_3&\sspeq&  \frac{r_1^2\sspp r_1\,r_2\sspm r_2\,r_3\sspp r_1\,r_3}{r_1\sspp r_2} \sspeq \frac{c_1\,c_2\sspp c_1\,c_3 \sspp c_2\,c_3\sspm c_1^2}{c_1\,c_3\,(c_1\sspp c_2)} \nonumber \\&&  \sspeq \frac{q^2\sspm c_1^2}{c_1\,c_3\,(c_1\sspp c_2)}\ ,\\ 
y_3&\sspeq&  2\,\frac{\sqrt{r_1\,r_2\,r_3\,(r_1\sspp r_2\sspp r_3) }}{r_1\sspp r_2} \sspeq \frac{2\,q}{(c_1\sspp c_2)\,c_3}\sspgr\frac{1}{c_3}\sspeq r_3\, .
\Eeqarray
The polar representation will be given later in {\sl Exercise 4}. For example, if $[c_1,\,c_2,\,c_3]\sspeq [1,\,2,\,3]$, \ie \dstyle{[r_1,\,r_2,\,r_3]\sspeq \left[1,\,\frac{1}{2},\,\frac{1}{3}\right]} (in some length unit), then $q\sspeq \sqrt{11}$ and \dstyle{M_3\sspeq \left(\frac{10}{9},\,\frac{2}{3}\,\sqrt{11}\right)}. See {\sl Figure 1}.\psn 
\parbox{16cm}{
{\includegraphics[height=8cm,width=.5\linewidth]{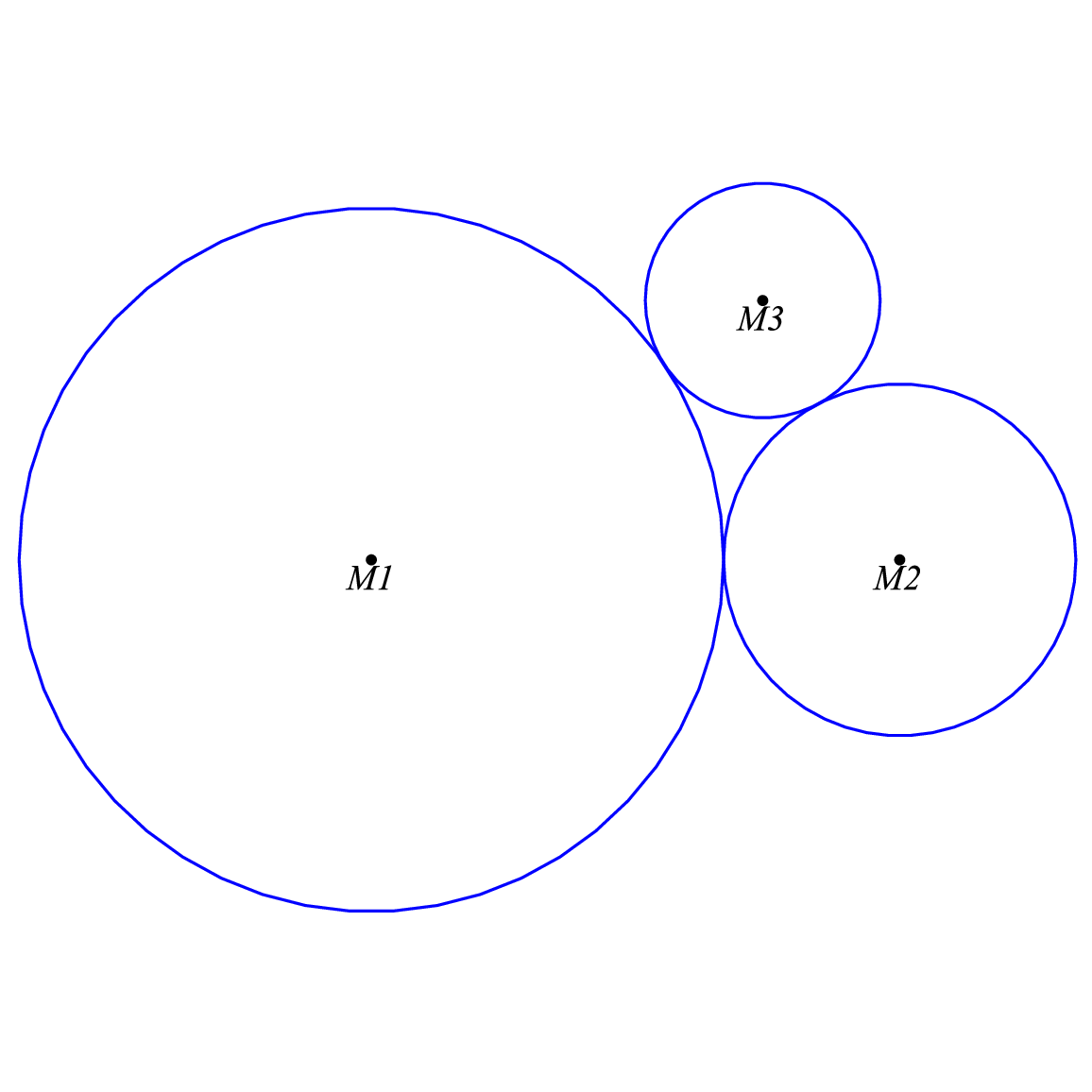}}
{\includegraphics[height=8cm,width=.5\linewidth]{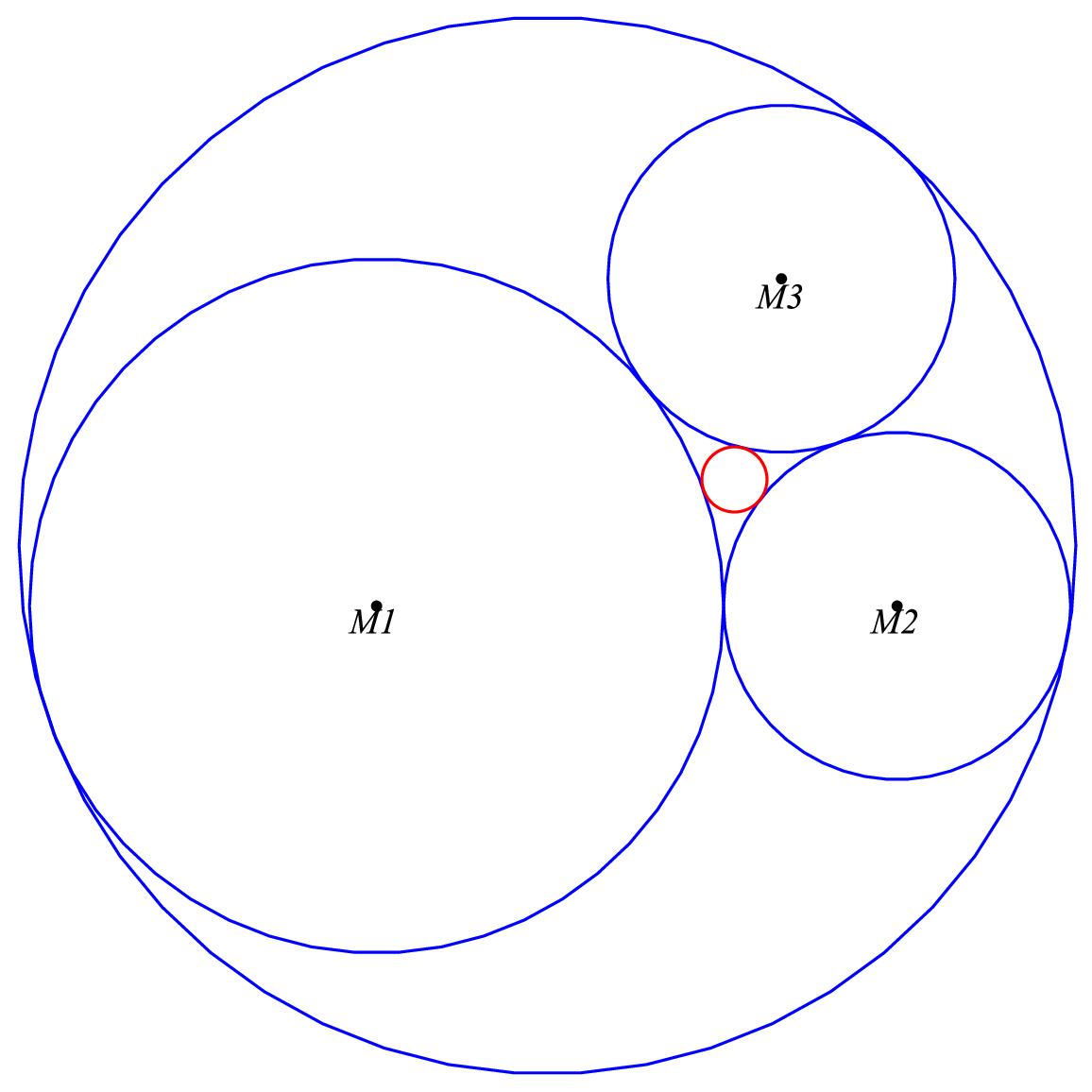}}}
\psn
{\bf Fig. 1:} Input circles for the curvature triple $[1,\, 2,\, 3]$.\ \ \ \  {\bf Fig. 2:} Five circles for the curvature triple \psn
\hskip 10cm \dstyle{\left[\frac{1}{2},\,1,\,1\right]}, and their outer and inner circles.
\psn
\pbn   
{\bf Exercise 2: Case $\bf c_{4,-} = 0$}\psn
A straight line appears as fourth touching `circle' for $c_{4,-}\sspeq 0$. In this case the three given circles with curvatures $[c_1,\,c_2,\,c_3]$, where $1\sspleq c_1\sspleq c_2\sspleq c_3$, have to satisfy the original eq. $(1)$ with only three curvatures. Therefore, the third curvature is given in terms of the two other ones by the quadratic equation $ c_3^2\sspm 2\,(c_1\sspp c_2)\,c_3\sspp (c_1\sspm c_2)^2$\sspeq 0. The result is
\Beq
c_3\sspeq c_1\sspp c_2 \sspp 2\,\sqrt{c_1\,c_2}\ .
\Eeq
The other solution with a negative square root does not qualify because $c_3 \sspgeq c_2\sspgeq c_1$, and for this solution $c_3\sspgeq c_2$ leads to $c_1\sspgeq c_2$ which is a contradiction.\pn
{\sl I.e.,} if $c_{4,-}\sspeq 0$ then \dstyle{\frac{c_3}{2}} is the sum of the arithmetic and geometric mean of the two other curvatures. \Eg $[1,\, 4,\, 9]$ has $c_{4,-}\sspeq 0$ (see {\sl Table 2}) and indeed $9 \sspeq 1\sspp 4 \sspp 2\sspcdot 2$, or \dstyle{\frac{9}{2}\sspeq \frac{5}{2}\sspp 2}.\pn
This case is later analyzed further in {\sl Exercise} $4, \bf ii)$ and {\sl Propositions} $5$ and $6$.
\psn
{\bf Exercise 3: A classical problem}\psn
Given a circle of radius $r_1$ touching two smaller touching circles of the same radius  $r_2 \sspkl r_1$, but not circumscribing them. See {\sl Figure 2}. Compute the area $F$ which is left from the circumscribing fourth touching circle of curvature $-c_{4,-}$ if the three interior circular disks are taken away. ($c_4 \sspkl 0$ because $c_1 \sspeq \frac{1}{r_1} \sspkl 4\,c_2 \sspeq \frac{4}{r_2}$, and $r_1\sspgr r_2 \sspgr \frac{r_2}{4}$). Note that in this exercise the curvatures ar not necessary integer.  
\psn
A similar problem is to find the area $F^{\prime}$ which remains after the fifth interior disk with curvature $c_{4,+}$ together with the three given circular disks is subtracted from the circumscribing circular disk.
\psn
In \cite{Strick}, p. 43, the first version of this problem is related to one of the $15$ from $150$ problems proposed by {\sl Sawaguchi Kazuyuki} in $1670$ (according to \cite{Seki} in $1671$) which he could not solve. {\sl Seki Takakazu} (also known as {\sl Seki Kowa}) \cite{Kowa} gave solutions to these $15$ problems in his book {\sl Hatsubi Samp$\overline o$} ($1674$), and they have later been reviewed by {\sl Tanaka Yoshizane} in his {\it Samp$\overline o$-meikai}.\psn
{\bf Solution:} Take as length unit ($l.u.$) the length of the small radius $r_2$. Then the given big circle has dimensionless radius \dstyle{\hat r_1\sspdef \frac{r_1}{r_2}}. The circumscribing fourth circle has dimensionless radius $\hat r\sspdef \frac{r}{r_2}$ with negative dimensionless curvature ${\hat c}_{4,-}\sspeq - \frac{1}{\hat r}$. Take $r_1\sspeq r_2\sspp a\ {l.u.} $, \ie $\hat r _1 \sspeq 1\sspp a$, with $a\sspgr 0$. Thus, \dstyle{\frac{1}{\hat r}\sspeq  -\hat{c}_{4,-}\sspeq 2\,\sqrt{1\sspp 2\,{\hat c}_1}\sspm (2\sspp {\hat c}_1)}, with ${\hat c}_1 \sspeq r_2 /r_1$. The negativity of $\hat c_{4,-}$ follows from $\hat{c}_1 \sspkl 4$, \ie  $a\sspgr -\frac{3}{4}$, which is satisfied because here $a\sspgr 0$. The result is \dstyle{ \frac{1}{\hat r} \sspeq -{\hat c}_{4,-}\sspeq 2\,\sqrt{\frac{3\sspp a}{1\sspp a}} \sspm \frac{3\sspp 2\,a}{1\sspp a}\sspgr 0}. The dimensionless area \dstyle{\hat F\sspdef \frac{F}{r_2^2}} satisfies therefore, 
\Beq
\frac{\hat F}{\pi}\sspeq \hat r^2\sspm (2\sspp {\hat r}_1^2)\sspeq \frac{1}{\hat c_{4,-}^{\,2}} \sspm \left(2\sspp \frac{1}{{\hat c}_1^{\,2}}\right)\ ,
\Eeq
which becomes, after some rewriting,
\Beq
\frac{\hat F(a)}{\pi}\sspeq  \frac{ 4\,(1\sspp a)^2\,(3\sspp 5\,a\sspp 2\,a^2)\,\sqrt{\frac{a\sspp 3}{a\sspp 1}} \sspm  2\,(1\sspp 2\,a)\,(3\sspp 4\,a \sspp 2\,a^2 \sspp 2\,a^3)}{(3\sspp 4\,a)^2}\ .
\Eeq
For the original {\sl Sawaguchi Kazuyuki} problem, according to \cite{Strick}, \dstyle{a\sspeq \frac{5}{2}} and the area $F$ was given as $120$ in some length unit $L.u.$ squared. Thus, the scaling factor $\lambda\sspeq \frac{l.u.\sspeq r_2}{L.u.}$ to be applied to the present computation is given by $120\sspeq \frac{F}{L.u.^2}\sspeq \lambda^2 \hat F\left(\frac{5}{2}\right)$, hence \dstyle{\lambda\sspeq \sqrt{\frac{120}{\hat F\left(\frac{5}{2}\right)}} \sspapprox 2.493}, due to \dstyle{{\hat F}\left(\frac{5}{2}\right)\sspeq\pi\,(\frac{1}{169}\, (196\,\sqrt{77}\sspm 681))\approx 19.312}. 
\pn
Thus, in the original problem the small radius is $\lambda L.u.$ and the radius of the large one is $\frac{7}{2}\,\lambda\,L.u.$.
\psn 
If the area of the inner touching circle with curvature \dstyle{{\hat c}_{4,+} \sspfed \frac{r_2}{r_{+}}} is also subtracted, the result becomes 
\Beq
\frac{\hat F^{\prime}(a)}{\pi}\sspeq  \frac{8\,(1\sspp a)^3\, (3\sspp 2\,a)\,\sqrt{\frac{a\sspp 3}{a\sspp 1}}}{(3\sspp 4\,a)^2} \sspm (2 \sspp (1\sspp a)^2)\ .
\Eeq 
For \dstyle{a\sspeq \frac{5}{2}} this becomes \dstyle{\hat F^{\prime}\left(\frac{5}{2}\right)\sspeq \pi\,\left(\frac{392}{169}\,\sqrt{77}\sspm \frac{57}{4}\right) \sspapprox 19.176}.
\hskip 5.5cm $\square$
\pbn
Starting with a curvature triple $[c_1,\,c_2,\,c_3]$ where $1\sspleq c_1\sspleq c_2\sspleq c_3$ (the non-degenerate circle cases), we ask for the conditions for which $c_1\,c_2\sspp c_1\,c_3  \sspp  c_2\,c_3\sspfed q^2$ with an integer $q$, and we choose $q\sspgr 0$, in fact, $q\sspgeq 3$. This is equivalent to $c_{4,\pm}$ being integers. Such curvature triples will be called {\bf Descartes-Steiner triples} or $DS$-triples. (This should not be confused with the {\sl Steiner} triples, or triple systems, from block design.) Of course, $[c_1^{\prime},\,c_2^{\prime},\,c_3^{\prime}] \sspeq [g\,c_1,\,g\,c_2,\,g\,c_3]$ with positive integer $g$ leads to  $q^{\prime}\sspeq g\,q$, and it is sufficient to study {\bf primitive  $\bf DS$-triples}.
\psn
{\bf Exercise 4: Centers of the fourth and fifth touching circle}\psn
{\bf i)\ Case $\bf c_{4,-}\sspneq 0$}\pn
The two circles which touch the given three touching circles can be found from the so-called complex {\sl Descartes} equation  \cite{Lagarias}:
\Beq
\left(\sum_{j=1}^4\,c_j\,z_j\right)^2 \sspm 2\,\left(\sum_{j=1}^4\, (c_j\,z_j)^2 \right)\sspeq 0\ ,
\Eeq
where in our convention the complex centers of the three given touching circles are $z_1\sspeq 0$, $z_2\sspeq r_1\sspp r_2$ and $z_3\sspeq x_3\sspp y_3\,i$, with $x_3$ and $y_3$ from eqs. $(3)$ and $(4)$. For the polar representation of $z_3\sspeq \rho_3\,e^{\varphi_3\,i}$ one has  \dstyle{\rho_3\sspeq \frac{q^2\sspp c_1^2}{c_1\, c_3\,(c_1 \sspp c_2)}} and \dstyle{\varphi_3\sspeq \arctan\left( \frac{2\,q\,c_1}{q^2\sspm c_1^2}\right)}. One can replace in eq. $(1)$ $c_j$ by $c_j\,z_j $, and therefore one obtains for the centers of the two touching circles $z_{4,-}$ and $z_{4,+}$
\Beq
z_{4,\pm} \sspeq \frac{1}{c_{4,\pm}}\,\left (c_1\,z_1\sspp c_2\,z_2\sspp c_3\,z_3 \ssppm 2\,\sqrt{c_1\,c_2\,z_1\,z_2\sspp c_1\,c_3\,z_1\,z_3\sspp c_2\,c_3\, z_2\,z_3}\,\right),
\Eeq
where now the square root of a complex number has to be taken. The square root will now be abbreviated as $Q$. With our choice for $z_j$, $j\sspeq 1,\,2,\,3$, one finds, for the center $M4p\sspeq (x_{4,+}\sspeq Re(z_{4,+}),\,y_{4,+}\sspeq Im(z_{4,+}))$, with  the real \dstyle{z_2\sspeq \frac{1}{c_1} \sspp \frac{1}{c_2}},
\Beqarray
x_{4,+}&\sspeq& \frac{1}{c_{4,+}}\, \left( c_2\,z_2 \sspp  c_3\,x_3 \sspp 2\,\sqrt{c_2\,c_3\,z_2\,\rho_3}\,\cos\left(\frac{\varphi_3}{2}\right)\right) \\
&\sspeq&  \frac{1}{c_{4,+}}\,  \left(\frac{c_1\sspp c_2}{c_1} \sspp c_3\,x_3 \sspp \frac{2\,q}{c_1} \right)\, , \\
y_{4,+}&\sspeq& \frac{1}{c_{4,+}}\, \left(c_3\,y_3 \sspp 2\,\sqrt{c_2\,c_3\,z_2\,\rho_3}\,\sin\left(\frac{\varphi_3}{2}\right)\right)  \\
&\sspeq& \frac{1}{c_{4,+}}\, \left(c_3\,y_3 \sspp 2 \right)\, .
\Eeqarray
The simplification occurs due to the identities \dstyle{\arctan\, x \sspeq \arccos\frac{1}{\sqrt{1\sspp x^2}}} for $x\sspgr 0$, and \pn
\dstyle{\cos\left(\frac{x}{2} \right) \sspeq \sqrt{\frac{1}{2}\,(1\sspp \cos\, x)}},  with the formula with a minus sign for \dstyle{\sin\left(\frac{x}{2} \right)}.\psn
The center $M4m\sspeq (x_{4,-},\,y_{4,-})$ is given by
\Beqarray
x_{4,-}&\sspeq&\frac{1}{c_{4,-}}\, \left( c_2\,z_2 \sspp  c_3\,x_3 \sspm 2\,\sqrt{c_2\,c_3\,z_2\,\rho_3}\,\cos\left(\frac{\varphi_3}{2}\right)\right)\, , \\
&\sspeq&\frac{1}{c_{4,-}}\,  \left(\frac{c_1\sspp c_2}{c_1} \sspp c_3\,x_3 \sspm \frac{2\,q}{c_1} \right)\, , \\
y_{4,-}&\sspeq&\frac{1}{c_{4,-}}\, \left(c_3\,y_3 \sspm 2\,\sqrt{c_2\,c_3\,z_2\,\rho_3}\,\sin\left(\frac{\varphi_3}{2}\right)\right) \\
&\sspeq& \frac{1}{c_{4,-}}\, \left(c_3\,y_3 \sspm 2 \right)\, .
\Eeqarray
The positive square roots are chosen. For example, for $[c_1,\,c_2,\,c_3] \sspeq [2,\, 3,\, 6]$ one finds $q\sspeq 6$, \dstyle{x_3\sspeq \frac{8}{15}}, \dstyle{y_3\sspeq \frac{2}{5}}, $c_{4,-}\sspeq -1$, $c_{4,+}\sspeq 23$ (see {\sl Table 2}), \dstyle{x_{4,+}\sspeq \frac{1}{23}\,\left(\frac{57}{10} \sspp 6\right)\sspeq \frac{117}{230} \sspapprox 0.509} and \dstyle{y_{4,+}\sspeq \frac{1}{23}\,\left(\frac{12}{5}\sspp 2\right)  \sspeq \frac{22}{115} \sspapprox  0.191}.  \dstyle{x_{4,-}\sspeq -\left(\frac{57}{10}\sspm 6\right) \sspeq \frac{3}{10} \sspeq 0.3} and  \dstyle{y_{4,-}\sspeq -\left(\frac{12}{5}\sspm 2\right)\sspeq -\frac{2}{5} \sspeq -0.4} \pn
For another example $[8,\,9,\,17]$ with $q\sspeq 19$ and $[c_{4,-},\,c_{4,+}]\sspeq [-4,\, 72]$ see {\sl Figure 4} and {\sl Table 2}.
\psn
{\bf ii)\ Case $\bf c_{4,-}\sspeq 0$}\pn
The center of the inner touching circle $M4p\sspeq (x_{4,+},\,y_{4,+})$ is identical with the one given in eqs. $(11)$ and $(13)$. Now a straight line will touch the three given touching circles with centers $M_1,\, M_2$ and $M_3$ chosen as explained above. To find the slope and the touching points $P$ and $Q$ with the first and second circle, respectively, one can employ a geometric construction shown in {\sl Figure 3}, found, \eg in \cite{Rechenduden}. This tangent will then automatically be tangent also to the third circle with center $M_3$, and the touching point is $R$. One draws two auxiliary circles: a circle around the origin $M1$ with radius $r_1\sspm r_2$ and a circle around \dstyle{MT\sspeq \left(\frac{r_1\sspp r_2}{2} ,0\right)}. The relevant intersection point of these two circles is $C\sspeq (x_C,\,y_C)$. By {\sl Thales} the angle $\angle(M_1,C,M_2)\sspeq \pi/2$. The touching point $P$ on the first circle  $(M1;\,r_1)$ is obtained by extending the radial line $\overline{M1,\, C}$. The tangent point $Q$ on the second circle $(M2;\,r_2)$ is obtained by the parallel to the line segment $\overline{C,\,M2}$ through $P$. The slope of the straight line through $P$ and $Q$ (the tangent) is given by $\tan(\pi\sspm \alpha)\sspeq -\tan\,\alpha$ with \dstyle{\tan\,\alpha\sspeq \frac{y_C}{r_1+r_2-x_C}}. The coordinates of $C$ as intersection point of the two auxiliary circles are
\Beq  
 x_C\sspeq  \frac{(r_1\sspm r_2)^2}{r_1\sspp r_2}\, , \hskip 1cm  y_C\sspeq  2\,\sqrt{r_1\,r_2}\,\frac{r_1\sspm r_2}{r_1\sspp r_2}\, .
\Eeq
Written in terms of curvatures this becomes
\Beq
x_C\sspeq  \frac{(c_1\sspm c_2)^2}{c_1\,c_2\,(c_1\sspp c_2)}\, , \hskip 1cm  y_C\sspeq 2\,\frac{1}{\sqrt{c_1\,c_2                                                                                                                                                                                                                                                                                                                                                                                                                                       }}\, \frac{c_2 - c_1}{c_2 + c_1} \, .
\Eeq
The coordinates of  $P$ satisfy (use an intercept theorem) $x_P \sspeq d\,\sqrt{r_1^2\sspm x_P^2}$ with \dstyle{d\sspeq \frac{x_C}{y_C}}, or  \dstyle{x_P\sspeq \frac{d\ r_1}{\sqrt{1\sspp d^2}}}  and $y_P\sspeq \sqrt{r_1^2\sspm x_P^2}$, which yields
\Beq  
 x_P\sspeq  r_1\,\frac{r_1\sspm r_2}{r_1\sspp r_2}\, , \hskip 1cm  y_P\sspeq 2\, \sqrt{r_1\,r_2}\,\frac{r_1}{r_1+r_2} \, .
\Eeq
In terms of curvatures this becomes
\Beq
x_P\sspeq \, \frac{c_2 - c_1}{c_1\,(c_2 \sspp c_1)} , \hskip 1cm  y_P\sspeq \frac{2}{\sqrt{c_1\,c_2}}\, \frac{c_2}{c_2 + c_1}\, .
\Eeq
The tangent through $P$ satisfies $y\sspeq -\tan(\alpha)\, x \sspp c$ \ie \dstyle{c\sspeq y_P\sspp \tan(\alpha)\, x_P\sspeq \frac{r_1\,(r_1\sspp r_2)}{2\,\sqrt{r_1\,r_2}}\sspeq \frac{1}{2\,\sqrt{c_1\,c_2}}\,\frac{c_1\sspp c_2}{c_1}}. With \dstyle{\tan\,\alpha\sspeq \frac{r_1\sspm r_2}{2\, \sqrt{r_1\,r_2}}\sspeq \frac{1}{2\,\sqrt{c_1\,c_2}}\,(c_2\sspm c_1)} this becomes
\Beq  
 y \sspeq \frac{\sqrt{r_1\,r_2}}{2\,r_2}\,\left (-\frac{r_1-r_2}{r_1}\, x \sspp  r_1\sspp r_2\right)\ .
\Eeq
In terms of curvatures this becomes 
\Beq  
 y \sspeq \frac{1}{2\,\sqrt{c_1\,c_2}}\,\left (-(c_2\sspm c_1)\,x \sspp  \frac{c_1\sspp c_2}{c_1}\right)\ .
\Eeq
The coordinates of the tangent point $Q$ are then obtained from \dstyle{\sin\,\beta\sspeq \frac{y_P}{r_1}\sspeq \frac{y_Q}{r_2}} and  $x_Q\sspeq r_1\sspp r_2\sspp \sqrt{r_2^2\sspm y_Q^2} $
\Beq 
x_Q\sspeq  r_1\,\frac{r_1\sspp 3\,r_2}{r_1\sspp r_2}\, , \hskip 1cm  y_Q\sspeq 2\, \sqrt{r_1\,r_2}\,\frac{r_2}{r_1\sspp r_2} \, .
\Eeq 
In terms of curvatures this becomes
\Beq 
x_Q\sspeq\, \frac{c_2\sspp 3\,c_1}{c_1\,(c_2\sspp c_1)}, \hskip 1cm  y_Q\sspeq \frac{2}{\sqrt{c_1\,c_2}}\,\frac{c_1}{c_2 \sspp c_1}\, .
\Eeq 
The coordinates of the touching point $R$ on the third circle $(M3;\,r_3)$ satisfy \dstyle{y_R\sspeq y_3 \sspp \frac{r_3}{r_1}\,y_P} and \dstyle{x_R\sspeq x_3 \sspp \sqrt{r_3^2\sspm (y_R\sspm y_3)^2}\sspeq x_3\sspp r_3\,\frac{r_1\sspm r_2}{r_1\sspp r_2 }}, hence, with {\it eqs.} $(3)$ and $(4)$,
\Beq 
x_R\sspeq r_1 + 2\,\frac{r_1\sspm r_2} {r_1\sspp r_2}\,r_3\, , \hskip 1cm  y_R\sspeq 2\, \frac{\sqrt{r_1\,r_2\,r_3}}{r_1\sspp r_2}\,(\sqrt{r_1\sspp r_2\sspp r_3}\sspp \sqrt{r_3}\,) \, ,
\Eeq 
where $r_3$ is from eq. (5)  
$ r_3 \sspeq  r_1\,r_2/(r_1\sspp r_2 \sspp 2\,\sqrt{r_1\,r_2})$.
\pn 
In terms of the curvatures \Eq{27} this becomes
\Beq
 x_R\sspeq \frac{1}{c_1}\sspp 2\,\frac{c_2\sspm c_1}{c_2\sspp c_1} \, \frac{1}{c_3}   \, , \hskip 1cm  y_R\sspeq  \frac{2}{c_3\,(c_1\sspp c_2)}\, \left(q\sspp \sqrt{c_1\,c_2}  \right)  \, ,
\Eeq
where $c_3$ is given by eq. (5). 
\pn
For the example $[c_1,\,c_2,\,c_3]\sspeq [1,\,4,\,9]$ we have $q\sspeq 7$ with $c_{4,-}\sspeq 0$ and $c_{4,+}\sspeq 28$ {\sl see Table 2} and {\sl Table 6}. $M3\sspeq \left(\frac{16}{15},\frac{14}{45}\right)\sspapprox (1.06,\,0.31)$, $\tan\,\alpha\sspeq \frac{3}{4}$, $\sin\,\beta\sspeq \frac{4}{5}$, $c\sspeq \frac{5}{4}\sspeq 1.25$, $C\sspeq \left( \frac{9}{20} ,\,\frac{3}{5}\right)\sspeq (0.45,\,0.6)$, $P\sspeq \left( \frac{3}{5} ,\,\frac{4}{5}\right)\sspeq (0.6 ,\,0.8)$, $ Q \sspeq \left( \frac{7}{5} ,\,\frac{1}{5}\right)\sspeq (1.4 ,\,0.2)$, $ R \sspeq \left( \frac{17}{15} ,\,\frac{2}{5}\right)\sspapprox (1.13,\,0.4)$.
\pbn
\parbox{16cm}{
{\includegraphics[height=8cm,width=.5\linewidth]{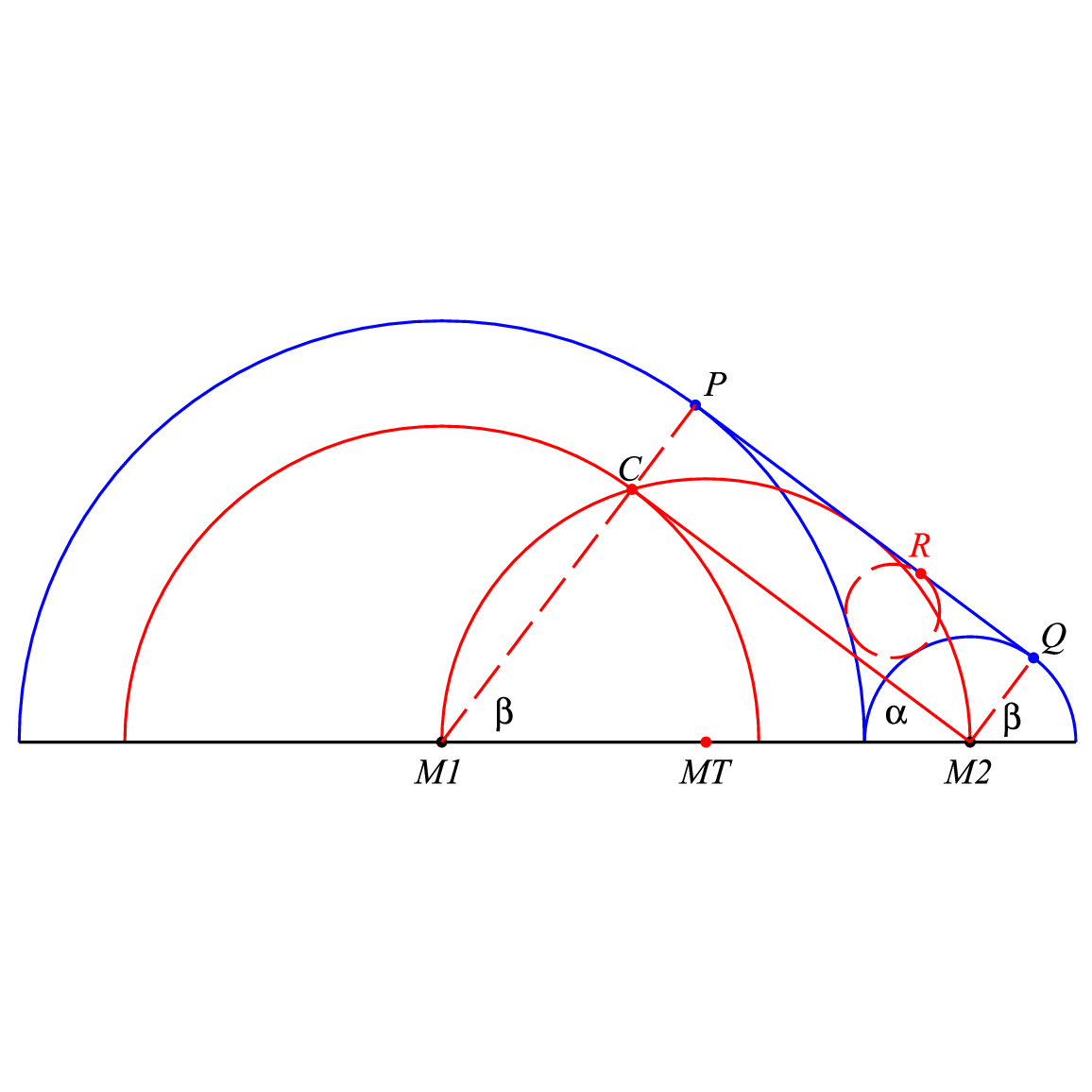}}
{\includegraphics[height=8cm,width=.5\linewidth]{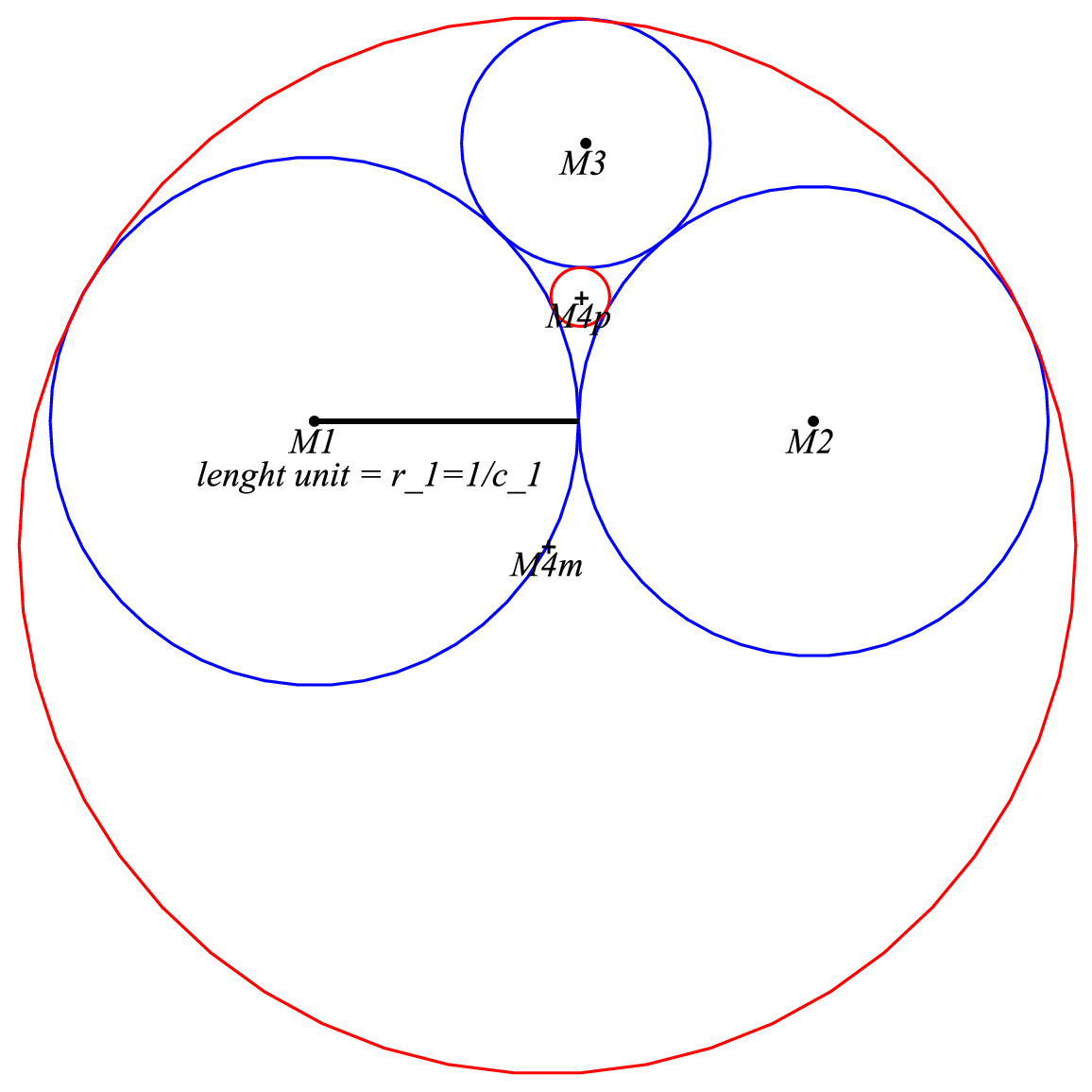}}
}
\psn
{\bf Fig. 3:} Geometric construction of tangent points \pn
$P$, $Q$ and $R$.\hskip 7cm   {\bf Fig. 4:} Five touching circles for $[8,\,9,\,17]$.
\psn
\pbn   
{\bf Definition 1: Primitive ${\bf DS}$-triples}\psn
{\it An integer curvature triple $ [c_1,\,c_2,\,c_3]$ with $0\sspkl c_1\sspleq c_2\sspleq c_3$, and $gcd(c_1,\,c_2,\,c_3)\sspeq 1 $ is called a primitive $DS$-triple if both curvatures $c_{4,+}$ and $c_{4,-}$ given in eq. ($2$) are integers.}
\psn
{\bf Proposition 1:} {\it There is no $DS$ triple of the type $[c,\,c,\,c]$, $c\sspin {\mathbb N}$}.\psn
{\bf Proof}: $\sqrt{3}$ is not an integer number.\psn
Primitive $DS$-triples obey a certain parity rule stated as follows.\psn
{\bf Proposition 2: Parity of primitive $\bf DS$-triples}\psn
{\it Two cases are possible for primitive $DS$-triples $ [c_1,\,c_2,\,c_3]$:\psn
{\bf i)} Two of the curvatures are even and the other is odd. Then the sum of the even ones is congruent to $0\,(mod\, 4)$.\psn
{\bf ii)} Two of the  curvatures are odd and the other is even. Then the sum of the odd ones is congruent to $2\,(mod\, 4)$.}
\psn
Before proceeding with the elementary proof some examples are given.\psn
{\bf Example 1:} 
\psn 
Case {\bf i)}: $[2,\,3,\,6]$  with $q\sspeq 6$ has $2\sspp 6\sspeq 8\sspequiv 0\, (mod\, 4)$. $[2,\,2,\,3]$  with $q\sspeq 4$ has $2\sspp 2\sspeq 4\sspequiv 0\, (mod\, 4)$. \psn
Case {\bf ii)}: $[2,\,3,\,15]$  with $q\sspeq 9$ has $3\sspp 15\sspeq 18\sspequiv 2\, (mod\, 4)$. $[25,\,25,\,48]$  with $q\sspeq 55$ has $25 \sspp 25\sspeq 50\sspequiv 2\, (mod\, 4)$. \psn
{\bf Proof:} Primitivity forbids that all three curvatures are even. They also cannot be all odd because always $q^2\sspequiv 0$ or $1\, (mod\, 4)$ ($q$ is even or odd) but then $q^2$ would become $3\,(mod\, 4)$. This leaves the two cases $\bf i)$ and $\bf ii)$. In the first case, with, say, $c_1\sspeq 2\,C_1\sspp 1$, $c_2\sspeq 2\,C_2$ and $c_3\sspeq 2\,C_3$ one finds $q^2\sspeq 4\,(C_1\,C_2 \sspp C_1\,C_3\sspp C_2\,C_3)\sspp 2\,(C_2\sspp C_3)$ but $q^2\sspequiv 0$ or $1\, (mod\, 4)$, hence $c_2\sspp c_3\sspequiv 0\,(mod\,4)$.  In the second case with, say, $c_1\sspeq 2\, C_1$, $c_2\sspeq 2\, C_2\sspp 1$ and $c_3\sspeq 2\,C_3\sspp 1$ one has $q^2\sspeq 4\,(C_1\,C_2 \sspp C_1\,C_3\sspp C_2\,C_3\sspp C_1)\sspp 2\,(C_2\sspp C_3)\sspp 1$, which is odd, therefore $q^2\sspequiv 1\,(mod\, 4)$. Thus, $C_2\sspp C_3$ has to be even and $c_2\sspp c_3 \sspequiv 2 (mod\, 4) $. \hskip  9cm $\square$
\pbn 
{\bf Proposition 3:} {\it For primitive $DS$-triples one has $gcd(c_1\sspp c_2,\, c_1\sspp c_3,\,c_2\sspp c_3)\sspeq 1$.}
\psn
{\bf Proof:} Note that this $gcd$ certainly has to be odd because otherwise this would require all curvatures to be either even or odd. But assuming that  $ c_1\sspp c_2\sspeq g\,X$, $ c_1\sspp c_3\sspeq g\,Y$ and $ c_2\sspp c_3\sspeq g\,Z$, with odd $g\sspgeq 3$, implies $2\,c_1\sspeq g\,(X\sspp Y\sspm Z)\sspfed g\,a_1$, $2\,c_2\sspeq g\,(X\sspm Y\sspp Z)\sspfed g\,a_2$  and $2\,c_3\sspeq g\,(-X\sspp Y\sspp Z)\sspfed g\,a_3$. Hence $a_j\sspeq 2\, A_j$, for $j\sspeq 1,\,2,\,3$, and then $c_j\sspeq g\,A_j$, but this contradicts the primitivity, hence $g$ has to be $1$. \hskip 15.1cm $\square$\psn
It is clear that also the converse is true. If for $x\sspdef c_1 \sspp c_2$,  $y\sspdef c_1 \sspp c_3 $ and $z\sspdef c_2 \sspp c_3$  $gcd(x,\,y,\,z)\sspeq 1$ holds, then also $gcd(c_1,\,c_2,\,c_3)\sspeq 1$. Assume that $gcd(c_1,\,c_2,\,c_3)\sspeq g\sspgr 1$. This contradicts $gcd(x,\,y,\,z)\sspeq 1$.
\pbn
\section{DS-triples with curvatures [c, c, d] or [c, d, d]}
For the next {\sl Proposition} we define the following set of ordered pairs $\cal M$.\psn
{\bf Definition 2: Set of ordered pairs $\bf (n,\,m)$ }\psn
{\it \Beq
{\cal M}\sspeq \left\{ (n,\,m)\sspin \mathbb N \times \mathbb N \,{\biggm \vert } \, n\sspgr m,\ gcd(n,\,m)\sspeq 1, (-1)^{n+m}\sspeq -1 \right\}\ .
\Eeq }
The last condition means that $n$ and $m$ have opposite parity. For $n\sspgeq 2$ there are $\varphi(2\,n)/2\sspeq$\seqnum{A055034}$(n)$ values for $m$, where $\varphi$ is {\sl Euler}'s totient function. For the characteristic triangle of $\cal M$ see \seqnum{A249866}.
\psn
{\bf Proposition 4: Primitive ${\bf DS}$-triples of type $\bf [c,\,c,\,d]$ or $\bf [c,\,d,\,d]$}\psn 
{\it Primitive $DS$-triples of the type $[c,\,c,\,d]$ or $[c,\,d,\,d]$ ($c \sspkl d$) come in two types, parametrized by $(n,\,m)\sspin {\cal M}$.\pn
Type {\bf I}: Either  $[(n-m)^2,\,(n-m)^2,2\,\,n\,m]$ or  $[2\,n\,m,\,(n-m)^2,\,(n-m)^2]$, depending on $n^2\sspp m^2\sspkl 4\,n\,m$ or $n^2\sspp m^2\sspgr 4\,n\,m$, respectively. In both cases $q \sspeq n^2\sspm m^2$. $c_{4,+}\sspeq 2\,n\,(2\,n\sspm m)$ and $c_{4,-}\sspeq 2\,m\,(2\,m\sspm n)$.\psn
Type {\bf II}: Either  $[2\,m^2,\,2\,m^2,\,n^2\sspm m^2 ]$ or  $[n^2\sspm m^2,\, 2\,m^2,\,2\,m^2]$, depending on $\sqrt{3}\,m\sspkl n $ or $\sqrt{3}\,m\sspgr n $, respectively. In both cases $q \sspeq 2\,n\,m$. $c_{4,+}\sspeq (n \sspp m)\,(n\sspp 3\,m)$ and $c_{4,-}\sspeq  (n \sspm m)\,(n\sspm 3\,m) $.} 
\psn
{\bf Proof:} This problem reduces to primitive Pythagorean triples as follows.
 Consider a $DS$-triple of the type $[c,\,c,\,d]$, where $1\sspleq c\sspkl d$. Now $q^2 \sspeq c^2\sspp 2\,c\,d\sspeq (c\sspp d)^2 \sspm d^2 $, \ie 
\Beq
(c\sspp d)^2  \sspeq q^2\sspp d^2\sspfed z^2, 
\Eeq
where we choose $z\sspeq c\sspp d$. This means that we have to consider the Pythagorean triple $(q,\,d,\,z)$.\pn
Similarly, a $DS$-triple of the type $[c,\,d,\,d]$, where $1\sspleq c\sspkl d$, leads to the Pythagorean triple $(q,\, c,\,z)$, with $(c\sspp d)^2 \sspeq q^2\sspp c^2\sspfed z^2$, with $z\sspeq d\sspp c$ like in the $[c,\,c,\,d]$ case.
\pn 
A Pythagorean triple $(x,\,y,\,z)$ is called primitive if the three positive integers $x,\, y$ and $z$ are pairwise relatively prime (see \eg  \cite{NivenZuckermanMontgomery}). 
\psn
We  first prove 
\psn
{\bf Lemma 1:} {\it  Each primitive $DS$-triple $[c,\,c,\,d]$ leads to a primitive Pythagorean triple $(q,\,d,\,z)$.\pn
Similarly, each primitive $DS$-triple $[c,\,d,\,d]$ leads to a primitive Pythagorean triple $(q,\,c,\,z)$. Here $z\sspdef c\sspp d $}. 
\psn
{\bf Proof:} It is sufficient to prove the $[c,\,c,\,d]$ case. The proof of other case runs similarly. \pn
With $\gcd(c,\,d) \sspeq 1$ one finds, by assuming that $\gcd(q,d)\sspeq g > 1$, that \dstyle{c^2 \sspeq g\,(g\,Q^2 \sspm 2\,c\,D)} with $q\sspeq g\,Q$ and $d\sspeq g\,D$. Thus, $g|c^2$, hence $c\sspeq g\,C$, in contradiction to $\gcd(c,\,d) \sspeq 1$. Therefore $g\sspeq 1$. It is clear that $\gcd(d,z\sspeq c+d)\sspeq 1$. Also $\gcd(q,\,z)\sspeq 1$ because $\gcd(q^2,z)\sspeq \gcd(c\,(c + 2\,d),\,z)$  and $\gcd(c,\,z)\sspeq d(c,\,d)\sspeq 1$ as well as $\gcd(c+2\,d,\,z) \sspeq 1$ due to $\gcd(c+2\,d,\,z) \sspeq \gcd(c+d,\,d)$, hence $\gcd(q^2,\,z)\sspeq \gcd(q,\,z) \sspeq 1$. \hskip 12.3cm$\square$ 
\pn
Note that it is sufficient to prove $\gcd(q,\, d)\sspeq 1$ because for Pythagorean triples $(q,\,d,\,z)$ this is equivalent to pairwise relatively prime $q,\, d$ and $z$. This fact is used in  \cite{HardyWright}, Theorem 225, p. 190. \pn 
\psn 
That each primitive Pythagorean triple leads to two primitive $DS$-triples, will be shown later as a {\sl Corollary 1} to this {\sl Proposition}.
\psn  
Using a well known theorem on primitive Pythagorean triples with $q^2 \sspp d^2 \sspeq z^2$  (see \eg \cite{NivenZuckermanMontgomery}, Theorem 5.5, p. 232, or \cite{HardyWright}, Theorem 225, p. 190) one proves that the structure of primitive $DS$-triples come in the two claimed types $I$ or $II$.
\pn
The primitive Pythagorean triples $(q,\,d,\,z)$ are characterized by the elements of the set $\cal M$ from {\sl Definition} $2$ either as case $\bf i)$: $q\sspeq n^2\sspm m^2$, $d\sspeq 2\,n\,m$ and $z\sspeq n^2\sspp m^2$ ($d$ is even, and then $q$ is odd) or, interchanging $q$ and $d$, as case $\bf ii)$: $q\sspeq 2\,n\,m$, $d\sspeq n^2\sspm m^2$ and $z\sspeq n^2\sspp m^2$ ($d$ is odd, and then $q$ is even). In both cases $c\sspeq z\sspm d$.\pn
For $c\sspkl d$, case $\bf i)$ needs $ (n\sspm m )^2 \sspkl 2\,n\,m$ \ie $n^2\sspp m^2\sspkl 4\,n\,m$, and the $DS$-triple is $[(n\sspm m )^2,\,(n\sspm m )^2,\,2\,n\,m]$. Case $\bf ii)$ needs $2\,m^2 \sspkl n^2\sspm m^2$ \ie $ \sqrt{3}\,m\sspkl n$, and the $DS$-triple is 
$[2\,m^2,\,2\,m^2,\,n^2\sspm m^2]$.\pn
Similarly, the primitive Pythagorean triples $(q,\,c,\,z)$ with case $\bf i)$: $q\sspeq n^2\sspm m^2$, $c\sspeq 2\,n\,m$ and $z\sspeq n^2\sspp m^2$, leads for $c\sspkl d$, \ie $n^2 \sspp m^2 \sspgr 4\,n\,m$, to the $DS$-triple $[2\,n\,m,\,(n\sspm m)^2,\,(n\sspm m)^2]$. The case $\bf ii)$: $q\sspeq 2\,n\,m$, $c\sspeq n^2\sspm m^2$ and $z\sspeq n^2\sspp m^2$, leads for $c\sspkl d$, \ie $n\sspkl \sqrt{3}\,m$, to the $DS$-triple $[n^2\sspm m^2,\,2\,m^2,\,2\,m^2]$.
\psn
What is left is to prove that these $DS$-triples, derived from primitive Pythagorean triples, are also primitive.\pn
For both primitive Pythagorean triples of type {\bf i)} one writes  $\gcd((n\sspm m)^2,\, 2\,n\,m) \sspeq \gcd( 2\,n\,m,\, n^2 \sspp m^2)\sspeq \gcd(n\,m,\, n^2 \sspp m^2)$. In the last step we used that $n^2 \sspp m^2$ is odd from {\sl Definition} $2$. Now $n^2\sspp m^2$ cannot be divisible by a factor $f_n\sspgr 1$ of $n$ because otherwise $m^2$, hence $m$, would also be divisible by $f_n$ but this contradicts the primitivity of the Pythagorean triple guaranteed by $gcd(n,\,m)\sspeq 1$. Similarly a factor $f_m\sspgr 1$ of $m$ cannot divide $n^2\sspp m^2$. Hence $\gcd((n\sspm m)^2,\, 2\,n\,m)\sspeq 1$.\pn
For both primitive Pythagorean triples of type {\bf ii)} one writes
$\gcd(n^2 \sspm m^2,\, 2\,m^2) \sspeq \gcd(n^2 \sspm m^2,\, m^2)$ (because 
$n^2 \sspm m^2$ is odd) $\sspeq \gcd(m^2,\, n^2) \sspeq 1$ because $\gcd(m,\, n) \sspeq 1.$
\psn
Therefore case {\bf I} of the Proposition is given by case $\bf i)$ from the primitive Pythagorean triples $(q,\,d,\,z)$ and by case $\bf i)$ from the primitive Pythagorean triples $(q,\,c,\,z)$, and case {\bf II} is given by the cases $\bf ii)$ of these two primitive Pythagorean triples.\pn 
Note that the parity of the primitive $DS$-triples are in accordance with {\sl Proposition}  $2$.\hskip 3cm $\square$ 
\pbn
{\bf Corollary 1:} {\it Each primitive Pythagorean triple leads to two primitive $DS$-triples, either both of type $[c,\,c,\,d]$  or one of type  $[c,\,c,\,d]$ and the other of type $[c,\,d,\,d]$. No pair of triples of the type $[c,\,d,\,d]$ is possible.}
\psn
This follows from the the conditions for cases $\bf I$ and $\bf II$.
A pair of primitive $DS$-triples $[c,\,d,\,d]$ would require $n^2\sspp m^2 \sspgr 4\,n\,m\sspgr 4\,m^2$, \ie $n^2\sspgr3\,m^2$, contradicting  $3\,m^2\sspgr n^2$. 
\psn
{\bf Example 2: Primitive $\bf DS$-triple  $\bf [c,\,d,\,d] $, and primitive Pythagorean triple $\bf (q,\,c,\,c \sspp d)$} \psn
The primitive $DS$-triple $[44,\,81,\,81]$ is of type $\bf I$ because $n\,m = 22$ and $n\sspm m \sspeq 9$, hence $n\sspeq 11$ and $m\sspeq 2$, and $q\sspeq 117\sspeq n^2\sspm m^2$ with the corresponding primitive Pythagorean triple $(117,\, 44, \,125)$. This satisfies indeed the $[c,\,d,\,d]$ type $\bf I$ inequality $n^2 \sspp m^2 \sspeq 125 \sspkl 88\sspeq 4\cdot 11\cdot 2$. Also $c_{4,-}\sspeq -28\sspeq 4\,(4\sspm 11)$ and $c_{4,+}\sspeq 440\sspeq 22\,(22\sspm 2)$.
\pn
Note that this primitive Pythagorean triple $(117, 44, 125)$ leads, because of the symmetry $(44,\, 117,\, 125)$, also to the type $\bf II$ primitive $DS$-triple of type $[c,\,c,\,d]$ namely $[8,\,8, 117]$ with $q\sspeq 2\,n\,m\sspeq 44$, $c_{4,-}\sspeq 45$ and $c_{4,+}\sspeq 221$. The inequality is here $2\,\sqrt{3}\sspkl 11 $. See {\sl Table 1} for both $DS$-triples with $[n,\, m]\sspeq [11,\,2]$.
\psn
{\bf Note 1:\  Special instances of Proposition 4}
\psn 
For type {\bf I} consider $(n,m\sspeq n-1)$, $n\sspgeq 2$, which is indeed an element  of $\cal M$ ($gcd(n,n-1)\sspeq 1$ is easy to prove by contradiction). The first inequality applies because  $2,n\,(n\sspm 1)\sspm 1 \sspgr 3\sspgr 0$, hence the  primitive $DS$-triple is $[1,\,1,\,4\,T(n-1)]$, where $n\sspeq 2$, with the triangular number $T(n-1)\sspeq$ \seqnum{A000217}$(n-1)$. Here $q\sspeq \sqrt{1\sspp 8\,T(n-1)} \sspeq 2\,n\,\sspm 1 \sspeq n^2\sspm (n-1)^2$. $c_{4,+}\sspeq 4\,T(n)$ and $c_{4,-}\sspeq 4\,T(n-2)$.  These are the triples $[1,\,1,\,4]$, $[1,\,1,\,12]$, $[1,\,1,\,24]$, $[1,\,1,\,40]$ ... with $q$ values $3,\,5,\, 7,\,9\,$ ... and $[c_{4,-},\,c_{4,+}]$ pairs $[0,\,12]$, $[4,\,24]$, $[12,\,40]$, $[24,\,60]$, ... .
\pn
For type {\bf II} take $m\sspeq 1$ , with $\sqrt{3} \sspkl n$, leading to primitive triples $[2,\,2,\,4\,k^2\sspm 1]$, if $n\sspeq 2\,k$ and $k\sspgeq 1$. Here $q\sspeq 2\,n\sspeq 4\,k$, $c_{4,+}\sspeq 4\,k\,(k\sspp 2)\sspp 3\sspeq (2\,k\sspp 1)\,(2\,k\sspp 3)$ and $c_{4,-}\sspeq (2\,k\sspm 1)\,(2\,k\sspm 3)$. This is $[2,\,2,\,3]$, $[2,\,2,\,15]$, $[2,\,2,\,35]$, ... with $q$ values $4,\,8, 12,\, ...$ and $[c_{4,-},c_{4,+}]$ pairs $[-1,\,15]$, $[3,\,35]$, $[15,\,63]$, ...\, .
\psn
The touching circles for $[1,\,1,\,4]$ with $[c_{4,-},\,c_{4,+}]\sspeq [0,\, 12]$, and for $[2,\,2,\,3]$ with $[c_{4,-},\,c_{4,+}]\sspeq [-1,\,15]$ are shown in {\sl Figure} $5$ and $6$, respectively.
\psn
The possible values for $z\sspeq c\sspp d$ give the length of the (dimensionless) hypotenuses of primitive Pythagorean triangles. They are given as a number triangle in \seqnum{A222946}, and in ordered form, with multiplicity, in \seqnum{A020882}. This is the $z$-sequence $5,\,13,\,17,\,25,\, 29,\, 37,\, 41,\, 53,\,61,\,65,\, 65,\, 73,\, 85,\, 85,\, 89,\, 101,\, ...$. In {\sl Table 1} one can find for the corresponding $(n,\,m)$ values the primitive  $DS$-triples of {\sl Proposition 4}.\pbn
\pbn
\parbox{16cm}{
{\includegraphics[height=8cm,width=.5\linewidth]{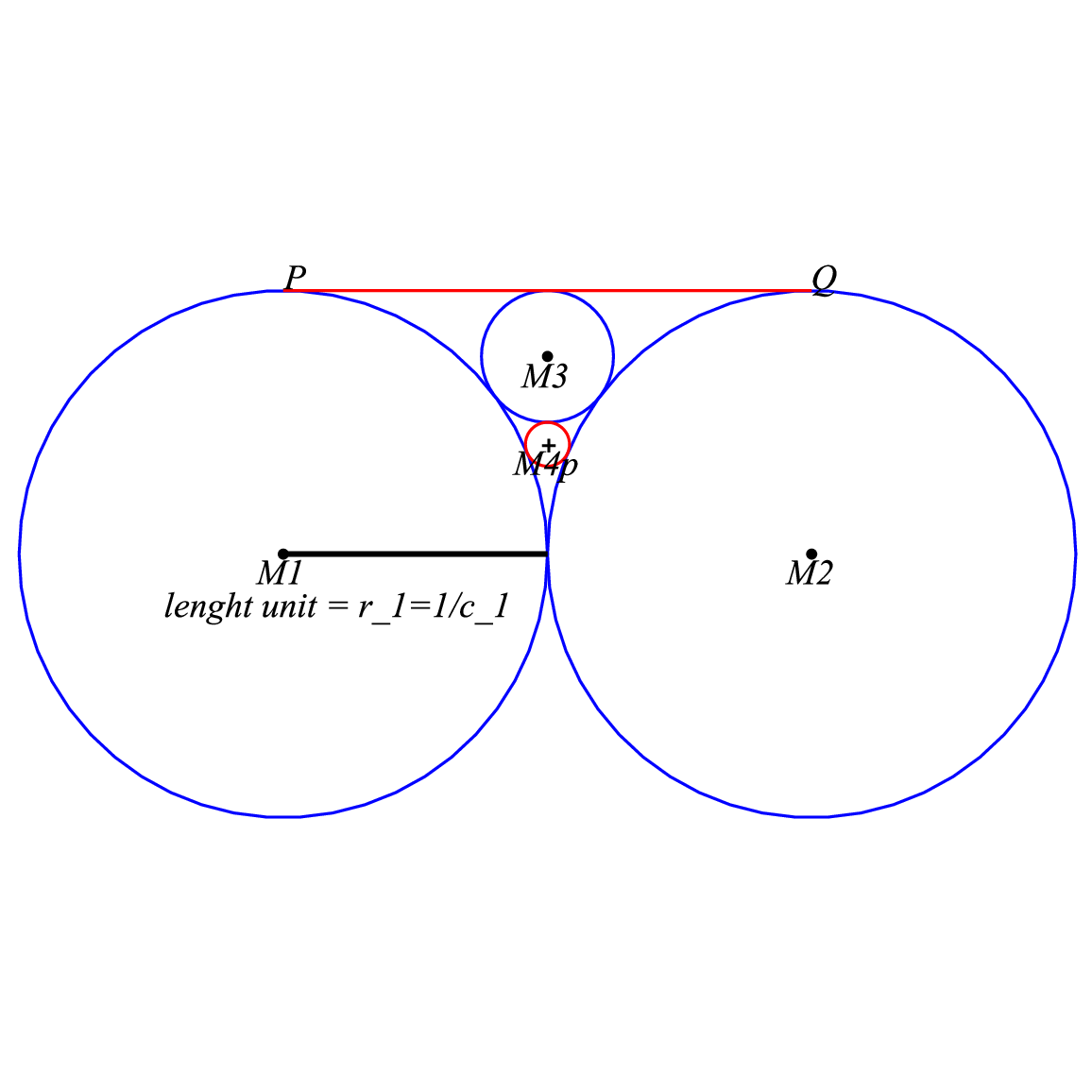}}
{\includegraphics[height=8cm,width=.5\linewidth]{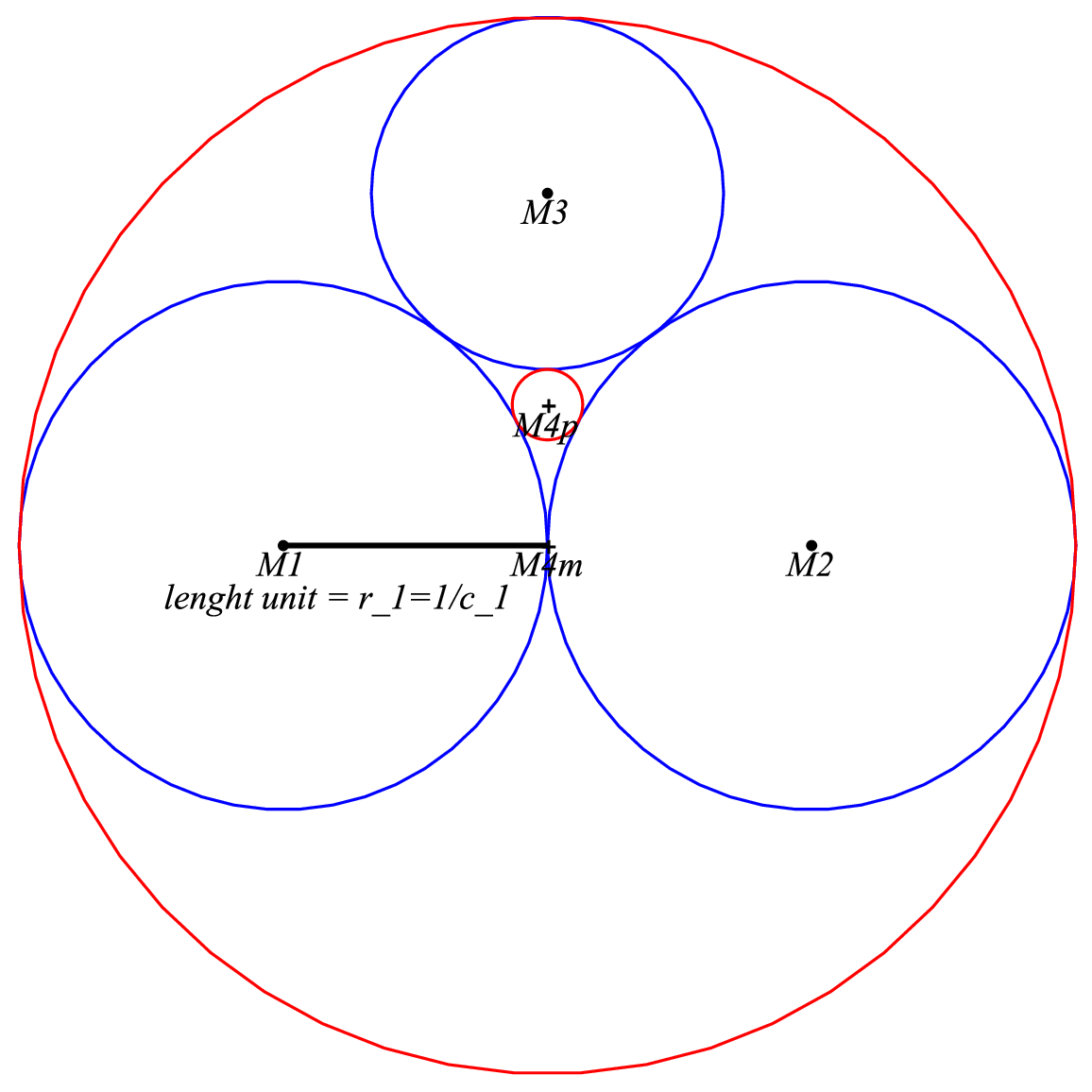}}
}
\psn
{\bf Fig. 5:} Four touching circles and a straight line for $[1,\, 1,\, 4]$.\ \ \ \  {\bf Fig. 6:} Five touching circles for $[2,\,2,\,3]$.
\pbn  
The next topic is to find all primitive $DS$-triples which belong to the special degenerated case when $c_{4,-}\sspeq 0$, \ie when a straight line touches the three touching circles. In order for this to happen there is a condition on the curvature of the third given circle in terms of the curvatures of the other two ones. \pn 
In {\sl Exercise} $2$, \Eq{5}, it was shown that $c_3\sspeq c_1\sspp c_2 \sspp 2\,\sqrt{c_1\,c_2}$\,.
This solution for $c_3$ is compatible with $\sum_{j=1}^3\,c_j^2 \sspeq 2\,q^2$, obtained from $c_{4,-}\sspeq 0$, with $q$ the square root defined after \Eq{2}. 
\pn
Due to the {\it AGM} inequality \cite{WikiAGM} in the case of two variables $4\,\sqrt{c_1\,c_2} \sspleq c_3$ (of course, $4\,\sqrt{c_1\,c_2} \sspleq 2\,(c_1\sspp c_2)\sspleq 4\,c_3)$. For the $DS$-triple $[1,\,4,\,9]$ (see {\sl Tables} {\sl 2} and {\sl 6}) one obtains $8\sspleq 9$).\pn
\pn
Because of the homogeneity with degree $1$ of the $c_3$ equation, a common factor of $c_1$ and $c_2$ will also multiply $c_3$. Therefore one imposes $gcd(c_1,\,c_2)\sspeq 1$. Then the only such $DS$-triple with $ c_1\sspeq c_2$ is $[1,\, 1,\, 4]$ with $q\sspeq 3$ and $c_{4,+}\sspeq 12$ (see {\sl Tables} {\sl 1} and {\sl 6}). All other solutions for $1\sspleq c_1\sspkl c_2$ with integer $\sqrt{c_1\,c_2}$, \ie $c_1\,c_2\sspeq h^2$, with integer $h\sspgr 0$, are obtained from $c_1\sspeq M^2$ and $c_2\sspeq N^2$ with positive integers $M$ and $N\sspgr M$. This follows from the fact that $c_1$ and $c_2$ have no common factor $\sspgr 1$, so each curvature has to be a square by itself. This is summarized in the following {\sl Proposition}.
\psn 
{\bf Proposition 5: All primitive $\bf DS$-triples with $\bf c_{4,-}\sspeq 0$}
\psn
{\it The primitive $DS$-triples which satisfy $c_{4,-}\sspeq 0$ are given by 
\Beq
[M^2,\, N^2,\, (M \sspp N)^2],\, with\, \gcd(M,\,N)\sspeq 1,\, 1\sspleq M\sspleq N,\,and\, 4\,q\sspeq c_{4,+} \sspeq 2\,(M^2\sspp N^2\sspp (M\sspp N)^2)\,.
\Eeq} 
{\bf Proof:}
\pn
Obviously $\gcd(M,\,N)\sspeq 1$ follows from $gcd(c_1,\,c_2)\sspeq 1$. The case $c_1\sspeq c_2$ needs $M\sspeq N\sspeq 1$, by primitivity, and $c_3\sspeq 4$, and $c_{4,+}\sspeq 4\,q\sspeq 2\cdot 6\sspeq 12$, from $c_{4,-}\sspeq 0$. If $c_1\sspkl c_2$ then $1\sspleq M\sspkl N$, and $c_3\sspeq (M\sspp N)^2$. $c_{4,+}\sspeq 4\,q\sspeq 2\,(M^2\sspp N^2 \sspp (M\sspp N)^2)\sspeq 4\,(M^2 \sspp N^2 +M\,N)$.\pn 
The $DS$-triples are primitive because also $\gcd((M \sspp N)^2,\,M^2)\sspeq 1\sspeq \gcd((M+N)^2,\,N^2)$ due to the same eqs. without the squares: $\gcd(M\sspp N,\,N)\sspeq \gcd(N,\,M\sspp N)\sspeq \gcd(N,\, M\sspp N\sspm N)\sspeq \gcd(N,\,M)\sspeq \gcd(M,\,N) \sspeq 1$. Similarly for the other $\gcd$. \hskip 7cm $\square$
\psn
For the characteristic triangular array $T(N,\,M)$ with $\gcd(N,\,M)\sspeq 1$ see \seqnum{A054521}. There a $0$ entry stands for $\gcd(N,\,M) \sspgr 1$.
\pbn
For the multiplicity $m_{3}(n)$ of primitive $DS$-triples with $c_{4,-}\sspeq 0$ and $c_3 \sspeq n^2$, for $n\sspgeq 2$, the following {\sl Proposition} holds.
\psn
{\bf Proposition 6: Multiplicity $\bf m_{3}(n)$ for $\bf c_3\sspeq n^2$, if $\bf c_{4,-}\sspeq 0$}\psn
{\it
\Beq
m_3(n) \sspeq \Caseszwei{1}{$ \rm{if}\ n\sspeq 2\, ,$}{\varphi$(n)/2$}{$ \rm{if}\ n\sspgeq 3\, .$}
\Eeq 
Here $\varphi$ is {\sl Euler}'s totient function \seqnum{A000010}, and $m_3(1)\sspeq 0$, of course. }
\psn
{\bf Proof:} The multiplicity of $c_3\sspeq (M \sspp N)^2\sspeq n^2$, for $ 1\sspleq M\sspleq N$, is the sane as the one for $N \sspp M\sspeq n$, for $n\sspgeq 2$, which is given by the number of partitions of $n$ with two relative prime parts, given in \seqnum{A023022}$\sspeq \{1,\, 1,\, 1,\, 2,\, 1,\, 3,\,2,\, 3,\,2,\, 5, ...\}$.\pn
This can be proved using the triangle $T(n,\,M)\sspeq$\seqnum{A054521}$(n,
\, M)$, for $n\sspgeq 2\,$ and $M\sspeq 1,\,2,\,...,\,n-1$ (the original triangle without its diagonal), because the sum over the anti-diagonals $T(n-M,\,M)$, for $M\sspeq 1,\,2,\,...,\,\floor{\frac{n}{2}}$, give $m_3(n)$, and due to $\gcd(n-M,\,M)\sspeq \gcd(n,M)$, this is \dstyle{m_3(n)\sspeq \sum_{m=1}^{\floor{\frac{n}{2}}}\,T(n,\,M)}. Then $m_3(2)\sspeq 1$, and for even $n\sspgeq 4$ the entries are symmetric around the middle $0$ (due to the $\gcd$ property), and for odd $n$ the first $(n-1)/2$ entries are repeated. Therefore, in both cases the sum is one half of the sum over row $n\sspgeq 3$ which is $\varphi(n)$. \hskip 12cm $\square$
\psn
{\bf Example 3:} $n\sspeq 5, c_3\sspeq 25$, with the two ($\varphi(5/2)\sspeq 2$) primitive 
$DS$-triples $[4,\,9,\,25]$ and  $[1,\,16,\,25]$, with $q\sspeq 19$ and $21$, respectively. 
\psn
The instances $c_3\sspeq n^2$, for $n\sspeq 2,\,3,\, ..., 24$, are shown in {\sl Table 6}.
\pbn
\section{DS-triples with distinct curvatures}
Now we have to analyze the primitive $DS$-triples with distinct (pairwise different) curvatures $[c_1,\, c_2,\, c_3]$ satisfying $0\sspkl c_1\sspkl c_2\sspkl c_3$ and $gcd(c_1,\,c_2,\,c_3)\sspeq 1$. See {\sl Table 2} for these triples with $c_3$ values from $6,\,7,\, ...,\, 38$ (except for $c_3\sspeq 8$). This table was found by testing all distinct primitive increasingly ordered positive integer triples for (positive) integer $q\sspeq \sqrt{c_1\,c_2\sspp c_1\,c_3\sspp c_2\,c_3}$.
\pn
The $c_j$'s are distinct, and $q$ can come in three versions: \psn
\Beq
{\bf i)}\ q\sspeq c_3,\ \ {\bf ii)} \sqrt{2\,c_1^2\sspp c_2^2}\sspkl q\sspkl c_3, \ \ \text{or}\ \ {\bf iii)}\  q\sspgr c_3.
\Eeq
In {\bf ii)} the first inequality is trivial. Examples for these three cases are found in {\sl Table 2}: $\bf i)$\, $c_3\sspeq q\sspeq 6$, $[2,\,3,\,6]$, $\bf ii)$ $c_3\sspeq 9 \sspgr 7\sspeq q\,(\sspgr \sspeq 3\,\sqrt{2}\,)$, $[1,\,4,\,9]$, and $\bf iii)$\, $c_3\sspeq 7\sspkl 9\sspeq q, [3,\,6,\,7]$. 
\pbn
{\bf Case i) $\bf q\sspeq c_3$}\psn
This part turns out to become quite lengthy. The result will be given finally in {\sl Proposition 7}.
\psn
Put $c_1\sspeq x$, $c_2\sspeq y$ and $c_3\sspeq q$, with $0\sspkl x\sspkl y$ and $\gcd(x,\,y,\,q)\sspeq 1$. Now $q^2\sspm (x+y)\,q \sspm x\,y\sspeq 0$ or, because integer $q\sspgr 0$,\ $2\,q\sspeq x\sspp y\sspp s$ with $s\sspdef \sqrt{(x\sspp y)^2\sspp 4\,x\,y}$, and $s\sspin \mathbb  N$. This results in an indefinite binary quadratic form representing $s^2$, \viz
\Beq
x^2 \sspp 6\,x\,y\sspp y^2 \sspeq s^2,\,\ \text{with}\ s\sspin \mathbb N,\,\ 0\sspkl x\sspkl y\,.  
\Eeq
This binary quadratic form is indefinite because the discriminant $D\sspeq 6^2\sspm 4\sspeq 32\sspgr 0$. Therefore, if for given $s\sspgeq 1$ an integer solution $(x,\, y)$ exists (not regarding the inequality) there will be infinitely many such solutions.
\pn 
Solutions satisfying also the inequality are found within the range $2\,\sqrt{2}\,x \sspkl s \sspkl 2\,\sqrt{2}\,y$. \Eg $(x,\,y)\sspeq (2,\,3)$ with $q\sspeq c_3 \sspeq 6$ needs $6\sspleq s\sspleq 8$, and there is fact only a solution for $s\sspeq 7$ with $(x,\,y)\sspeq (2,\,3)$. Or the $DS$-triple $[3,\,10,\,15]$ qualifies with $s\sspeq 17$ and $9\sspkl 17\sspkl 28 $. See {\sl Table 2}.\pn
Note that this quadratic form (not considering the inequality) is symmetric in $x$ and $y$, and under a common sign flip of $x$ and $y$.
\pn
We are interested in primitive $DS$-triples, and the following Lemma shows that this means $\gcd(x,\,y)\sspeq 1$, hence one is interested only in proper solutions of eq. $(34)$.
\pbn
{\it {\bf Lemma 2:}  $\gcd(x,\,y,\,q)\sspeq 1$ if and only if $\gcd(x,\,y)\sspeq 1$.}
\pbn
{\bf Proof:}  $\Rightarrow$)  Assume that $\gcd(x,\,y)\sspeq g\sspgr 1$, then $x\sspeq g\,\hat x$ and $y\sspeq g\,\hat y$, with $\gcd(\hat x,\,\hat y)\sspeq 1$, From \Eq{34} $s^2\sspeq g^2\,\hat s^2$, hence $s\sspeq g\,\hat s$. Therefore \dstyle{q\sspeq g\, \frac{\hat x\sspp \hat y\sspp \hat s }{2}}, and because $\hat x^2\sspp 6\,\hat x\,\hat y\sspp \hat y^2 \sspeq \hat s^2 $, one finds, independent of the parity of $\hat s$, that $\hat x\sspp \hat y\sspp \hat s$ is always even. Thus, $q\sspeq g\,\hat q$ with $\tilde q\sspin N$, and 
$\gcd(x,\,y,\,q)\sspeq g\sspgr 1$, in contradiction to $\gcd(x,\,y,\,q)\sspeq 1$.\psn
$\Leftarrow$) $\gcd(x,\,y,\,q)\sspeq \gcd(\gcd(x,\,y),\,q)\sspeq \gcd(1,\,q)\sspeq 1$. \hskip 8cm $\square$ 
\psn
It is convenient to diagonalize the binary quadratic form \Eq{34} in order to obtain a (generalized) {\sl Pell} equation. This is straightforward because $x^2 \sspp 6\,x\,y\sspp y^2 \sspeq s^2$ can be written as $(x \sspm y)^2 \sspm 2\,(x\sspp y)^2\sspeq -s^2$. Defining $X\sspdef y\sspm x\sspgr 0$ and $Y\sspdef x\sspp y \sspgr 2\,x \sspgr 0$ one has, 
\Beq
X^2 \sspm 2\, Y^2\sspeq -s^2\ ,\  {\text with} \ s\sspin \mathbb N, \ \text{and}\ \ 0 \sspkl  X\sspkl Y\,, 
\Eeq 
\ie $s\sspeq \sqrt{2\,Y^2\sspm X^2}$.\pn
Note that there are no integer solutions of this equation (without the inequality requirement) with vanishing $X$ or $Y$.
Every integer solution of \Eq{34} with $y\sspgr x\sspgr 0$ produces an integer  solution $(X,\,Y)$ with  $Y\sspgr X\sspgr 0$ of this $\sl Pell$ type equation. Conversely, given any integer solution  $(X,\,Y)$ of \Eq{35} with $Y\sspgr X\sspgr 0$  leads to \dstyle{x\sspeq \frac{Y\sspm X}{2}} and \dstyle{y\sspeq \frac{Y\sspp X}{2}} with $y\sspgr x\sspgr 0$. One has to check whether $x$ and $y$ are integers. $X\sspp Y$ (hence $Y\sspm X$) is indeed even due to the following {\sl Lemma}.
\pbn
{\it {\bf Lemma 3:} If $(X,\,Y)$ is an integer solution of $X^2 \sspm 2\, Y^2\sspeq -s^2$, for $s\sspin \mathbb N$, then $X\sspp Y$ is even.}
\pbn
{\bf Proof:} We show that in fact $X$ and $Y$ have the same parity as $s$, hence $X\sspp Y$ is even. From \Eq{36} follows that if $s$ is even then $X$ is even. Then also $Y$ is even, because with $s\sspeq 2\,\tilde s$ and $X\sspeq 2\,\tilde X$ one has $Y^2\sspeq 2\,(\tilde X^2 \sspp \tilde s^2)$. If $s$ is odd then $X$ is odd, and with $s\sspeq 2\,\tilde s\sspp 1$ and $X\sspeq 2\,\tilde X\sspp 1$ one has $Y^2\sspm 1\sspeq 4\,(T(\tilde X)\sspp T(\tilde s))$ with the triangular numbers \dstyle{T(n)\sspeq n\,(n+1)/2} (given in \seqnum{A000217}). Therefore, $Y^2\sspequiv 1 (mod\, 4)$, and $Y$ is also odd. \hskip 16.3cm $\square$
\psn
The next question is whether proper solutions $(x,\,y)$ of \Eq{34} correspond to proper solutions $(X,\,Y)$ of \Eq{35}. This is indeed the case as the next {\sl Lemma} shows.
\pbn
{\it {\bf Lemma 4:}  $\gcd(x,\,y)\sspeq 1$ if and only if $\gcd(X,\, Y)\sspeq 1$.}
\pbn
{\bf Proof:} $\Rightarrow$) Suppose that $\gcd(X,\, Y)\sspeq g\sspgr 1$. Hence $y\sspm x\sspeq X\sspeq g\,\tilde X$ and $x\sspp y\sspeq Y\sspeq g\,\tilde Y$ and $\gcd(\tilde X,\,\tilde Y)\sspeq 1$. Hence $2\,x\sspeq g\,\tilde x$ and $2\, y\sspeq g\,\tilde y$ with $\tilde x\sspeq \tilde Y\sspm \tilde X$ and $\tilde y\sspeq \tilde Y\sspp \tilde X$.\pn
Case $g$ odd, $\sspgeq 3$: $\tilde x$ and $\tilde y$ are both even and $x\sspeq g\,\frac{\tilde x}{2}$ and $y\sspeq g\,\frac{\tilde y}{2}$ and $\gcd(x,\,y)\sspgeq g\sspgr 1$ contradicting $\gcd(x,\, y)\sspeq 1$.\pn
Case $g$ even, $\sspgeq 2$: $g \sspeq 2\,\tilde g$ with $x\sspeq \tilde g\,\tilde x$  and $y\sspeq \tilde g\,\tilde y$. If $\tilde g \sspgeq 2$ then $\gcd(x,\,y)\sspgeq \tilde g\sspgeq 2$, contradicting $\gcd(x,\, y)\sspeq 1$. If $\tilde g \sspeq 1$, \ie $\gcd(X,\, Y)\sspeq 2$, hence $X\sspeq 2\,\tilde X$ and $Y\sspeq 2\,\tilde Y$, with $\gcd(\tilde X,\,\tilde Y)\sspeq 1$, satisfying $\tilde X^2 \sspm 2\,\tilde Y^2\sspeq - \tilde s^2$, with $s\sspeq 2\,\tilde s$. But $\tilde X\sspp \tilde Y$ has to be even by applying the argument from {\sl Lemma 3}. But both numbers can only be odd because their $\gcd$ equals $1$, and then \dstyle{\gcd(x,\,y)\sspeq \gcd(\tilde Y\sspm \tilde X,\,\tilde Y\sspp \tilde X) = 2\gcd\left(\frac{\tilde Y\sspm \tilde X}{2},\,\frac{\tilde Y\sspp \tilde X}{2}\right) \sspgeq 2 \sspgr 1}, (because $\tilde X$ and $\tilde Y$ are both odd), contradicting $\gcd(x,\, y)\sspeq 1$.  
\psn
$\Leftarrow$) Suppose that $\gcd(x,\, y)\sspeq g\sspgr 1$, then $x\sspeq g\,\tilde x$ 
and $y\sspeq g\,\tilde y$ with  $\gcd(\tilde x,\, \tilde y)\sspeq 1$. Thus, $X\sspeq g\, (\tilde y\sspm \tilde x)\sspgr 0 $ and $Y\sspeq g\,(\tilde y\sspp \tilde x)\sspgr 0$, and $\gcd(X,\,Y)\sspgeq g\sspgr 1$ contradicting $\gcd(X,\, Y)\sspeq 1$. \hskip 4cm $\square$
\psn
For proper solutions $(X,\,Y)$ of eq. $(35)$ the even $s$ case is excluded because then $X$ and $Y$ would both be even (as shown in the proof of {\sl Lemma 3}). Therefore we only have to consider odd $s$ with odd $X$ and odd $Y$. In this case eq. $(35)$ is equivalent to the following problem involving triangular numbers $T(n) \sspeq $\seqnum{A000217}$(n)$.
\Beq 
T(\widehat X) \sspp T(\widehat s)\sspeq 2\,T(\widehat Y)\, ,
\Eeq 
where $s\sspeq 2\, \hat s\sspp1$, $X\sspeq 2\, \widehat X\sspp 1 $, $Y\sspeq 2\, \widehat Y\sspp 1$, and  $ 0\sspleq \widehat X\sspkl \widehat Y$.\psn
{\bf Example 4:} A proper integer solution of eq. $(35)$ with $s\sspeq 7$ is  $(X,\,Y)\sspeq (1,\,5)$ (one of the two proper positive fundamental solutions of this equation satisfying also  $X \sspkl Y$; the other being $(17,\,13)$ with $Y\sspkl X$). This corresponds to the solution of eq. $(36)$ with $(\hat X,\,\hat Y)\sspeq (0,\,2)$ and $\hat s\sspeq 3$. This shows that the relevant solutions $(\hat X,\,\hat Y)$ of eq $(36)$ do not have to be proper even if $\gcd(X,\,Y)\sspeq 1$. (The other solution leads to $T(8)\sspp T(3)\sspeq 2\, T(6)$, \ie  $36\sspp 6\sspeq 2\cdot 21$. Again, this is not a proper solution of eq. $(36)$, and it does not satisfy the inequality). \pn
The solution $(X,\,Y)\sspeq (1,\,5)$ corresponds to one of the proper solutions of eq. $(35)$, \viz\, $(x, y)\sspeq (2,\,3)$ with $0 \sspkl x\sspkl y$. The other solution $(17,\,13)$ leads to $(x,\, y)\sspeq (-2,\,15)$ which belongs to the family with the positive fundamental solution $(3,\, 2)$, which does not satisfy the inequality.
\pbn
After this digression we give the {\sl Proposition} with the solution of \Eq{35} with the restrictions on $X$ and $Y$.
\pbn
\vfill\eject\noindent
{\bf Proposition 7: Solution of $\bf X^2 \sspm 2\,Y^2\sspeq -s^2$, with $\bf 0\sspkl X\sspkl Y$}
\psn
{\it Proper solutions ($\gcd(X,\,Y)\sspeq 1$) of the generalized {\sl Pell} equation \Eq{35} with restriction are possible if and only if $s\sspeq s(n)\sspeq$\seqnum{A058529}$(n)$, with $n \sspgeq 2$, \ie $s$, $X$ and $Y$ are odd, and there are exactly $P(n)\sspeq 2^{P_1(n)\sspp P_7(n)-1}$ solutions, where $P_1(n)$ and $P_7(n)$ are the numbers of distinct prime divisors of $s(n)$ congruent to $1$ and $7$ modulo $8$, respectively.\psn
The (here irrelevant) improper solutions arise only from composite $s$ with divisors that are primes congruent to $1$ or $7$ modulo $8$.\pn
The case $s\sspeq 1$ does not qualify because there is only one family of solutions (the ambiguous case) with positive proper fundamental solution $(X,\, Y)\sspeq (1,\,1)$, and positive solutions have $X\sspgr Y\sspgr 0$.  
}
\psn
Before giving the proof the following example of such a generalized {\sl Pell} equation, but not representing a negative square, may be useful, and then {\it Lemmata 5} to {\it 13} will follow.
\psn
{\bf Example 5: $X^2 \sspm 2\, Y^2 \sspeq -7\cdot 17\cdot  23 \sspeq -2737$} 
\psn
With the restriction $0\sspkl Y$ (always obtainable by a sign flip in $X$ and $Y$), but both signs for $X$, one has $8$ families of proper solutions (and no improper ones). The {\it positive proper fundamental solutions} ({\it ppfs}) $(X,\, Y)$ are $(1,\,37)$ and $(145,\,109)$,\ $(25,\,41)$ and $(89,\,73)$,\ $(31,\,43)$ and $(79,\,67)$,\ $(41,\,47)$ and $(65,\,59)$. Each of these pairs generates an infinite conjugate family, consisting of pairs related by sign flip in the $X$ entries. Only the first members of the four pairs are found as the unique solutions satisfying also $X\sspkl Y$. There are precisely $4 \sspeq 2^{1+2-1}$ solutions with $P_1\sspeq 1$, from $17$ and $P_7\sspeq 2$, from $7$ and $23$.
\pbn
For the proof of {\sl Proposition 7} the method employing {\it representative parallel primitive forms} ({\it rpapf}s) for solving generalized {\sl Pell} equations, adapted to the present case, are first recalled. This method has been described in \cite{WLI} and \cite{WLII}, where also references are given.
\psn
A necessary ingredient is the first {\it reduced form} obtained from a given unreduced one. For the definition of reduced forms see \cite{WLI}, also for references.
\psn
{\bf Lemma 5: Half-reduced right neighbor form}
\psn 
{\it  The half-reduced right neighbor of a form $F\sspeq [a,\,b,\,c]$ of discriminant $Disc(F)\sspeq b^2\sspm 4\,a\,c$ is \pn $\widetilde F\sspeq [c,\, -b\sspp 2\,c\,t,\,c']$, with $c'$ such that $Disc(\widetilde F)\sspeq Disc(F)$, and a unique \dstyle{t \sspeq \ceil{\frac{f(Disc(F))+b}{2\,c}\sspm 1}}  if $c\sspgr 0 $ and \dstyle{ t \sspeq \floor{1\sspm \frac{f(Disc(F))+b}{2\,|c|}}}  if $c\sspkl 0 $, where $f(Disc(F))\sspdef \ceil{\sqrt{Disc(F)}}$. Vanishing $c$ is excluded because it leads to a factorization of $F$. For the present case, \Eq{35}, $Disc(F) \sspeq 8$ and $f(8)\sspeq 3$.}
\psn
{\bf Proof:} For the definition of a {half-reduced neighbor form} of $F$ see \cite{WLI}, p. 2 (where $\widetilde F$ is denoted by $F^R$, and $Disc$ by $D$). Note that the transformation from $F$ to $\widetilde F$ is accomplished by the so called ${\bf R}(t)$-transformation, with a matrix occurring in the next {\it Lemma}. Therefore half-reduced right neighbor transformations will simply be called ${\bf R}$-transformations in the following. \hskip 6cm $\square$
\psn
{\bf Lemma 6: Pell form, its reduced principal form, cycle, and automorphic matrix}
\psn
{\it {\bf a)} The {\sl Pell} form $F_{Pell}({\bf A},\vec X)\sspeq \vec X^{\top} {\bf A} \vec X$ with $\vec X^{\top}\sspeq (X,Y)$ ($\top$ for transposed) and ${\bf A}\sspeq Matrix([1,0],\,[0,-2])$, \ie $F_{Pell}\sspeq X^2\sspm 2\,Y^2$, also  denoted by $[1,\,0,\,-2]$, has first reduced form, called principal (reduced) form $F_p\sspeq [1,\,2,\,-1]$. This is obtained after two proper (Determinant $+1$) equivalence transformation, \ie ${\bf A}_p\sspeq {\bf B}^{\top}\,{\bf A}\,{\bf B}$, where ${\bf B}\sspeq {\bf R}(0)\,{\bf R}(1)$, where ${\bf R}(t)\sspeq Matrix([0,-1],\,[1,t])$. Thus, ${\bf B}\sspeq -Matrix([1,1],\,[0,1])$. The discriminant of these two forms is $Disc \sspeq -4\,Det({\bf A})\sspeq 8$, implying that these are indefinite binary quadratic forms. 
\pn
{\bf b)} The cycle generated by the principal form $F_p\sspeq [1,\,2,\,-1]$ is a $2$-cycle with (reduced) partner $F^{\prime}\sspeq [-1,\,2,\,1]$, obtained from $F_p$ by an ${\bf R}(-2)$ transformation. A further ${\bf R}(+2)$ transformation leads from $F^{\prime}$ back to $F_p$. 
\psn
{\bf c)} The automorphic matrix $\bf Auto$ related to this $2-$cycle is (first $t\sspeq -2$ then $t\sspeq +2$)
\Beqarray
&&{\bf Auto}\sspeq {\bf R}(-2)\,{\bf R}(2) \sspeq -Matrix([1,2],\,[2,5])\,, \nonumber \\
&&{\bf Auto}^{-1}\sspeq -Matrix([5,-2],\,[-2,1])\,.
\Eeqarray
\psn
{\bf d)} Later becomes important:  
\Beqarray
&&{\bf Auto}^{\prime}\sspdef {\bf B}\,(-{\bf Auto})\,{\bf B}^{-1}\sspeq Matrix([3,4],\,[2,3])\,, \nonumber \\
&&({\bf Auto}^{\prime})^{-1}\sspeq Matrix([3,-4],\,[-2,3])\,.
\Eeqarray
}
\pbn
{\bf Proof:} 
\pn
{\bf a)} Two half-reduced right neighbor transformations are needed to obtain the first reduced form $F_p$ from $F_{Pell}$, The $t$ values are $0$ and $1$. This leads to the transformation matrix {\bf B} with the two ${\bf R}-$transformations.   
 \psn
 {\bf b)} ${\bf A}_p\sspeq Matrix([1,1],\,[1,-1]$. The matrix for $F^{\prime}$ is ${\bf A}^{\prime}\sspeq {\bf R}(-2)^{\top}\,{\bf A}_p\, {\bf R}(-2)\sspeq Matrix([-1,1],\,[1,1])$.
 Then ${\bf A}_p\sspeq {\bf R}(2)^{\top}\,{\bf A}^{\prime}\,{\bf R}(2) $, This leads to the automorphic matrix $\bf Auto$ given in part {\bf c)}. 
 \psn
The other parts {\bf c)} and {\bf d)} are  obvious. \hskip 10.5cm $\square$ 
\pn
For some application of automorphic matrices for ordinary {\sl Pell} forms with $D(n) \sspeq$\seqnum{A000037}$(n)$ representing $k\sspeq +1$ see \cite{WLII}, section $2$. 
 \pbn
 {\bf Lemma 7:\ Representative parallel primitive forms, and solutions of generalized Pell equations}
\psn
{\it {\bf a)} For the solution of the generalized Pell equation $X^2\sspm D\,Y^2\sspeq k$, with non-vanishing integer $k$, $D\sspin \mathbb N$, not a square, so-called representative parallel primitive forms ({\it rpapf}s) are used. They are determined by $Disc$ and $k$. For even discriminant Disc\sspeq 4\,D, (later $D\sspeq 2$), the \rpapfs $F_{pa}(D,\,k;\,j)$ are obtained from solving the congruence $j^2\sspm D\sspequiv \modan{0}{k}$, for $j\sspin \{\range{0}{\abs{k}-1}\}$, and restricting to primitive forms $F_{pa}(D,\,k;\,j)\sspeq [k,\,2\,j,\, c(D,\,k;\,j)]$, with \dstyle{c(D,\,k;\,j)\sspeq \frac{j^2\sspm D}{k}}. The case of odd $Disc$ is here not relevant.
\psn
{\bf b)}\, All families of infinitely many proper solutions of $X^2 \sspm D\,Y^2\sspeq k$ are found from the {\it rpapf}s that reach the principal form $F_p(D)$ by half-reduced right neighbor transformations. For each such {\it rpapf} the infinitude of proper solutions, indexed by $i$, is
\Beq
\LTvez{X_i(j(D,\,k))} {Y_i(j(D,\,k))}\sspeq {\bf B}(D)\,({\bf Auto}(D))^i\,{\bf Rpa}(\vec t(j(D,\,k)))^{-1}\,\LTvez{1}{0}, \ \ {\text for}\ i\sspin \mathbb Z,     
\Eeq
with $j(D,\,k)$ a quadratic residue $j$ modulo $k$ from part {\bf a)} that gives a primitive form, and the definition ${\bf Rpa}(\vec t(j(D,\,k)))\sspdef {\bf R}(t_1)\cdots {\bf R}(t_L)$ , with the $L$-tuple $\vec t\sspeq (t_1,\,...,\,t_L)$ determined from reaching for the first time the reduced principal form $F_p(D)$ from the \rpapf $F_{pa}(j(D,\,k))$ via half-reduced right neighbor ${\bf R}(t)$-transformations. The number of such families is denoted by $f\sspeq f(D,\,k)$ (not to be confused with $f\sspeq \ceil{\sqrt{Disc}}$ in the definition of reduced forms in \cite{WLI}).
\pn
Note that if \Eq{39} produces a solutions with negative $Y_i$ then one applies a common sign change in $X_i$ and $Y_i$. Only solutions with non-negative $Y$ will then be considered. Therefore one could as well use $-{\bf Auto}$ or $-{\bf B}$ in \Eq{39}. 
}
\psn
This method using {\it rpapf}s for $D$ and $k$ to find all proper solutions of a generalized {\sl Pell} equation of discriminant $4\,D$, representing a non-vanishing $k$, produces automatically some fundamental solution for $i\sspeq 0$. But $X$ may be negative. Here we define fundamental solution as positive ones with smallest $X$ (but sometimes the $i\sspeq 0$ solutions with negative $X$ are also called fundamental, like in the later {\sl Example 7}). This method differs from {\sl Nagell} \cite{Nagell} where inequalities are used in {\it Theorem 108} and in {\it Theorem 108a}, pp. 205-207, for positive and negative $k$ (there called $N$), respectively, to find fundamental solutions with $|X| \sspgr 0$, $Y\sspgeq 0$ and $|X|\sspgeq 0$, $Y\sspgr 0$, respectively. Here one has to find first all {\it rpapf}s and discard those not related to the fundamental form $F_p(D)$. \psn
\psn
{\bf Example 6:} For $D\sspeq 2$ and $k\sspeq -7\cdot 17\cdot 23 \sspeq -2737$ the \rpapfs come from the $f \sspeq 8$ solutions of $j(2,2737)$, \ie $74,\, 465,\, 1145, \,  1201,\, 1536,\, 1592,\, 2272,\, 2663$, with corresponding $c(2,2737;j)$ values $2,\, 79,\, 479,\, 527,\, 862,\, 926,\, 1886,\, 2591$. All $8$ \rpapfs reach the $F_p(2)\sspeq [1,\,2,\,-1]$. See the later proof of {\sl Lemma 8 b)}.
\pbn
{\bf Proof:} \psn
{\bf a)} The notion of representative parallel primitive forms (\rpapfs), here for even discriminants, is found in \cite{WLII}, {\it section 3}, with references to {\sl Scholz-Sch\"oneberg} \cite{SchSch} and {\sl Buell} \cite{Buell}. 
\psn
{\bf b)} Observe first that the primitivity of forms is respected by each ${\bf R}(t)$ transformation. With $\gcd(a,\, b,\, c)\sspeq 1$ and $a^{\prime} \sspeq c$, $b^{\prime}\sspeq -b\sspp 2\,c\,t$ and $c^{\prime}\sspeq a\sspm b\,t\sspp c\,t^2$ (from the Disc invariance, with $c\sspneq 0$) assume that $x^{\prime} \sspeq g\,\tilde x^{\prime}$, for $x\sspin \{a,\,b,\,c\}$, where $g \sspgeq 2$ and  $\gcd(\tilde a^{\prime},\,\tilde b^{\prime},\, \tilde c^{\prime}) \sspeq 1$. Then the definition of  $a^{\prime}$ implies $c\sspeq g\,\tilde a^{\prime}$. From $b^{\prime}$ follows then that $b\sspeq g\,\tilde b$ with  $\tilde b \sspeq 2\,\tilde a^{\prime}\,t \sspm \tilde b^{\prime}$.  From $c^{\prime}$ follows that $a\sspeq g\,\tilde a$ with $\tilde a\sspeq \tilde c^{\prime}\sspp (\tilde b \sspm \tilde a^{\prime}\,t \sspm)\,t$. Thus, the contradiction $1\sspeq gcd(a,\,b,\,c)\sspgeq g \sspgeq 2$ follows.
\psn
Also a proper solution $\vec X\sspeq (X,\,Y)^{\top}$, \ie $\gcd(X,\,Y)\sspeq 1$, of a form $F$ representing some $k$ is transformed to $\vec X^{\prime} \sspeq {\bf R}^{-1}(t)\,\vec X\sspeq (t\,X\sspp Y, -X)^{\top}$ (see the next paragraph). And this implies that $\gcd(t\,X\sspp Y, -X)\sspeq \gcd(-X,\, (t\,X\sspp Y) -t\,X)\sspeq \gcd(X,\,Y)\sspeq 1$. Thus, a proper solution transforms into a proper solution.
\psn 
A parallel form $[k,\,2\,j,\,c(D,\,k;\,j)]$ has the special trivial non-negative proper fundamental solution $(X,\, Y)\sspeq(1,\,0)$. If forms are related by ${\bf A}^{\prime}\sspeq {\bf R}^{\top}\, {\bf A}\, {\bf R}$, where $F\sspeq F({\bf A},\,\vec X)$ and $F^{\prime}\sspeq F^{\prime}({\bf A}^{\prime},\,\vec {X^{\prime}})$, then the solutions of each representation of $k$ of these two forms are related by $\vec {X^{\prime}}\sspeq  {\bf R}^{-1}\, \vec X$. This explains the product of inverse ${\bf R-}$transformations needed to reach a fundamental solution of the principal form $F_p$ representing $k$. Then one has to reach a fundamental solution of $F_{Pell}$ representing $k$ by application of ${\bf B}$, because ${\bf A}_{Pell}\sspeq {\bf B}^{-1,\top}\ {\bf A}_p\,{\bf B}^{-1}$, with $F_p\sspeq F_p({\bf A}_p,\,\vec X_p)$. Looping around the reduced cycle determined by $F_p$ in both directions leads to all solutions of the family corresponding to the chosen \rpapf\   reaching  $F_p$. (Note that what is called here family is called class by {\sl Nagell} \cite{Nagell}, p. 204. The notion class is here reserved for the number of reduced cycles for given discriminant, called $h(Disc)$). \pn
That every unreduced primitive form reaches a reduced form after repeated application of $\bf R$-transformations is the content of of \cite{SchSch}, Satz 70 on p. 113. 
\hskip 7.1 cm $\square$  
\pbn
{\bf Lemma 8: Conjugated rpapfs for $\bf X^2 \sspm 2\,Y^2\sspeq -s^2$.}
\psn
{\it 
{\bf a)} Each representative parallel primitive form (rpapf) $F_{pa}(D\sspeq 2,\, k=-s^2,\, j)$, $s\sspgr 1$, of $X^2 \sspm 2\,Y^2\sspeq -s^2$, for $j\sspin \{j_1,\,...,\,j_{f(2,-s^2)}\}$, is not reduced. 
\psn
{\bf b)} Each $F_{pa}(2,\,-s^2,\, j)$ reaches the cycle of reduced forms containing the principal form $F_p\sspeq [1,\,2,\,-1]$ after applying $\bf R$- transformations. Hence $f(2,-s^2)$, the number of families of proper solutions with $Y \sspgr 0$, equals the number of rpapfs.
\psn
{\bf c)} If $s$ is an odd positive integer number $\sspneq 1$ then the rpapfs come in conjugate pairs \dstyle{\left[-s^2,\, 2\,j,\, -\frac{j^2\sspm 2}{s^2}\right]} and \dstyle{\left[-s^2,\, 2\,(s^2\sspm j),\,-\frac{(s^2\sspm j)^2\sspm 2}{s^2}\right]} for each solution $j$ of the congruence $j^2 \sspm 2\sspequiv \modan{0}{s^2}$, with $j\sspin\{0,\,1,\,...,\,s^2\sspm 1\}$. Hence $f(2,\,-s^2)$ is a positive even integer number.
}\pbn
{\bf Proof:}
\psn
{\bf a)} For the definition of reduced indefinite binary forms $F\sspeq [a,\,b,\,c]$ of discriminant $Disc$ see the two equivalent definitions in \cite{WLI} eqs. $(1)$ and $(2)$ (where the discriminant is named $D$). 
This implies (see eq. $(3)$ there) that $4|a\,c| \sspkl Disc = 4\,D $, \ie $|a\,c|\sspkl 2$. \pn 
But from {\sl Lemma 7} each rpapf $F_{pa}$ satisfies $a\,c \sspeq j^2\sspm 2$ (where $a\sspeq k$). Hence $|a\,c| \sspeq |j^2\sspm 2|$, for $j\sspin \{0,\,1,\,...,\,|k|\sspm 1\}$ (here $|k|\sspm 1\sspeq s^2 - 1\sspgr 3$). Thus, $|j^2\sspm 2|\sspgr 2$ for $j\sspgeq 2$, and only $j\sspeq 1$ satisfies $|1 \sspm 2|\sspkl 2$, but $1\sspm 2\sspnotequiv 0 (mod -s^2)$, for $s \sspgr 1$.
\psn
{\bf b)} Because the number of reduced cycles of discriminant $Disc\sspeq 4\,D \sspeq 8$, called the class number $h(Disc)$, that is here $h(8)\sspeq 1$ (\eg \cite{Buell}, Table\, 2B,\, p. 241), there is only the reduced 2-cycle from {\sl Lemma 6 b)}, namely the principal form $F_p\sspeq [1,\,2, -1]$, and its ${\bf R}(-2)$-transform $F^{\prime}\sspeq [-1,\, 2,\, 1]$. These two reduced forms can also be found directly by using the definition II of reduced forms given in \cite{WLI}, eq. $(3)$ for $Disc\sspeq 8$ (there called $D$).\pn 
As mentioned above in the proof of {\sl Lemma 7} b) every indefinite primitive form of discriminant $Disc$ which is not reduced reaches a reduced form after a chain of half-reduced neighbors, by applying ${\bf R}(t)$-transformations. This is the content of \cite{SchSch}, Satz 70 on p. 113. 
\psn
{\bf c)} Remember from the proof of {\sl Lemma 3} that for proper solutions $s$ has to be odd. That there are two different solutions for each $j$ follows from the degree $2$ of the solvable congruence. These two given solutions are different because $s^2\sspeq 2\,j$ is impossible because $s$ is odd. If $j^2\sspm 2 \sspequiv 0\,(mod\, s^2)$, for $s\sspgr 1$ and $j\sspin\{0, 1,\,...,s^2-1\}$ then also $(s-j)^2 \sspm 2 \sspequiv 0\,(mod \,s^2)$. Most $j$s do not qualify, and not all odd $s\sspgr 1$ are possible, as will be discussed later. These two forms are indeed primitive because $j\sspeq 0$ is not a solution of the congruence for $s\sspgr 1$, and with odd $s$ one has $\gcd(-s,\, 2\,j) \sspeq \gcd(s,\, j)\sspeq 1$ because the congruence shows that $j$ cannot have any odd prime divisor of odd $s \sspgeq 3$, otherwise this divisor would also divide $2$. Finally $\gcd(a,\,b,\,c) = \gcd(\gcd(a,\,b),\,c)$ is used.\ \hskip 9.3cm $\square$
\pbn
{\bf Example 7: Two conjugate families for} $\bf k = - 119^2\sspeq -14161,\,s \sspeq 119\sspeq 7\cdot 17$.
\psn
The two pairs of {\it rpapfs} are $[-14161,\,4136,\,-302]$, $[-14161,\, 24186,\, -10327]$ and $[-14161,\, 6448,\,-734]$, $[-14161,\,21874,\,-8447]$. With $-{\bf B}(D) = -{\bf B}(2)$ the corresponding fundamental solutions for $i\sspeq 0$ of \Eq{39} are  $[41,\,89]$, $[-41,\,89]$ and $[79,\,101]$, $[-79,\,101]$. The corresponding $\vec t$ values are $(-6,\, 2,\, 2,\, 2,\,2,\, 2,\, 3,\, 1)$, $(-1, 6, 7, 2, 2)$ and $(-4,\, 3,\, 3,\, 2,\, 2,\, 2, \,1) $, $(-1,\, 4,\, 2,\, 3,\, 5,\, 2)$.
\pbn
{\bf Lemma 9: Recurrence for solutions of generalized Pell equations}
\psn
{\it Eq.\,$(39)$, using -{\bf Auto}(D) instead of {\bf Auto}(D), implies a recurrence for non-negative powers and non-positive powers $i$ with input $\LTvez{X_0(j(D,\,k))} {Y_0(j(D,\,k))}$, with $Y_0(j(D,\,k))\sspgeq 0$, which involves ${\bf Auto}^{\prime}(D)$ like in \Eq{38}, but with arguments $D$ (not the instance $D\sspeq 2$).
\Beqarray
\LTvez{X_i(j(D,k))} {Y_i(j(D,k))} &\sspeq& {\bf Auto}^{\prime}(D)\,\LTvez{X_{i-1}(j(D,\,k))} {Y_{i-1}(j(D,\,k))}, \ \ \text{for} \ \ i\sspin \mathbb N,\\
\LTvez{X_{-|i|}(j(D,\,k))} {Y_{-|i|}(j(D,\,k))} &\sspeq& 
({\bf Auto}^{\prime}(D))^{-1}\,\LTvez{X_{-(|i|-1)}(j(D,k))} {Y_{-(|i|-1)}(j(D,k))}, \ \ \text{for} \ \ |i|\sspin \mathbb N, 
\Eeqarray
for every solution $j\sspeq j(D,\,k)$ satisfying the congruence given in  {\it Lemma 7}, a).
}
\pbn
{\bf Proof:} This is obvious. Again, if the result does not produce a non-negative $Y_i$, for $i\sspin \mathbb Z$, then a common sign flip in $(X_i,\,Y_i)$ will be applied. \hskip 11.5cm $\square$
\pbn
{\bf Lemma 10: Powers of $(-{\bf Auto})$ and ${\bf Auto}^{\prime}$ for $\bf D\sspeq 2$}
\psn
{\it Integer powers of the special $2\ssptimes 2$ matrices $-{\bf Auto}$, from {\it Lemma 6, c)}, \Eq{37}, and  ${\bf Auto}^{\prime}$, from {\it Lemma 6, d)}, \Eq{38}, are given by
\Beqarray
({\bf -Auto})^n&\sspeq& {\LTmatz{S(n-1,\,6)\sspm S(n-2,\,6)}{2\,S(n-1,\,6)}{2\,S(n-1,\,6)}{S(n,\,6)\sspm S(n-1,\,6)}},\ \ \text{for}\ \ n\sspin \mathbb Z\, , \\
({\bf Auto}^{\prime})^n&\sspeq& {\LTmatz{S(n,\,6)\sspm 3\,S(n-1,\,6)}{4\,S(n,\,6)}{2\,S(n,\,6)}{S(n,\,6)\sspm 3\,S(n-1,\,6)}},\ \ \text{for}\ \ n\sspin \mathbb Z\, ,
\Eeqarray
where the {\sl Chebyshev} polynomial system 
$\{S(n,\,x\sspeq 6)\}_{n=-\infty}^{+\infty}$, with recurrence $S(n,\,6)\sspeq 6\,S(n-1,\,6) \sspm S(n-1,\,6)$ entered, with inputs $S(-1,\, 6)\sspeq 0$ and $S(0,\,6)\sspeq 1$.
}
\psn
{\bf Proof:} This is a standard application of the {\sl Cayley-Hamilton} theorem (see \eg \cite{CourantHilbert}). The relevant recurrence is the one of the {\sl Chebyshev} $S-$polynomial system \{S(n,\, 6)\} given by \seqnum{A001109}$(n+1)$, for integer numbers $n \sspgeq -1$. For other integer indices this is given \via the recurrence by  $S(-|n|, 6) \sspeq -S(|n|-2,\, 6)$, for $|n|\sspgeq 2$. For $S(n,\,6)\sspm S(n-1,\,6)$ see \seqnum{A001653}$(n+1)$, for $n\sspin \mathbb N_0$, and for $S(n,\,6)\sspm 3\,S(n-1,\,6)$ see \seqnum{A001541}$(n)$, for $n\sspin \mathbb N_0$. \hskip 13cm $\square$  
\pbn
{\bf Corollary 2:}
\Beq
({\bf -B})\,({\bf -Auto})^n \sspeq {\LTmatz{3\,S(n-1,\,6)\sspm S(n-2,\,6)}{ S(n-1,\,6)\sspp S(n,\,6)}{2\,S(n-1,\,6)}{S(n,\,6)\sspm S(n-1,\,6)}},\ \ \text{for}\ \ n\sspin \mathbb Z\, .
\Eeq
\psn
Note that $3\,S(n-1,\,6)\sspm S(n-2,\,6)\sspeq$\seqnum{A001541}$(n)$, for $n\sspgeq 0$.
\pbn
{\bf Lemma 11: $\bf |X_i|/Y_i$ bound, and proper positive fundamental solution}
\psn
{\it 
{\bf a)} If there is a solution $(X_i,\, Y_i)$ of the generalized {\sl Pell} equation $X^2\sspm 2\,Y^2 \sspeq -s^2$, for $s\sspin \mathbb N$, with $Y_i\sspgr 0$, then
\Beq
\frac{|X_i|}{Y_i} \sspeq \sqrt{2}\,\sqrt{1\sspm \frac{1}{2}\,\left(\frac{s}{Y_i}\right)^2} \sspeq \sqrt{2}\bigg / \sqrt{1\sspp \left(\frac{s}{X_i}\right)^2} \sspkl \sqrt{2}\,.
\Eeq
{\bf b)} The positive proper fundamental solution ({\it ppfs}) $(\widehat X,\,\widehat Y)$ for one of the members of each pair of conjugate families of solutions of this {\sl Pell} equation representing $s \sspgr 1$ satisfies 
\Beq
0\sspkl \widehat X\sspkl \widehat Y\sspkl s\,.
\Eeq }
\psn
{\bf Proof:} \pn
{\bf a)} This follows immediately from this {\sl Pell} equation, and the fact that neither a solution with $Y_i\sspeq 0$ nor with $X_i\sspeq 0$ is possible.\pn
{\bf b)} Here {\sl Nagell}'s \cite{Nagell} inequalities from {\it Theorem 108a, pp. 206-7} are useful. The positive proper fundamental solution of $X^2\sspm 2\,Y^2\sspeq +1$, with only one family of solutions with $Y\sspgr 0$, called an ambiguous family, is $(\widetilde X,\, \widetilde Y)\sspeq (3,\,2)$. Hence, for $s\sspgr 1$, with $\widehat Y \sspgr 0$  and the member of a conjugate pair of (proper or improper) solutions (related by $X\mapsto -X$, $Y\mapsto Y$), with $\widehat X\sspgr 0$ ($\widehat X\sspeq 0$ is impossible), one has
\Beqarray
0\sspkl \widehat Y &\sspleq& \frac{2\,s}{\sqrt{2\,(3\sspm 1))}}\sspeq s,\nonumber \\
0\sspkl \widehat X &\sspleq& {\sqrt{\frac{1}{2}\,(3\sspm 1)}}\,s\sspeq s\,.
\Eeqarray
Therefore, if there is a proper solution to $X^2\sspm 2\,Y^2\sspeq -s^2$, where $s\spin \mathbb N\setminus \{1\}$, then from part {\bf a)}, 
\Beq
0\sspkl \frac{\widehat X}{\widehat Y}\sspeq \sqrt{2}\bigg/ \sqrt{1\sspp \left(\frac{s}{\widehat X} \right)^2}\sspkl \sqrt{2}\big /\sqrt{2} \sspeq 1,
\Eeq
because if $\widehat X\sspeq s$ then also $\widehat Y\sspeq s$, and the solution is not proper, hence ${\widehat X}\sspkl s$.\pn
For $s\sspeq 1$ the proper positive fundamental solution of the ambiguous family is $(1,\,1)$. Therefore this case has been excluded. \hskip 13.8cm $\square$
\pbn
\vfill\eject\noindent
{\bf Definition 3: Conjugated family pairs $\bf I$ and $\bf II$}
\psn
{\it For each conjugate family pair of proper solutions of $X^2 \sspm 2\,Y^2\sspeq -s^2$, for solvable odd $s\sspgeq 3$ (later shown in {\sl Proposition 7} to be $s\sspeq s(n)\sspeq$\seqnum{A058529}$(n)$, for $n\sspgeq 2$) family $I$ is defined as the one with the positive proper fundamental solution ({\it ppfs}) $X_{i^{\prime}}\sspeq \widehat X$ and $Y_{i^{\prime}}\sspeq \widehat Y$ of {\sl Lemma 11 b)} for a certain index $i\sspeq i^{\prime}$ in \Eq{39}, with $D\sspeq 2,\, k\sspeq -s^2$, after a certain sign choice for ${\bf B}(2)\sspeq {\bf B}$. 
\pn
The conjugate family $II$ has then, for a certain index  $i\sspeq i^{\prime\prime}$, $X_{i^{\prime\prime}}\sspeq -X_{i^{\prime}} \sspeq -\widehat X\sspkl 0$ and $ Y_{i^{\prime\prime}}\sspeq Y_{i^{\prime}}\sspeq \widehat Y\sspgr 0$.
}
\pbn 
{\bf Example 8: Proper fundamental solutions for families $\bf I$ and $\bf II$} 
\psn
For $D\sspeq 2$, $s\sspeq 17$, family  $I$ has the {\it ppfs} $(\widehat X, \,\widehat Y)\sspeq  (7, 13)$ with $i\sspeq i^{\prime} = 0$ in \Eq{39} with ${\bf B}(2)\sspeq - {\bf B}\sspeq Matrix([1,1],\,[0,1])$ and ${\bf Auto}(2)\speq -{\bf Auto}\sspeq Matrix([1,2],\,[2,5])$. For family $II$ $(-\widehat X, \,\widehat Y)\sspeq  (-7, 13)$ for $i\sspeq i^{\prime\prime}\sspeq 0$. For $i\sspeq 1$ this becomes the fundamental positive solution $(31,\, 25)$ (later not qualifying because $X\sspgr Y$). The corresponding  conjugate \rpapfs are $[-289,\, 488,\, -206]$ and $[-289,\, 90,\, -7]$ with $\vec t$ tuples $(-1,\, 6,\, 2,\, 1)$ and $ (-6,\, 2,\, 2)$, respectively.
\pbn
{\bf Lemma 12: Relations between solutions of conjugate families} 
\psn
i) The solution $(X_{i^{\prime\prime}},\,Y_{i^{\prime\prime}})^{II}$ \sspeq $(-X_{i^\prime},\,Y_{i^\prime})^{I}\sspeq (-\widehat X,\,\widehat Y)$,\psn
ii) for $ (X_{i^{\prime\prime}\pm k},\,Y_{i^{\prime\prime}\pm k})^{II} \sspeq (-X_{i^\prime\mp k},\,Y_{i^\prime\mp k})^{I}$, for $k\sspin \mathbb N$, with correlated signs. 
\pbn
{\bf Proof}\pn
{\bf i)} From {\sl Definition 3}.\psn
{\bf i)} Induction over $k$ is used, with the help of the recurrences from {\sl Lemma 9}, with \Eq{40} or \Eq{41} (adapted with $D\sspto 2,\, k\sspto -s^2$) for positive or negative signs, respectively. The base case $k\sspeq 0$ is given in part {\bf i)}. Assuming that for positive signs the claim holds for $k\sspeq 1,\,2,\,...,\,j$  then (we omit the superscripts) $X_{i^{\prime\prime}+j+1}\sspeq 3\, X_{i^{\prime\prime}+j} \sspp 4\,Y_{i^{\prime\prime}+j}\sspeq -(3\,X_{i^\prime - j} \sspm 4\, Y_{i^\prime - j})$  and $Y_{i^{\prime\prime}+j+1}\sspeq 2\, X_{i^{\prime\prime}+j} \sspp 3\,Y_{i^{\prime\prime}+j}\sspeq -2\,X_{i^\prime-j} \sspp 3\,Y_{i^\prime-j}$ Using the recurrence \Eq{41} in family $I$ this implies $X_{i^{\prime\prime}+j+1}\sspeq -X_{i^\prime-(j+1)}$. Similarly for $Y_{i^{\prime\prime}+j+1}\sspeq Y_{i^\prime-(j+1)}$.\pn
The proof for negative signs runs along the same line. \hskip 8cm $\square$
\pbn
{\bf Lemma 13: Relations between $\bf X_i$ and $\bf Y_i$, signs of $\bf X_i$, and uniqueness}\psn
{\it Consider only proper solutions of $X^2\sspm 2\,Y^2\sspeq -s^2$, with odd $s\sspgeq 3$, and a fixed family of conjugate pairs.
\pn 
{\bf a)} Family $I$\pn
{\bf i)} There is only one solution satisfying $0\sspkl X\sspkl  Y$, namely $(\widehat X,\, \widehat Y)$ from {\sl Lemma 11} b), with $\widehat X\sspeq X_{i^{\prime}}$ and $\widehat Y\sspeq Y_{i^{\prime}}$, for a certain $i\sspeq i^{\prime}$ in \Eq{39}.
\pn
{\bf ii)} The solutions $(X_i,\,Y_i)$ for $i\sspgr i^{\prime}$ satisfy $0\sspkl Y_i\sspkl X_i$. \pn
{\bf iii)} The solutions $(X_i,\,Y_i)$ for $i\sspkl i^{\prime}$ satisfy $X_i\sspkl 0$ and  $|X_i|\sspgr Y_i\sspgr 0$.
\pbn
{\bf b)} Family $II$ \pn
{\bf i)} The solution $(X_i,\,Y_i)$ for $i\sspeq i^{\prime\prime}$ satisfies 
$X_{i^{\prime\prime}} \sspkl 0$ and $ 0 < |X_{i^{\prime\prime}}| \sspkl Y_{i^{\prime\prime}}$.\pn
{\bf ii)} The solutions $(X_i,\,Y_i)$ for $i\sspgr i^{\prime\prime}$ satisfy $0\sspkl Y_i\sspkl X_i$. \pn
{\bf iii)} The solutions $(X_i,\,Y_i)$ for $i\sspkl i^{\prime\prime}$ satisfy $X_i\sspkl 0 $ and  $0\sspkl Y_i\sspkl |X_i|$.
} 
\psn
{\bf Note 2:} As a necessary condition for a final solution $(X_i,\,Y_i)$ one needs later $0\sspkl X_i$, Therefore if $X_i\sspkl 0$ a common sign flip in $(X_i,\,Y_i)$ can be used, which is a proper equivalence transformation.
\pbn
{\bf Proof:} This uses the recurrences eqs. $40$ and $41$ for $D\sspeq 2,\, k\speq -s^2$ (only certain odd $s \sspgeq 3$ will later be shown to qualify) and ${\bf Auto}^{\prime}(2)\sspeq {\bf Auto}^{\prime}\sspeq Matrix([3,4],\,],[2,3])$ from \Eq{38} together with its inverse. In these recurrences the input $(X_0,\,Y_0)$ is then taken as $(X_{i^{\prime}},\, Y_{i^{\prime}})$ or  $(X_{i^{\prime\prime}},\, Y_{i^{\prime\prime}})$ from {\sl Definition 3}.
\pn
The proofs will be done by induction starting with $i^{\prime}\sspp 1$, respectively with $i^{\prime\prime}\sspp 1$, using ${\bf Auto}^{\prime}$, and starting with $i^{\prime}\sspm 1$, respectively with $i^{\prime\prime}\sspm 1$, using $({\bf Auto}^{\prime})^{-1}$. 
\pn
{\bf  a) ii)} Base case of the induction: $X_{i^{\prime}\sspp 1}\sspeq 3\,X_{i^{\prime}}\sspp 4\,Y_{i^{\prime}}\sspeq 3\,{\widehat X}\sspp 4\,{\widehat Y}$ and  $Y_{i^{\prime}\sspp 1}\sspeq 2\,X_{i^{\prime}}\sspp 3\,Y_{i^{\prime}}\sspeq 2\,{\widehat X}\sspp 3\,{\widehat Y}$ by the recurrence Eq(40) for $i\sspeq i^{\prime}\sspp 1$. Hence $0\sspkl Y_{i^{\prime}\sspp 1} \sspkl X_{i^{\prime}\sspp 1}$ because $\widehat X \sspp \widehat Y \sspgr 2\,\widehat X \sspgr 0$ from {\sl Lemma 11} b).\pn
Assuming that $ 0\sspkl Y_{i} \sspkl X_i$ for $i\sspeq i^{\prime}\sspp k$, for $k\sspeq 1,\ 2,\, ...,\,j$, it follows that $ 0\sspkl Y_{j+1} \sspeq 2\,X_j\sspp 3\,Y_j \sspkl 3\,X_j\sspp 4\,Y_j \sspeq X_{j+1}$, because from the assumption $X_j \sspp Y_{j}\sspgr 2\,Y_j \sspgr 0$.
\psn 
{\bf a) iii)} Base case of the induction: $X_{i^{\prime}-1} \sspeq 3\,X_{i^{\prime}}\sspm 4\,Y_{i^{\prime}}\sspeq 3\,\widehat X\sspm 4\,\widehat Y\sspkl -\widehat Y\sspkl 0$, and $|X_{i^{\prime}-1}| \sspeq -3\,\widehat X\sspp 4\,\widehat Y\sspgr -2\,\widehat X\sspp 3\,\widehat Y \sspeq Y_{i^{\prime}-1}$ because $-\widehat X\sspp \widehat Y \sspgr-\widehat X\sspp \widehat X\sspeq 0$.
\pn
Assuming that $X_i\sspkl 0$ and  $|X_i|\sspgr Y_i\sspgr 0$ for $i\sspeq i^{\prime}\sspm k$, for $k\sspeq 1,\ 2,\, ...,\,j$, it follows that $X_{j+1} \sspeq 3\,X_j\sspm 4\,Y_j\sspkl -4\,Y_j \sspkl 0$ and $|X_{j+1}| \sspeq  -3\,X_j\sspp 4\,Y_j \sspgr - 2\,X_j \sspp 3\, Y_j \sspeq Y_{j+1}$ because $-X_j\ \sspp Y_j\sspeq |X_j|\sspp Y_j \sspgr Y_j\sspgr 0$. 
\psn
{\bf a) i)} For solvable odd $s \sspgeq 3$ there exists the proper solution  $(X_{i^{\prime}},\, Y_{i^{\prime}})$ from {\sl Definition 3} satisfying $0\sspkl X_{i^{\prime}} \sspkl Y_{i^{\prime}}$. This is the only solution in family I as shown in parts {\bf ii)} and {\bf iii)}.  
\pbn
{\bf b) ii)} Base case of the induction: $X_{i^{\prime\prime}\sspp 1}\sspeq 3\,X_{i^{\prime\prime}}\sspp 4\,Y_{i^{\prime\prime}}\sspeq -3\,{\widehat X}\sspp 4\,{\widehat Y}$ and  $Y_{i^{\prime\prime}\sspp 1}\sspeq 2\,X_{i^{\prime\prime}}\sspp 3\,Y_{i^{\prime\prime}}\sspeq -2\,{\widehat X}\sspp 3\,{\widehat Y}$. Hence $0\sspkl Y_{i^{\prime\prime}\sspp 1} \sspkl X_{i^{\prime\prime}\sspp 1}$ because $\widehat X \sspm \widehat Y \sspkl \widehat Y\sspm \widehat Y \sspeq 0.$
\pn
Assuming that $ 0\sspkl Y_{i} \sspkl X_i$ for $i\sspeq i^{\prime\prime}\sspp k$, for $k\sspeq 1,\ 2,\, ...,\,j$, it follows that $ 0\sspkl Y_{j+1} \sspeq 2\,X_j\sspp 3\,Y_j \sspkl 3\,X_j\sspp 4\,Y_j$, because from the assumption $X_j \sspp Y_{j}\sspgr 2\,Y_j \sspgr 0$.
\pbn
{\bf  b) iii)} Base case of the induction: $X_{i^{\prime\prime}-1} \sspeq 3\,X_{i^{\prime\prime}}\sspm 4\,Y_{i^{\prime\prime}}\sspeq -3\,\widehat X\sspm 4\,\widehat Y\sspkl -7\widehat X\sspkl 0$, and $|X_{i^{\prime\prime}-1}| \sspeq 3\,\widehat X\sspp 4\,\widehat Y$. $Y_{i^{\prime\prime}-1} \sspeq -2\,X_{i^{\prime\prime}} \sspp 3\, Y_{i^{\prime\prime}} \sspeq  2\,\widehat X\sspp 3\,\widehat Y\sspgr 0$. Hence $|X_{i^{\prime\prime}-1}|\sspgr Y_{i^{\prime\prime}-1}$  because $\widehat X \sspp \widehat Y \sspgr 0$.\pn
Assuming that $X_i\sspkl 0$ and  $|X_i|\sspgr Y_i\sspgr 0$ for $i\sspeq i^{\prime\prime}\sspm k$, for $k\sspeq 1,\ 2,\, ...,\,j$, it follows that $X_{j+1} \sspeq 3\,X_j\sspm 4\,Y_j\sspkl -4\,Y_j \sspkl 0$ and $|X_{j+1}| \sspeq -3\,X_j\sspp 4\,Y_j \sspgr Y_{j+1} \sspeq -2\,X_j \sspp 3 Y_j$ because $X_j \sspm Y_j \sspkl -Y_j\sspkl 0$. 
\psn
{\bf b) i)} This follows from the {\sl Definition 3}.
\hskip 10.5cm $\square$ 
\pbn
{\bf Corollary 3: Lower bound for $\bf Y_i$ for both families}
\psn
\Beq
Y_i\sspgeq \ceil{\frac{s}{\sqrt{2}}}\,,
\Eeq
From $Y_i \sspeq \frac{1}{\sqrt{2}}\sqrt{X_i^2 \sspp s^2}\sspgr s/sqrt{2}$, because 
$\abs{X_i}\sspgr 0$\pn
For some $s\sspeq s(n) \sspeq$ \seqnum{A058529}$(n)$ this bound is attained; \eg., $n\sspeq (1),\,7,\, 17,\, 23,\, 41,\, 103,\,137,\, 193,\,...$. 
\pbn
{\bf Lemma 14: $\bf X_i$ decreases with decreasing $\bf i$}
\psn
{\it In both members of each pair of conjugate proper solutions of the generalized {\sl Pell} equation $X^2\sspm 2\,Y^2 \sspeq -s^2$, representing positive odd $s\sspgeq 3$, and positive $Y$ values, the values $X_i$ decrease with decreasing index $i$. 
}
\psn
{\bf Proof:} 
\pn
For each family $I$, \psn
1) if $i\sspgr i^{\prime}$, with $i^{\prime}$ from {\sl Definition 3},  $\triangle X_i\sspdef X_i\sspm X_{i-1}\sspeq X_i\sspm (3\,X_i\sspm 4\,Y_i)$, because by inversion of \Eq{40}, with $({\bf Auto}^{\prime})^{-1}$ from \Eq{38}. Hence    
\dstyle{\triangle X_i\sspeq 2\,Y_i\left(2\sspm \frac{X_i}{Y_i}\right)\sspgr 0} because $X_i\sspgr Y_i\sspgr 0$ from {\sl Lemma 13 a) ii)}, and \dstyle{\frac{X_i}{Y_i} \sspkl \sqrt{2}\sspkl 2.} from \Eq{45}. \psn
2) If $i\sspkl i^\prime$ then $\triangle X_i\sspdef X_i\sspm X_{i-1}\sspeq X_i\sspm (3\,X_i\sspm 4\,Y_i)\sspeq 2\,(-X_i \sspp 2\,Y_i)$ like above, but now $X_i\sspkl 0$ from {\sl Lemma 13 a) iii)}. Hence  $\triangle X_i\sspgr 0$.\psn
3) The case $i\sspeq i^{\prime}$, with $\triangle X_{i^{\prime}}\sspdef X_{i^{\prime}}\sspm X_{i^{\prime}-1}$ works like in part 1) with \dstyle{\frac{X_{i^{\prime}}}{Y_{i^{\prime}}} \sspeq \frac{\widehat X}{\widehat Y} \sspkl 1\sspkl 2} from {\sl Lemma 13 a) i)}.
\psn
For each family $II$, the claim follows then with $X^{II}_i\sspeq -X^I_{-i}$\, and $Y^{II}_i\sspeq +Y^I_{-i} $, for $i\sspgr i^{\prime\prime}$ as well as for $i\sspkl i^{\prime\prime}$, and $X^{II}_{i^{\prime\prime}}\sspeq -X^I_{i^\prime}$. \hskip 14.6cm $\square$ 
\pbn
Now the existence and uniqueness proof for each family of {\sl Proposition 7} follows.  
\pbn
{\bf Proof of Proposition 7:}\psn
If there exist proper solutions of $X^2\sspm 2\,Y^2 \sspeq -s^2$, for odd $s \sspgr 1$ then {\sl Lemma 13} shows that for each conjugate pair of families of solutions there is exactly one such solution with $0\sspkl X\sspkl Y\sspkl$, namely the positive proper fundamental solution of family I. The family pairs I and II are given in {\sl Definition 3}. The {it w.l.g.) positive $s$ has to be odd $\sspgeq 3$ as proved in {\sl Lemma 3} because even $s$ leads only to improper solutions. The case $s\sspeq 1$ where only one family with $Y\sspgr 0$ exists does not qualify because all positive proper solutions, obtained from $(X\,Y)\sspeq (1,\,1)$, by applying powers of ${\bf Auto}^{\prime}$ from \Eq {38}, satisfy $X\sspgr Y\sspgr 0$, as is seen by induction.  
\psn
An aside: One can use the \rpapf method \Eq{39} to find this trivial $(1,\,1)$ solution for $s\sspeq 1$. From the parallel form $F_{pa}(j\sspeq 0)\sspeq [-1,\,0,\,2]$ one needs three $\bf R$-transformations with $\vec t\sspeq (0,\,-1,\,2)$ to reach $F_p\sspeq [1, 2,\,-1]$, leading with $i\sspeq 1$ and $-{\bf B}$ and $-{\bf Auto}$, to the positive  {\sl Pell} solution $(1,\, 1)$ representing $-1$.
\psn
What is left is to investigate for which odd $s\sspgeq 3$ values proper solutions of $X^2\sspm 2\,Y^2\sspeq -s^2$ exist. For this all solutions of the \rpapfs of {\sl Lemma 7} for the present $D\sspeq 2$ case have to be found, which means to solve the congruence $j^2\sspm 2\sspequiv \modan{0}{-s^2}$, \ie $j^2\sspm 2\sspequiv \modan{0}{s^2}$, for $s\sspeq \prod_{i=1}^k p_{i}^{\beta_i}$, with odd primes $p_i$ and $\beta_i\sspgeq 1$. 
\pn
This quadratic congruence is reduced to the congruence with odd prime power moduli $p_i^{\alpha_i}$, with $ \alpha_i\sspeq 2\,\beta_i$. The number of \rpapfs is then $\nu(s^2)\sspeq \prod_{i=1}^k \nu(p_i^{\alpha_i})$ with the number of solutions to the congruence $j^2\sspm 2 \sspequiv \modan{0}{p_i^{\alpha_i}}$. See \eg \cite{Apostol}, {\it Theorem 5.28}, pp. 118, 119. In the second step this congruence is treated by the lifting theorem \cite{Apostol}, {\it Theorem 5.30}, pp. 121 - 122.
\pn
One starts with finding the odd primes for which $2$ is a quadratic residue modulo $p$. This is easy because for odd  $p$ the {\sl Legendre} symbol $L(2,\,p) \sspeq (-1)^{\frac{p^2-1}{8}}$. See, \eg \cite{Apostol}, {\it Theorem 9.5}, p. 181. Thus, $L(2,\,p)\sspeq +1$ for $p$ congruent to $\modan{-1}{8}$ or $\modan{+1}{8}$ \seqnum{A001132}. From the degree of the congruence for each of these primes there are two solutions $j\sspeq r$ and  $j\sspeq p\sspm r$, with $r\sspin\{0,\, 1,...,\, p-1\}$. 
\pn
Now the lifting of the congruence to a power $p^k$ is obtained inductively from $p^{k-1}$ to $p^k$, for k \sspgr 2. In our case the lifting is unique for each solution; part a) of the lifting theorem. This is because if there is a solution $r$ of $j^2\sspm 2\sspequiv \modan{0}{p^{k-1}}$ for $0\sspleq r\sspkl p^{k-1}$ then $2\,r\sspnotequiv\modan{0}{p}$ because otherwise $(2\,r)^2\sspeq (K\,p)^2$, for $K\sspin \mathbb Z$, but then $8\sspeq p\,(K^2\,p - 4\,L\,p^{k-2})$, with $L\sspin \mathbb Z$, \ie $8\sspequiv \modan{0}{p}$; a contradiction for odd $p$. Thus, the two solutions for each $p\sspequiv \modan{1}{8}$ or $p\sspequiv \modan{7}{8}$ are lifted to two solutions for any modulus $p^k$. \pn
In our case, for $s^2$, only even powers $\alpha\sspeq 2\,\beta$ appear. \pn
In the proof of {\sl Lemma 8} part b) it is shown that all \rpapfs reach the (only) fundamental $2$-cycle. Therefore the number of solutions corresponds to the number of \rpapfs for the generalized {\sl Pell} equation for {$Disc\sspeq 8$ and $k\sspeq -s^2$. This number is  $2^{P_1\sspp P_7}$, where $P_1$ and $P_7$ are the number of prime divisors of $s$ congruent to $1$ and $7$ modulo $8$, respectively. The sequence of these solvable $s$ is \seqnum{A058529}$(n)$, for $n \sspgeq 2$.
\pn 
Each \rpapf leads to a solutions by \Eq{39} for $D\sspeq 2$, where only solutions with  $Y\sspgr 0$ are taken into account. Of course, there are also the obvious solutions with $Y\sspeq -0$, not considered here.  
\pn
For general quadratic congruences see \cite{Apostol} chapter 9.1, pp. 178 ff. 
\pn
Because the \rpapfs come in conjugate pairs, leading to families of solutions ({\sl Lemma 8} and {\sl Definition 3}), and only family I contains the unique solution with $0\sspkl X\sspkl Y$ there are only $2^{P_1\sspp P_7\sspm 1}$ solutions to \Eq{35}.
} 
\hskip 15.5cm $\square$
\pbn 
For these $s\sspeq s(n)\sspeq$\seqnum{A058529}$(n)$, with $n \sspgeq 2$, see {\it Table 3}, for the $s-$ values $7,\,17,...,\,313$, also for $(x\sspeq c_1,\, y\sspeq c_2)$ solutions , together with $q\speq c_3$, and the $t$-tuples for the {\it rpapf}s needed in \Eq{39}. The $t$-tuples come in conjugate pairs, but only one of them leads to family I solutions with the required unique $0\spkl X\sspkl Y$ solution. The other $t$-tuple is given within solid brackets.
Composite $s$ values (with also improper solutions) are underlined.
\pbn
\pbn
{\bf Case ii) $\bf q\sspkl c_3$}
\psn
We use $q\sspeq c_3\sspm k$, with $1\sspleq k \sspkl \sqrt{2\,x^2 \sspp y^2}$, where $x\sspeq c_1 < c_2 \sspeq y\sspkl q\sspp k$, and  $gcd(x,\,y,\,q\sspp k) \sspeq 1$.\psn
This case will lead to two binary quadratic fields, one indefinite with $ Disc\sspeq 8$, the other a positive definite form with $Disc\sspeq -8$, representing a negative integer $-a$ and the positive number $a$, respectively. The proper solutions will be found by the method of the representative parallel primitive forms ({\it rpapfs}) like above. Improper cases are reduced to proper ones.
\psn
{\bf Lemma 15: Quadratic equation for $\bf (x,\, y)$} 
\psn
{\it 
\Beq 
x^2 \sspp y^2 \sspp 6\,x\,y\sspp  4\,(x\sspp y)\,k\sspeq t^2, 
\Eeq
{\it w.l.g.}\ $t\sspgeq 0$.
}
\pn
Only solutions with $0\sspkl x\sspkl y$ are of interest.
\psn
{\bf Proof:} 
From the definitions
\Beq
q^2 \sspeq x\,y \sspp (q\sspp k)\,(x+y),
\Eeq
leading to the positive solution (the negative one is excluded because $q\sspgr 0$)
\Beq
2\,q\sspeq x\sspp y \sspp \sqrt{x^2 \sspp y^2 \sspp 6\,x\,y\sspp 4\,k\,(x\sspp y)} \sspfed x\sspp y + t,
\Eeq
Integer solutions of $q$ need the claimed formula. \hskip 9cm $\square$
\psn
Diagonalization of \Eq{50} leads, with $X\sspeq y\sspm x$, $Y\sspeq y\sspp x$ and $\widehat Y\sspeq Y\sspp k$, to
\psn
{\bf Lemma 16: Quadratic form for $\bf (X,\, \widehat Y)$}
\psn
\Beq
X^2 \sspm 2\,\widehat Y^2\sspeq -a(t,\,k)\sspdef -(t^2\sspp 2\,k^2)\, .
\Eeq
Only solutions satisfying $0\sspkl X\sspkl \widehat Y\sspm k$ are of interest, because $2\,x \sspeq Y\sspm X$, $2\,y\sspeq Y\sspp X$ and $Y\sspeq \widehat Y \sspm k$.   
\psn  
{\bf Proof:} Elementary. \hskip 13.7cm $\square$
\psn
{\bf Corollary 4:} \psn
\Beq
2\,q \sspeq \widehat Y\sspp t\sspm k,\ {\rm and} \  2\,c_3 \sspeq 2\,(q\sspp k) \sspeq \widehat Y\sspp t\sspp k.  
\Eeq
The left-hand side of \Eq{53} shows an indefinite binary quadratic form of discriminant $Disc \sspeq 8$ representing a negative integer $-a$, and the right-hand side is a positive definite form of $Disc \sspeq -8$ representing a positive integer $a$.
\psn
From the use of $t$ in \Eq{52} it is clear that $t\sspgr 0$. One can show that $t\sspeq 0$ is indeed excluded in \Eq{53}.
\psn
{\bf Lemma 17: Exclusion of instance $\bf t\sspeq 0$.}
\psn
{\it There is no solution for $(X,\,\widehat Y)$ and $0\sspkl X\sspkl Y$ if $t\sspeq 0$}. 
\psn
{\bf Proof:} If $t\sspeq 0$ then $X$ is even, $X\sspeq 2\,\widetilde X$, \ie $\widetilde Y^2\sspm  k^2 \sspeq 2\,\tilde X^2$, $2\,(\widehat Y\sspp k)\,(\widehat Y\sspm k)\sspeq X^2$, \ie $\sqrt{2\,(Y\sspp 2\,k)\,Y}\sspeq X$, hence $X\sspgr \sqrt{2\,Y\,Y} \sspgr Y\sspgr 0$, contradicting the required inequality.
\hskip 1.5cm $\square$
\psn
{\bf Note 3: Symmetry between $\bf X$ and $\bf t$}

\psn
In the equation  $X^2\sspm 2\,\widehat Y^2 \sspp t^2 \sspp 2\,k^2 \sspeq 0$ there is a symmetry between $X$ and $t$. However, if a solution satisfies $0\sspkl X \sspkl Y\sspeq \widehat Y \sspm k$ then $0\sspkl t \sspkl Y$ cannot be satisfied for $k\sspgr 0$. Because then $0\sspeq X^2\sspm 2\,\widehat Y^2 \sspp t^2\sspp 2\,k^2 \sspkl (\widehat Y\sspm k)^2 \sspm 2\,\widehat Y^2 \sspp (\widehat Y\sspm k)^2 \sspp 2\, k^2 \sspeq -4\,\widehat Y\,k \sspp 4\,k^2\sspeq  -4\,k\,Y\ \sspkl 0$, a contradiction.
\psn
The possible values $a$ for the two quadratic forms are now determined one after the other, and then the intersection of these $a$ values are of interest as candidates for testing the required $(X,\,Y)$ inequalities.
\psn
{\bf A) The indefinite binary quadratic form}
\psn
{\bf Lemma 18: All values $\bf a$ for proper solutions of $\bf X^2\sspm 2\,\widehat Y^2\sspeq -a$}
\psn
{\it {\bf a)} The indefinite quadratic form $X^2\sspm 2\,\widehat Y^2$ represents $-a$, for $a\sspgeq 2$, properly if and only if,
\Beq
a\sspeq 2^{e_2}\,\prod_{i\sspeq 1}^{P_1} p_{1,i}^{e_{1,i}}\,\prod_{i\sspeq 1}^{P_7} p_{7,i}^{e_{7,i}},
\Eeq
with $e_2\sspin \{0,\,1\}$, $e_{1,i} \sspgeq 0$, $e_{7,i} \sspgeq 0$, $P_1\sspgeq 0$ and $P_7 \sspgeq 0$, but not all exponents vanish. The products are over powers of odd primes with $p_{1,i}$ congruent to \modan{1}{8}, \seqnum{A007519}, and $p_{7,i}$ are the primes congruent to \modan{7}{8}, \seqnum{A007522}. The combination of these odd primes is shown in \seqnum{A001132}, together with prime $2$ in \seqnum{A038873}. \pn
$-a\sspeq -1$ is also represented by one family (the ambiguous case) with $\widehat Y\sspgr 0$ with the fundamental solution $(1,\,1)$. This will not be of interest because $0\sspkl X\sspkl Y\sspeq \widehat Y\sspm k$ is required later.
\psn
{\bf b)} The number of {\it rpapfs}, \ie the number of families of solutions with positive integer $\widehat Y$, is $f(-a)\sspeq 2^{P_1 \sspp P_7}$.
\psn
{\bf c)} The sequence of these numbers $a$ is given in \seqnum{A057126}.
}
\psn
{\bf Proof:} \psn
{\bf a)} It is sufficient to find the $a$ values for proper representations $(X,\,\widehat Y)$. The improper solutions reduce to these solutions. Improper solutions are later also of interest because $\gcd(X,\,\widehat Y\sspeq Y\sspp k)$ is not necessarily $1$. \pn
All proper solutions are found using the representative parallel primitive forms ({\it rpapfs}), See {\sl Lemma 7} for $D\sspeq 2$ and $k\sspeq -a$. They are given by solving the congruence $j^2 \sspm 2\sspequiv \modan{0}{-a} \sspequiv \modan{0}{a}$ for $j\sspin \{0,\,1,\,...,\,a\sspm 1\}$ as $F_{pa}(-a,\,j)\sspeq [-a,\,2\,j,\,(2\sspm j^2)/a]$. Each solution $j$ of this congruence leads by $\bf R$-transformations always to the principal form $F_p\sspeq [1,\,2,\,-1]$ (see the proof of {\sl Lemma 8} b)). Therefore each such $F_{pa}(-a\, j)$ produces a family of proper representations of $-a$ by $X^2\sspm 2\,\widehat Y^2$. These infinite $j$-families are obtained by \Eq{39} (with $D\sspeq 2$ and $k\sspeq -a$), restricting to positive $\widehat Y$ by a possible overall sign flip in $(X,\, \widehat Y)$, The number of such families is called f(-a) (there are then also solutions with $\widehat Y\sspto -\widehat Y$, of course).
\pn 
The congruence for $j$ is solved by using the prime number decomposition $a\sspeq 2^{e_2}\,\prod_{i=1}^K p_i^{e_i}$, with odd primes $p_i$. The solution for each prime number is then determined. 
\pn 
If $a$ is even the solution of $j^2\sspm 2\sspequiv \modan{0}{(p\sspeq 2)}$ is $j\sspeq \pm 0$. A lifting to powers of prime $2$ is not possible by the lifting theorem \cite{Apostol}, {\it Theorem 5.30}, pp. 121 - 122, because there case $(b), (b_2)$ applies inductively for each power $\alpha \sspgr 2$, since $2\cdot 0\sspequiv \modan{0}{2}$. Therefore $e_2\sspin\{0,\,1\}$.
\pn 
For each odd prime $p$ the congruence $j^2\sspm 2\sspequiv \modan{0}{p}$ has two solutions $ j\sspeq r$ and $j\sspeq p\sspm r$, for a certain $r\sspin \{0,\,1,\,...,\,p-1\}$, only if $p\sspequiv \modan{\pm 1}{8}$ from the {\sl Legendre} symbol $(2,\,p)\sspeq +1$. In this case there is an inductive proof of a unique lifting of each of the two solutions to any power $k\sspgeq 2$ of such a $p$. Adapt this proof to the one given in the proof of {\sl Proposition 7}. 
\psn
{\bf b)} The number $f(-a)$ of solutions of $j$, \ie the number of {\it rpapfs}, is 
\pn
$f(-a)\sspeq \nu (2)\,\prod_{i\sspeq 1}^{P_1} \nu(p_1,\,i)\,\prod_{i\sspeq 1}^{P_7}\nu(p_7,\,i)\sspeq 2^{0\sspp P_1\sspp P_7}$, where $\nu(p)$ is the number of solutions of the congruence for $j$ for prime $p$, a factor of $a$, viz $\nu(2)\sspeq 1$ and $\nu(p)\sspeq 2$ for each odd $p$. 
\psn
{\bf c)} This follows from the definition of \seqnum{A057126}. See the first comment there, and the following conjecture which is true because  $X^2 \sspm 2\,\widehat Y \sspeq -a$ has the same numbers $a$ as proper representations of the form $2\,x^2\sspm y^2$ (see a Nov 09 2009 comment in \seqnum{A035251} where also the $a$ values with only uninteresting improper representations are listed). \hskip 13.5cm $\square$
\pbn
\vfill\eject\noindent
{\bf Example 9: Proper solutions for $\bf a\sspeq 2\cdot7\cdot 17\sspeq 238$}
\psn
The number of {\it rpapfs} $F_{pa}(-238,\, j)$ is $f(-238)\sspeq 2^{1\sspp 1}\sspeq 4$, viz $[-238, \,148,\,-23],\,[-238,\,216,\,-49]$\, $[-238,\,260,\,-71],\,[-238,\,328,\,-113]$, for $j\sspeq 74,\,108,\,130,\,164$, respectively. The corresponding $\vec t$ values are $(-3,\,4,\,2),\, (-2,\,5,\,1),\,(-1,\,2,\,2,\,2,\,2,\,2,\,2)$ and $(-1,\,3,\,2,\,2,\,2,\,1)$. The corresponding fundamental  solutions $(X_0,\,\widehat  Y_0)$ are from \Eq{39}, with $Y\sspeq \widehat Y$, $D\sspeq 2$, $k\sspeq -238$ and ${\bf B(2)}\sspeq -{\bf B}\sspeq Matrix([1,1],\,[0,1])$ and ${\bf Auto}(2)\sspeq -{\bf Auto}$ from \Eq{37}, $(-10,13),\,(2,\,11),\,(-2,\,11),\,(10,\,13)$. The only solutions satisfying $0\sspkl X\sspkl \widehat Y$ in these four families are $(2,\,11)$ and $(10,\,13)$. \pn
However, these solutions will not survive as final solution of \Eq{53} because the solutions for $(t,\,k)$ will have no solution for $a\sspeq 238$. See the later {\sl Lemma 18}, which shows that there is no solution with prime $p\sspequiv \modan{7}{8}$.
\psn
{\sl Lemma 11} can be rewritten with $Y\sspto \widehat Y$, and $s\sspeq \sqrt a$, with the positive proper fundamental solution (\it ppfs) $(X^{\prime},\, \widehat Y^{\prime})$ as
\psn
{\bf Corollary 5:}\psn
\it{\bf a)} If there is a solution $(X_i,\,{\widehat Y}_i)$ for the generalized {\sl Pell} equation $X^2\sspm 2\,\widehat Y^2 \sspeq -a$, for $a\sspin \mathbb N$, with ${\widehat Y}_i\sspgr 0$, then
\Beq
\frac{|X_i|}{{\widehat Y}_i} \sspeq \sqrt{2}\,\sqrt{1\sspm \frac{1}{2}\,\left(\frac{\sqrt{a}}{\widehat Y_i}\right)^2} \sspeq \sqrt{2}\bigg / \sqrt{1\sspp \left(\frac{\sqrt{a}}{X_i}\right)^2} \sspkl \sqrt{2}\,.
\Eeq
{\bf b)} The positive proper fundamental solution ({\it ppfs}) $(X^{\prime},\,\widehat Y^{\prime})$ for one of the members of a pair of conjugate families of solutions of this {\sl Pell} equation representing $-a$ satisfies 
\Beq
0\sspkl X^{\prime}\sspkl \widehat Y^{\prime}\sspkl \sqrt{a}\,.
\Eeq }
The definition of conjugate family pairs $I$ and $II$ is taken over from {\sl Definition 3} with $s\sspto \sqrt{a}$.
Hence family $I$ is the one which has the {\it ppfs} solution $(X^{\prime},\,\widehat Y^{\prime})$. The conjugate family has then the solution $(-X^{\prime},\,\widehat Y^{\prime})$.
\pn
In {\sl Example 9} there are two conjugate family pairs. The first pair has in family $I$ the {\it ppfs} $(2,\, 11)$, and $(-2,\,11)$ in family $II$. The other pair has {\it ppfs} $(10,\, 13)$ in family $I$, and $(-10,\,13)$ in family $II$.
\psn
Adapting {\sl Definition 3}, {\sl Lemma 13} and {\sl Lemma 14} to the present case, with $s\sspto \sqrt{a}$ and $Y\sspto \widehat Y$, $\widehat X\sspto X^{\prime}$ and $ \widehat Y\sspto \widehat Y^{\prime}$,shows that there is in each conjugate pair of solutions precisely one solution, namely ($X^{\prime},\,\widehat Y^{\prime})$ in family $I$, for each represented $a$ value of \Eq{55}, satisfying  the necessary $0\sspkl X\sspkl \widehat Y$ requirement.\pn
Therefore the number of solutions of $X^2\sspm 2\,\widehat Y^2\sspeq -a$, $a$ from \Eq{55}, also satisfying $0\sspkl X\sspkl \widehat Y$ has precisely $2^{P_1\sspp P_7\sspm 1}$ solutions.
\psn
But in order to obtain the final result for \Eq{53} one has also to solve the definite quadratic form for $(t,\,k)$ of \Eq{53} representing certain positive integers $a$ (not the one of \Eq{55}).
\pbn
{\bf B) The positive binary quadratic form}
\psn
Before treating the special $(t,\,k)$ case some general remarks on definite binary quadratic forms are collected. Compare with \cite{SchSch}, \S 29 and \S 30, 104 - 112.
\pn
A definite quadratic form $F\sspeq [a,\,b,\,c;x,\,y]$, with discriminant $Disc \sspeq -(4\,a\,c\sspm b^2)\sspkl 0$ is positive or negative definite, \ie $F\sspgr 0$ or $F\sspkl 0$ for all integers $(x,\, y)\sspneq (0,\,0)$, if $a \sspgr 0$  or  $a \sspkl 0$. This implies that $c\sspgr 0$ or $c\sspkl 0$, respectively. This follows from the identity (see \cite{SchSch}, p. 102, \Eq{126}). $4\,a\,F\sspeq (2\,a,x\sspp b\,y)^2 \sspp \abs{Disc}\,y^2)$.
\pn
An aside: Note that $Disc\sspeq -q^2\sspeq (i\,q)^2$, for $q\sspgeq 0$, is not considered because $a\,F$ factorizes into a complex conjugate pair: $a\,F\speq K\,
\overline K$, with $K\sspeq (a\,x\sspp \frac{b}{2})\sspp \frac{q}{2}\,y\,i\sspdef A\sspp B\,i$. (Compare \cite{SchSch} \Eq{127}, p. 103, for the indefinite $Disc\sspeq q^2\sspgr0$ case). For the positive definite $F$ representing $k\sspin \mathbb N$ this leads to the representations of two squares $A^2\sspp B^2\sspeq a\,k\sspin \mathbb N$. For this see \cite{Grosswald}, \S 3, Theorem 3, p. 15, or \cite{MorenoWagstaff}, Theorem 2.5, p. 27.
\pn
A positive definite Form is called {\it reduced} if $0\sspleq\abs{b}\sspleq a\sspleq c$. This implies that $0\sspkl a\sspleq \floor{\sqrt{\frac{-Disc}{3}}}$, hence there is only finite number of solutions for $a$ hence for $b$ ($-a\sspleq b\sspleq a)$ with given $Disc$ (determining $c$). Thus, the number of these reduced forms for given $Disc \sspkl 0$ is finite, see \cite{Buell}, {\sl Theorem 2.2.}, p. 13.
\psn
We call a positive definite quadratic form $F$ {\it half-reduced} if $0\sspeq \abs{b}\sspgeq a$. Hence such an $F$ is then not reduced if $c\sspgr a$.
\psn
A $t$-family of {\it parallel forms}, each improperly (with a $Det\sspeq -1$ transformation) equivalent to $F\sspeq [a,\,b,\,c]$ is defined by $F^{\prime}(t) \sspeq [a,\, b^{\prime}(t),\, c^{\prime}(t)]$, with $b^{\prime}(t)\sspeq -b \sspp 2\,a\,t$ and $c^{\prime} \sspeq c\sspm b\,t \sspp a\,t^2$, for $t\sspin \mathbb Z$. This corresponds to ${\bf A}^{\prime}\sspeq ({\bf N}\,\bf R(t))^\top\,{\bf A}\,({\bf N}\,\bf R(t))$, with ${\bf N}\sspeq Matrix([0,1],\,[1,0])$, $Det\,{\bf N}\sspeq -1$, and ${\bf R}(t)\sspeq  Matrix([1,-1],\,[0,t])$, $Det\,{\bf R}\sspeq +1$, as used earlier in {\sl Lemma 6}.
\psn
Each primitive positive definite form $[a,\, b,\,c]$ representing properly, with $\vec x\sspeq (x,\, y)^{\top}$, an integer $k\sspgr 0$ (for the positive definite case) is improperly equivalent to a $t$-family of parallel forms $F^{\prime\prime}(t)$ obtained after two steps. \pn 
Start with $F\sspeq F_1$ and $\vec x_1$, then $F^{\prime}\sspeq [k,\,b^{\prime},\, c^{\prime}]$, with $Disc\, F \sspeq Disc\, F^{\prime}$ and  ${\bf A^{\prime}\sspeq M^\top {\bf A}\, {\bf M}}$, where ${\bf M}\sspeq Matrix([x_1,v_1],\,[y_1,w_1])$, with $Det\, {\bf M}\sspeq +1\sspeq x_1\,w_1\sspm y_1\,v_1$. $b^{\prime}$ and $c^{\prime}$ look complicated, not needed explicitly here. Thus, $\vec x^{\,\prime}\sspeq {\bf M}^{-1} \vec x_1 \sspeq (1,\,0)^\top $. \pn
In a second step a $t$-family of parallel forms $F^{\prime\prime}(t)\sspeq [k,\,b^{\prime\prime},\,c^{\prime\prime}]$ is obtained improperly \via \pn 
${\bf A}^{\prime\prime}\sspeq {\bf T}^\top(t)\,{\bf A^{\prime}}\,{\bf T}(t)$, where from above ${\bf T}(t)\sspeq  {\bf N}\,{\bf R}(t)\sspeq Matrix([1,t],\,[0,-1])$. Thus, $\vec x^{\,\prime\prime}\sspeq {\bf T}^{-1}\,\vec x^{\,\prime} \sspeq  {\bf T}\,\vec x^{\,\prime} \sspeq (1,\,0)^{\top}$. Now $b^{\,\prime\prime}(t)\sspeq -b^{\prime} \sspp 2\,k\,t$ and $c^{\,\prime\prime}(t)\sspeq c^{\prime}\sspm b^{\prime}\,t\sspp k\,t^2$.
\pn 
Because $b^{\prime\prime}\sspequiv = -\modan{b^{\prime}}{2\,k}$, a representative parallel form $F^{\prime\prime}\sspeq F^{\prime\prime}(t_0)$ with $b^{\,\prime\prime}(t_0)\sspin [0,\,2\,k)$ can be chosen. \pn
This leads to a program for all representative parallel positive definite forms forms ({\it rpaf})s $F$ (forget the superscripts $\prime\prime$) of given $Disc\sspkl 0$ representing an integer $k\sspgr 0$. We are interested later only on even $Disc$, \ie $Disc \sspeq -4\,D$. This becomes then a $j-$family of {\it rpaf}s, \ie \dstyle{F(j)\sspeq \left [k,\, 2\,j,\, \frac{j^2\sspp D}{k}\right]}, for $j\sspin [0,\,k-1]$. This is because from the Disc formula $b$ is even, $b\sspeq 2\,j$ and $b\sspin [0,2\,k)$ ($b^{\,\prime\prime}(t_0)$ from above).\pn
A representative parallel form $F^{\prime\prime}$ that is improperly equivalent to $F$ 
 may not be primitive, even if $F$ was primitive. Therefore one discards in the $j$-family of the {\it rpaf}s the imprimitive forms and stays with a set of primitive ones the {\it rpapf}s.
\psn
{\bf Example 10: Representative parallel forms for some Disc and $\bf k$}
\psn
1) The representable $k$ values for $Disc\sspeq -32,\, D\sspeq 8$ are, including also imprimitive forms,   
$1,\, 2,\, 3,\, 4,\, 6,\, 8,\, 9,\, 11,\, 12,\,\, 17, \,18,\, 19,\, 22, \,24,\, 27,\, 33,\, 34,\, 36,\, 38,\,...$. Those with primitive forms (\rpapfs) are $1,\, 3,\, \,4,\,8,\,9,\,11,\,12,\,17,\,19,\,24,\, 27,\, 33,\, 36,\,...$\,.\pn
The number of $j$-values for these \rpapfs are $1,\, 2,\, 1,\, 2,\, 2, \,2,\, 2,\, 2,\, 2, \,2,\, 2,\, 2,\, 4,\, 2,\,...$\,. \pn
For $k\sspeq 2$ there is only the imprimitive $j\sspeq 0$ {\it rpaf} $[2, 0, 4]$. \pn
For $k\sspeq 4$ the two {\it rpaf}s with $j\sspeq 0$ and $2$ are $[4, 0, 2]$ and $[4, 4, 3]$, but only $[4, 4, 3]$ is primitive, the \rpapf\, for $k\sspeq 4$. 
\psn
2) The representable $k$ values for $Disc\sspeq -20,\, D\sspeq 5$ are given in \seqnum{A343238}: $1,\, 2,\, 3,\, 5,\, 6,\, 7,\, 9,\, 10,\, 14,\, 15,\,$
\pn $18,\,...$\,. The number of $j$-solutions is given in \seqnum{A343240}. It seems that for each such $k$ value all $j$ solutions lead to primitive parallel forms (tested for the first $70$ of these $k$-values).\pn
For $k\sspeq 18$ the two $j$-values are $7$ and $11$ with corresponding \rpapfs $[18,\, 14,\, 3]$ and $[18,\, 22,\, 7]$.
\psn
A $t$-family of {\it right neighbor forms} is defined by $\widetilde F(t) \sspeq [c,\,-b\sspp 2\,c\,t,\,a\sspm b\,t\sspp c\,t^2]$ this is obtained like above in the parallel case but omitting $\bf N$, hence exchanging $a$ and $c$. This is a proper equivalence transformation $\bf R(t)$. 
\pn
To obtain a unique {\it half-reduced right neighbor form} ({\it hredrnf}) one chooses \dstyle{t\sspeq \ceil{\frac{b\sspm c}{2\,c}}}, from the half-reduced definition.
\psn
Now we treat the reduced positive definite quadratic form  $t^2 \sspp 2\,k^2$, with $Disc\sspeq -8$, representing an integer $a\sspgeq 1$. This is the reduced form $[1,\,0,\,2]$ which is the only one for $Disc\sspeq -8$. See the Example table in \cite{Buell}, p. 18 (where $\triangle = Disc$; class number $h(-8)\sspeq 1$). Therefore there is only this $1$-cycle of reduced forms. \pn
The \rpapfs are determined from the congruence $j^2\sspp 2\sspequiv \modan{0}{a}$ for $j\sspin\{0,\,1,\,...,\,a-1\}$ by $F_{pa}(a,\,j)\sspeq [a,\,2\,j,\,(j^2\sspp 2)/a]$. Every {\it rpapf} leads \via {\bf R}-transformations to this primitive form $[1,\,0,\,2]$. That every definte form leads via a chain of right neighbor forms to a reduced one is proven in \cite{SchSch} Satz 75, \S 38, pp. 108 - 109. As shown above in the proof of {\sl Lemma 7}, part {\bf b)}, {\bf R}-transformations respect primitivity of forms.
The formula for the $(t,\,k)$ solutions for given $a\sspgeq 1$, obtained from the {\it rpapfs} with the $\vec t$ values for these ${\bf R}$-transformations, is here (compare this with \Eq{39})
\Beq
\LTvez{t(j(a))} {k(j(a))}\sspeq \pm{\bf Rpa}(\vec t(j(a)))^{-1}\,\LTvez{1}{0}, 
\Eeq
with the solutions of the $j-$congruence determining the {\it rpapfs}. A sign is chosen to obtain $k\sspgr 0$.
\psn
{\bf Example 11: Chain of $\bf R$-transformations}
\psn
For $a\sspeq 17$ the two \rpapfs are $[17,\,14,\,3]$ and $[17,\,20,\,6]$, for $j\sspeq 7$ and $10$. The $\vec t$-tuples are $(2,\,-1)$ and $(2,\,2)$ with the chains $\{[17, 14, 3],\,[3, -2, 1],\, [1, 0, 2]\}$  and $\{[17, 20, 6],\,[6, 4, 1],\,[1, 0, 2]\}$.
\psn
{\bf Lemma 19: All values $\bf a\sspgr 0$ for proper solutions of $\bf t^2\sspp 2\,k^2\sspeq a$}
\psn
{\it {\bf a)} The positive definite quadratic form $t^2\sspp 2\,k^2$, with $Disc\sspeq -8$, represents an integer $a$, for $a\sspgeq 2$, properly if and only if,
\Beq
a\sspeq 2^{e_2}\,\prod_{i\sspeq 1}^{P_1} p_{1,i}^{e_{1,i}}\,\prod_{i\sspeq 1}^{P_3} p_{3,i}^{e_{3,i}},
\Eeq
with $e_2\sspin \{0,\,1\}$, $e_{1,i} \sspgeq 0$, $e_{3,i} \sspgeq 0$, $P_1\sspgeq 0$ and $P_3\sspgeq 0$, but not all exponents vanish. The products are over powers of odd primes with $p_{1,i}$ congruent to \modan{1}{8}, \seqnum{A007519}, and $p_{3,i}$ are the primes congruent to \modan{3}{8}, \seqnum{A007520}. The combination of these odd primes is shown in \seqnum{A033200}, together with prime $2$ in \seqnum{A033203}. 
\pn
The case $a\sspeq 1$ is also represented properly but the solutions $(\pm 1,\,0)$ is not of interest because $k$ should be positive.\psn
{\bf b)} The number of \rpapfs, \ie the number of families of solutions with positive integer $k$, is $f(a)\sspeq 2^{P_1 \sspp P_3}$ (for $a\sspgeq 2$). These \rpapfs appear as conjugate pairs for values $j$ and $a\sspm j$, differing in the sign of $t\sspgr 0$. 
}
\psn
{\bf Proof:}
\psn
{\bf a)} One has to solve the congruence $j^2\sspp 2\sspequiv \modan{0}{a}$. See the proof of {\sl Lemma 17 a)} where now, besides for prime $2$ treated there and applying here too, one has to solve for odd primes $p$ the congruence $j^2\sspp 2\sspequiv \modan{0}{p}$, \ie to compute the {\sl Legendre} symbol $L(-2,\,p)\sspeq +1$. By multiplicativity of this symbol this means  $L(-1,\,p)\,L(2,\,p) \sspeq +1$. The solution for signs $++$ is $p\sspequiv \modan{1}{4}$ and $p\sspequiv \modan{\pm 1}{8}$, \ie $p\sspequiv \modan{+1}{8}$, because the other case means $p\sspeq 1\sspp 4\,K\sspeq 7\sspp 8\,L$, with integers $K$ and $L$, leading to the contradiction  $3\sspeq 2\,(K\sspm 2\,L)$. The other solution with signs $--$ is $p\sspequiv \modan{3}{4}$ and $p\sspequiv \modan{\pm 3}{8}$, \ie $p\sspequiv \modan{+3}{8}$, because the other case means $p\sspeq 3\sspp 4\,K\sspeq 5\sspp 8\,L$, with integers $K$ and $L$, leading to the contradiction  $1\sspeq 2\,(K\sspm 2\,L)$.
Therefore one find $L(-2,\,p)\sspeq +1$ for $p\sspequiv \modan{1}{8}$ or $p\sspequiv \modan{3}{8} $.\pn
The unique lifting to powers of these primes for each of the two solutions of the $j-$congruence is proved inductively like in the proof of Proposition 7.
\psn
{\bf b)} Follows like in {\sl Lemma 18 b)} where also the \rpapfs appeared in conjugate pairs. \hskip 11.8cm $\square$
\psn
{\bf Example 12: Proper solutions for $a\speq 3^2\,17\speq 153$}
\psn
The number of {\it rpapfs} $F_{pa}(153, j)$ is $f(153)\sspeq 2^{1\sspp 1}\sspeq 4$, viz $[153,\, 82,\, 11],\,[153,\, 116,\, 22]$\, $[153,\, 190,\, 59]$ 
\pn $[153,\, 224,\, 82]$, for $j\sspeq 41,\,58,\,95,\,112$, respectively. Observe the two conjugate $j$-pairs $(41,\, 112)$ and $(58\, 95)$, and the corresponding pairs of \rpapfs. The $\vec t$ values are $(4,\,3)$, $(3,\,3,\,1)$, $(2,\,3,\,2,\,0)$, $(1,\,-3,\,-3)$, respectively. The solutions are  $(t, k\sspgr 0)$ are $(-11,\,4)$, $(-5,\,8)$, $(5,\,8 )$, $(11,\,4)$, respectively. The pairs of conjugate solutions differ by the sign of $t$, Later only the solution $(11,\,4)$ will survive if also the $(X,\,Y)$ values for the indefinite form representing $-a$ with $0\sspkl X\sspkl Y$ are determined.
\pbn
{\bf C) Combining solutions of sections A) and B)}
\psn 
The proper representations of $X^2\sspm 2\,{\widehat Y}^2$ of certain $-a$ of section {\bf A)} ({\sl Lemma} 17) and those of $t^2\sspp 2\,k^2$ for certain $a$ of section {\bf B)} for the case {\bf ii)} with $q\sspeq c_3\sspm k$, are combined to find 
the representations of \Eq{53} ({\sl Lemma} 15), and then the requirement $ 0\sspkl X\sspkl Y\sspeq \widehat Y \sspm k$, for $1\sspleq k\sspkl \sqrt{X^2\sspp Y^2}$ is enforced. Because $Y\sspeq \widehat Y\sspm k$ this uses solutions of the form for $(t,\,k)$. Finally, $c_1\speq x\sspeq \frac{Y\sspm X}{2}$. $c_2\sspeq y\sspeq \frac{Y \sspp X}{2}$, and $c_3\sspeq \frac{\widehat Y\sspp t\sspp k}{2}$ have to be checked for $\gcd(c_1,\,c_2,\,c_3)\sspeq 1$.
\pbn  
{\bf Proposition 8: Combined solutions}
\psn
{\it The solutions of $X^2\sspm 2\,\widehat Y^2\sspeq -a$ and $t^2\sspp 2\,k^2\sspeq a$, for $ a\sspgeq 2$, qualifying for the final solutions come in three versions. The trivial case, if both congruences are solved by improper solutions (obtainable from  proper ones) is not considered because the curvature triples will not be primitive.
\psn
{\bf Type a)} 
\psn
If both binary quadratic forms are solved properly then
\Beq
a\sspeq 2^{e_2}\,\prod_{i\sspeq 1}^{P_1} p_{1,i}^{e_{1,i}},
\Eeq
with $e_2\sspin \{0,\,1\}$, $e_{1,i} \sspgeq 0$ and $P_1\sspgeq 0$, but not all exponents vanish. The products are over powers of distinct odd primes with $p_{1,i}$ congruent to \modan{1}{8}, \seqnum{A007519}. For the sequence $a$ see \seqnum{A192453} without member $1$.\pn
It is guaranteed  (by {\sl Corollary 5 b)}) that there is exactly one solution satisfying the necessary condition $0\sspkl X\sspkl \widehat Y$, later to be tested with $Y\sspeq \widehat Y\sspm k$ for $X\sspkl Y$.
\psn
The first instance for a final solution appears for $a\sspeq 89$. See {\sl Table 4} for $c_3\sspeq 9$.
\psn
{\bf Type b)}
\psn  
If the solutions $(X,\widehat Y)$ are improper and those of $(t,\,k)$ are proper, then
\Beqarray
a &\sspeq& A^2\,\tilde a,\nonumber \\
A &\sspeq& \prod_{i\sspeq 1}^{P_1^{\prime}} p_{1,i}^{e_{1,i}}\,\prod_{i\sspeq 1}^{P_3} p_{3,i}^{e_{3,i}},\, \text{and} \nonumber \\
\tilde a &\sspeq& 2^{e_2}\,\prod_{i\sspeq 1}^{P_1^{\prime\prime}} p_{1,i}^{e_{1,i}},
\Eeqarray
where $P_1^{\prime} \sspp P_1^{\prime\prime}\sspeq P_1$, the number of distinct prime factors of $a$ congruent to \modan{1}{8}, \seqnum{A007519}, and $P_3$ is the number of distinct primes congruent to \modan{3}{8}, \seqnum{A033200}. Also $P_1^{\prime}\sspp P_3 \sspgr 0$ in order that $A\sspgeq 2$, \ie the exponents $e_{1,i}$ and $e_{3,i}$ are nonnegative and do not both vanish. $e_2\sspin {0,\,1}$, and $\tilde a \sspgr 2$.
\pn 
$X\sspeq A\,\tilde X$, $\widehat Y\sspeq A\,\tilde {\widehat Y}$, and there is exactly one proper solution for $(\tilde X,\, \tilde{\widehat Y})$ representing $-\tilde a$ satisfying the necessary condition $0\sspkl\tilde X\sspkl \tilde{\widehat Y}$ for each pair of conjugate families of solutions (therefore $\tilde a\sspeq 2$ is excluded). Later $Y\sspeq \widehat Y\sspm k$ has to be tested for $X\sspkl Y$ for each unique solution from each family pair, with each of the $2^{P_1\sspp P_3\sspm 1}$ solutions of $k$ from the $(t,\, k)$ equation for positive $t$ and $k$. Finally the curvature triples have to be tested for primitivity.\pn
For the sequence $A$ see \seqnum{A225771} without $1$. For the sequence $\tilde a$ see \seqnum{A192453} without $1$ and $2$.
\psn
The first instance for a final solution appears for $a\sspeq 153\sspeq 3^2\cdot 17$. See {\sl Table 4} for $c_3\sspeq 12$.
\psn
\vfill\eject\noindent
{\bf Type c)} 
\psn
If the solutions $(X,\widehat Y)$ are proper and those of $(t,\,k)$ are improper, then
\Beqarray
a &\sspeq& A^2\,\tilde a,\ \text{with} \nonumber\\
A &\sspeq& \prod_{i\sspeq 1}^{P_1^{\prime}} p_{1,i}^{e_{1,i}}\, \prod_{i\sspeq 1}^{P_7} p_{7,i}^{e_{7,i}},\ \text{and} \nonumber\\
\tilde a &\sspeq& 2^{e_2}\,\prod_{i\sspeq 1}^{P_1^{\prime\prime}} p_{1,i}^{e_{1,i}},
\Eeqarray
where $e_2\sspin \{0,\,1\}$, $P_1^{\prime} \sspp P_1^{\prime\prime}\sspeq P_1$, the number of distinct prime factors of $a$ congruent to \modan{1}{8}, \seqnum{A007519}, and $P_7$ the number of distinct prime factors of $a$ congruent to \modan{7}{8}, \seqnum{A007522}. For the primes congruent to $\{1, 7\} (mod\,8)$ see \seqnum{A001132}, or \seqnum{A038873} with the member $2$.  Also $P_1^{\prime}\sspp P_7 \sspgr 0$ in order that $A\sspgeq 2$; $e_2\sspeq 1$ with $P_1^{\prime\prime}\sspeq 0$ is not allowed because the forbidden $t\sspeq 0$ would follow. The case $e_2\sspeq 0$ and $P_1^{\prime\prime}\sspeq 0$ is also not allowed because this leads to $k\sspeq 0$.
\psn
For the sequence $A$ see \seqnum{A058529} without member $1$. For the sequence $\tilde a$ see \seqnum{A192453} without $1$ and $2$.
\psn
$t\sspeq A\,\tilde t$, $k\sspeq A\,\tilde k$, with the $2^{P_1^{\prime\prime}\sspm 1}$ proper solutions $(\tilde t,\,\tilde k)$ representing $\tilde a$, with positive $t$ and $k$. Then $Y\sspeq \widehat Y \sspm k$ has to be tested for $X\sspkl Y$ for each of the $2^{P_1\sspp P_7 \sspm 1}$ solutions with $Y\sspgr 0$ and $0\sspkl X\sspkl \widehat Y$ for each of the $2^{P_1^{\prime\prime}\sspm 1}$ solutions for $\tilde k$. Finally the primitivity of all $(c_1,\,c_2,\,c_3)$ has to be checked.
\psn
The first instance for a final solution appears for $a\sspeq 7^2\cdot 2\cdot 17 sspeq 1666$ (larger than the values shown in {\sl Table 4}).
with $[X,\,\widehat Y]\sspeq [4,\,29]$, $[t,\,k]\sspeq [28,\,21]$, $[X,\,\widehat Y]\sspeq [4,\,8]$,  $[c_1,\,c_2,\,c_3]\sspeq [2,\, 6,\, 39])$ and $[c_{4,-},\,c_{4,-}]\sspeq [11,\,83]$.
}
\psn
Before giving the proof several examples will illustrate this {\sl Proposition}.
\psn
{\bf Example 13}
\psn
Necessary conditions for the final solution are $0\sspkl X\sspkl Y$ and $c_2\sspkl c_3$, \ie $X\sspkl t\sspp 2\,k$.
\psn
{\bf Type a)} \pn
{\bf i)}: $a\sspeq 2\cdot 17^2\sspeq 578$, $(X,\,\widehat Y)\sspeq (12,\,19)$, $(t,\,k)\sspeq (24,\,1)$. $Y\sspeq \widehat Y\sspm k\sspeq 18$. Hence $0\sspkl X\sspkl Y$. $c_3\sspeq (\widehat Y \sspp t\sspp k)/2\sspeq 44/2\sspeq 22$. $q\sspeq (\widehat Y \sspp t\sspm k)/2\sspeq 21$. $c_1\sspeq x\sspeq (Y\sspm X)/2\sspeq 3$ and $c_2\sspeq x\sspeq (Y\sspp X)/2\sspeq 15$. Thus, there is the primitive $DS$-triple $(3,\,15,\,22)$, with $c_{4,\pm}\sspeq c_1\sspp c_2\sspp c_3 \,\pm\, 2\,q$ given by $c_{4,-}\sspeq -2$ and $c_{4,+}\sspeq 82$. See {\it Table 4}. 
\psn
{\bf ii)}: $a\sspeq 17\cdot 41\sspeq 697$. Two solutions $(X,\,\widehat Y)_1\sspeq (5,\, 19)$ and $(X,\,\widehat Y)_2\sspeq (19,\, 23)$. Two solutions  $(t,\,k)_1 \sspeq (25,\,6)$ and  $(t,\,k)_2 \sspeq (7,\,18)$. Four possibilities: $Y_{1,1}\sspeq  \widehat Y_1 \sspm k_1\sspeq 13$, $Y_{1,2}\sspeq \widehat Y_1\sspm k_2\sspeq 1$, and $Y_{2,1}\sspeq \widehat Y_2\sspm k_1\sspeq 17$, $Y_{2,2} \sspeq \widehat Y_2 \sspm k_2\sspeq 5$. Only $(X_1,\, Y_{1,1})\sspeq (5,\,13)$ qualifies with $X\sspkl Y$. $c_3\sspeq (19\sspp 25\sspp 6)/2\sspeq 25$, $q\sspeq 19$, $c_1\sspeq x\sspeq (13\sspm 5)/2\sspeq 4$ and $c_2\sspeq y\sspeq (13\sspp 5)/2 \sspeq 9$. The {\sl DS}-triple $(4,\,9,\,25)$ is primitiv. $c_{4,-}\sspeq 0$ and $c_{4,+}\sspeq 76$. This is a degenerate case. See {\sl Table 4} and {\sl Table 6} under $n\sspeq 5$.
 \psn
{\bf Type b)} \pn
{\bf i)} $a\sspeq (3^2\cdot 17^2)\,17\sspeq 44217$. Two pairs of conjugate families. $A\sspeq 3\cdot 17$, $\tilde a\sspeq 17$. $X\sspeq A\,\tilde X$, $\widehat Y\sspeq A\,\widetilde{\widehat Y}$. Positive proper solutions representing $3^2\,17^3$: $(t,\,k)_1\sspeq (107,\,128)$ and $(t,\,k)_2\sspeq (197,\,52)$. Proper positive solution representing $-17$: $(\tilde X,\,\widetilde{\widehat Y})\sspeq (1,\,3)$, $(X,\,\widehat Y)\sspeq (51,\,153)$. $Y_1\sspeq \widehat Y\sspm k_1\sspeq 153\sspm 128\sspeq 25$, $Y_2\sspeq \widehat Y\sspm k_2\sspeq 153\sspm 52\sspeq 101$, $(X,\,Y_1)\sspeq (51,\,25)$ not qualifying, and the solution is $(X,\,Y_2)\sspeq (51,101)$. $c_3\sspeq (\widehat Y\sspp t_2\sspp k_2)/2\sspeq 201$ and $q\sspeq 149 $. $x\sspeq (101\sspm 51)/2\sspeq 25$ and $y\sspeq (101\sspp 51)/2\sspeq 76$, The final primitive $DS$-triple is $(25,\,76,\,201)$ with $[c_{4,-},\,c_{4,+}]\sspeq [4, 600]$.
\psn 
{\bf ii)} $a\sspeq 3^2\cdot 17^3\sspeq 44217$. Two pairs of conjugate families. $A\sspeq 3$, $\tilde a\sspeq 17^3$. $X\sspeq A\,\tilde X$, $\widehat Y\sspeq A\,\widetilde{\widehat Y}$. Proper solutions representing $3^2\cdot 17^3$ like above: $(t,\,k)_1\sspeq (107,\,128)$ and $(t,\,k)_2\sspeq (197,\,52)$. Proper solution representing $-17^3$: $(\tilde X,\,\widetilde{\widehat Y})\sspeq (55,\,63)$, $(X,\,\widehat Y)\sspeq (165,\,189)$. $Y_1\sspeq \widehat Y\sspm k_1\sspeq 189\sspm 128\sspeq 61$, $Y_2\sspeq \widehat Y\sspm k_2\sspeq 189\sspm 52\sspeq 137$, $(X,\,Y_1)\sspeq (165,\,61)$, not qualifying, and $(X,\,Y_2)\sspeq(165,\,137)$, also not qualifying.\pn
Conclusion: different splittings of $a$ with improper solutions for $(X,\,\widehat Y)$ lead to different answers. Here only the final solution of part {\bf i)} is possible.\psn 
{\bf iii)} A similar situation appears for $a\sspeq 3^2\cdot 17^2\cdot 41\sspeq 106641$. The splitting $a\sspeq A^2\,\tilde a$ with $A\sspeq 3\cdot 17$ and $\tilde a\sspeq 41$ does not lead to a final solution, whereas $A\sspeq 3$ and $\tilde a\sspeq 17^2\cdot 41\sspeq 11849$ leads to only two (not eight) final solutions with $(X,\,Y)_1\sspeq (9,\,85)$ and $(X,\,Y)_2\sspeq (9,\,65)$. $c_{3,1}\sspeq 315$ and $c_{3,2}\sspeq 312$.
$q_1\sspeq 169$ and $q_2\sspeq 146$. $[c_1,\,c_2,\,c_3]_{1}\sspeq [38,\, 47,\,315]$ and $(c_1,\,c_2,\,c_3)_{2}\sspeq (28,\, 37,\,312)$. $[c_{4,-},\,c_{4,+}]_1\sspeq [62,\,738]$ and $[c_{4,-},\,c_{4,+}]_2\sspeq [85,\, 669]$.
\pn
A splitting with $A\sspeq 17$ and $\tilde a\sspeq 3^2\cdot 17$ is not possible because $\tilde a \notequiv \modan{1}{8}$.\pn
In all these examples neither type {\bf a)} nor type {\bf c)} is possible because there is no proper solution for $(X,\,\widehat Y)$ due to the factor $3$ of $a$.
\psn
{\bf Type c)} \pn
{\bf i)} $a\sspeq 17^2\,(2\cdot 17)\sspeq 9826$. $A\sspeq 17$, $\tilde a\sspeq 2\cdot 17$.
Positive proper solution representing $-2\cdot 17^3$: $(X,\,\widehat Y)\sspeq (16,\,71)$. Positive proper solution representing $2\cdot 17\sspeq 34$: $(\tilde t,\,\tilde k)\sspeq (4,\,3)$. $(t,\,k)\sspeq (68,\,51)$. $Y\sspeq \widehat Y\sspm k\sspeq 20$. $(X,\,Y)\sspeq (16,\,20)$. $c_3\sspeq (71\sspp 68\sspp 51)/2\sspeq 95$, $q\sspeq 44$. Final $DS$-triple $(c_1,\,c_2,\,c_3)\sspeq (2,\,18,\,95)$, and  $[c_{4.-},\,c_{4,+}]\sspeq[27,\,203]$.\pn
No other splitting of $a$ is possible.
\pn
Note that there is no final solution of type {\bf a)} for $a\sspeq 2\cdot 17^3$ because $(t,\,k)\sspeq (76,\, 45)$, hence $(X,\,Y)\sspeq (16,\, 71\sspm 45)\sspeq 2\cdot(8,\,13)$ is improper. 
\pn
There is also no final solution of type {\bf b)} with splitting $A\sspeq 17$ and $\tilde a\sspeq 2\,\cdot 17$ because with $(X,\, \widehat Y)\sspeq 17\cdot(4,\,5)\sspeq (66,\,85)$ and the just given $(t,\,k)\sspeq (76,\, 45)$ one obtains the improper solution $(X,\,Y)\sspeq 2\cdot(33,\,20)$. Thus, for $a\sspeq 2\cdot 17^3$ there is only the one type {\bf c)} solution.  
\psn
{\bf ii, a)}  $a\sspeq (7\cdot 17)^2\,17\sspeq 240737$. $A\sspeq 7\cdot 17$, $\tilde a\sspeq 17$. Two positive proper solution representing $-240737$: $(X,\,\widehat Y)_1\sspeq (9,\,347)$ and $(X,\,\widehat Y)_2\sspeq (115,\,363)$. Positive proper solution representing $17$: $(\tilde t,\,\tilde k)\sspeq (3,\,2)$, $(t,\,k)\sspeq (357,\,238)$. $Y_1\sspeq \widehat Y_1\sspm k \sspeq 109$, and $Y_2\sspeq \widehat Y_2\sspm k \sspeq 98$ is not qualifying. The solution $(X,\,Y)\sspeq (9,\,109)$ leads to $c_3\sspeq (347\sspp 357\sspp 238)/2 \sspeq 471$ and $q\sspeq 233$. The final primitive $DS$-triple is $(c_1,\,c_2,\,c_3)\sspeq (50,\,59,\,471)$, with $[c_{4.-},\,c_{4,+}]\sspeq [114,\,1046]$.
\psn 
{\bf ii, b)}  $a\sspeq 7^2\,17^3\sspeq 34391$. $A\sspeq 7$, $\tilde a\sspeq 17^3$. Like in the previous example the two proper solutions $(X,\,\widehat Y)_1\sspeq (9,\,347)$ and $(X,\,\widehat Y)_2\sspeq (115,\,363)$. Positive proper solution representing $17^3\sspeq 4913$: $(\tilde t,\,\tilde k)\sspeq (45,\,38)$. $(t,\,k)\sspeq (315,\,266)$. $Y_1\sspeq \widehat Y_1 \sspm k \sspeq 81$, and $Y_2\sspeq \widehat Y_2 \sspm k \sspeq 97\sspkl 115$ does not qualify. The solution $(X,\,Y)\sspeq (9,\,81)$ leads to $c_3\sspeq (347\sspp 315\sspp 266)/2\sspeq 464$ and $q\sspeq 198$. The only final primitive $DS$-triple is $(c_1,\,c_2,\,c_3)\sspeq (36,\,45,\,464)$, with $[c_{4.-},\,c_{4,+}]\sspeq [149,\,941]$.
\pn
In these examples $\bf ii)$ neither type {\bf a)} nor type {\bf b)} applies because of the factor $7$ of $a$.
\pbn
{\bf Proof of Proposition 8}:
\psn
{\bf a)} This follows from the intersection of {\sl Lemma 18, a)} \Eq{55} and 
{\sl Lemma 19, a)} \Eq{59}.\psn
{\bf b)} The proper solution $(t,\, k)$ of $t^2\sspp 2\,k^2\sspeq A^2\,\tilde a$ allows distinct primes $2,\, p_{1,i},\, p_{3,i}$ by {\sl Lemma 19} a). The improper solution $(X,\,\widehat Y)\sspeq (A\,\tilde X,\,A\,\widetilde{\widehat Y})$, with $\gcd(\tilde X,\,\widetilde{\widehat Y})\sspeq 1$, of $X^2\sspm 2\,\widehat Y^2 \sspeq -A^2\,\tilde a$, allows for $(\widetilde X,\, \widetilde{\widehat Y})$ only distinct primes $2,\, p_{1,i}$ because the previous restriction for $(t,\, k)$ implies $P_7\sspeq 0$. $A$ allows distinct $p_{1,i}$ and $p_{3,i}$ but not $2$, because $e_2\sspin \{0,\,1\}$.
\pn
For the fact that for $a\sspgeq 2$ pairs of conjugate families, of types $I$ and $II$, appear for proper solution the of $(X,\,\widehat Y)$ see {\sl Corollary 5 b)}. This also proves the existence of a unique positive solution $0 \sspkl X^{\prime}\sspkl \widetilde Y^{\prime}$ in each family $I$. Therefore the number of solutions is $2^{P_1\sspp P_7\sspm 1}$. The  proper $(t,\, k)$ solutions come for $a\sspgr 2$ in a pair of conjugate solutions with $k\sspgr 0$, and there is one solution with also $t\sspgr 0$ ($t\sspeq 0$ is forbidden by {\sl Lemma 17}). Thus, the number of positive proper solutions $(t,\, k)$ is $2^{P_1\sspp P_3\sspm 1}$.
\psn
{\bf c)} Only primes $2,\, p_{1,i}$, and $p_{i,7}$ can appear because of the proper solutions of $(X,\,\widehat Y)$ by {\sl Lemma 18}. Hence the improper solutions $(t,\, k)\sspeq(A\,\tilde t,\,A\, \tilde k)$ can have $A$ with primes of type $p_{1,i}$ and $p_{7,i}$, and $\tilde a$ can only have primes of type $2$ and $p_{1,i}$ for the solutions of $(\tilde t,\,\tilde k)$. The splitting of the primes $p_{1,i}$ in $A$ and $\tilde a$ is  $P_1^{\prime}\sspp P_1^{\prime\prime}\sspeq P_1$.
\pn
The case $\tilde a\sspeq 2$ ($P_1^{\prime\prime}$\sspeq 0) is forbidden because the $(\tilde t,\, \tilde k)$ solution with $k\sspgr 0$ would be $(0,\, 1)$, with forbidden $t\sspeq 0$. Also $\tilde a\sspeq 1$ is excluded because the solutions $(\tilde t,\, \tilde k)\sspeq (\pm 1,\,0)$ lead to $k\sspeq 0$. \hskip 1.1cm $\square$ 
\pbn
{\bf Case iii)} {$\bf q\sspgr c_3$}
\psn
The case $0\sspkl q\sspeq c_3\sspp k$, with $k\sspgeq 1$, is treated by solving \Eq{50} of the previous case {\it ii)} after substituting $k$ by $-k$. The quadratic equations \Eq{53} are the same but the variable $\widehat Y$ is now $Y\sspm k$, thus $0\sspkl X\sspkl Y\sspeq \widehat Y\sspp k$ has to be satisfied. The inequalities for the curvatures are 
\Beq
c_1\sspeq \frac{\widehat Y \sspp k\sspm X}{2} \sspkl c_2\sspeq \frac{\widehat Y \sspp k\sspp X}{2} \sspkl c_3\sspeq q\sspm k\,.
\Eeq
\psn 
Here two cases follow from \Eq{51} with $k\sspto -k$, and 2\,x \sspeq Y\sspm X and 2\,y\sspeq Y\sspp X,  namely:
\psn
Case a) $ 0\sspkl 2\,q\sspeq Y\sspp t$, with $t\sspgeq 0$ and $t\sspdef \sqrt{2\,Y^2 \sspm X^2\sspm 4\,Y\,k} \sspeq \sqrt{2\,\widehat Y^2 \sspm X^2\sspm 2\,k^2}$. 
\psn
Case b) $ 0\sspkl 2\,q\sspeq Y\sspm t$, with $t\sspgeq 0$ defined in case a).
\psn
{\bf Lemma 20: Upper bound for $\bf k$}
\psn
{$1\sspleq k\sspleq \floor{\sqrt{\frac{a}{2}}}$, with $a\sspeq 2\,Y^2 \sspm X^2 \sspeq t^2 \sspp 2\,k^2$.}
\psn
{Proof:} In order that $t$ is real one needs $a\sspeq 2\,\widehat Y^2\sspm X^2 \sspgeq 2\,k^2$. Thus, $ k\sspleq \sqrt{\frac{a}{2}}$. \hskip 4cm $\square$
\psn
The next Lemma shows that for both cases the instance $t\sspeq 0$ does not qualify.
\psn
{\bf Lemma 21:} {\it There is no solution for $(X,\,\widehat Y)$ and $c_2\sspkl c_3$ if $t \sspeq 0$}.
\psn
{\bf Proof:} If $t\sspeq 0$ then $0\sspkl X^2\sspeq 2\,\widehat Y^2 \sspm 2\,k^2\sspeq 2\,(\widehat Y\sspm k)\,(\widehat Y \sspp k)$. In both cases $2\,q\sspeq \widehat Y \sspp k$, hence\pn
$ 0\sspkl \widehat Y\sspm k \sspeq 2\,(q\sspm k) \sspeq 2\,c_3$. One shows that, in fact, $c_2 \sspgr c_3$.  $2\,c_2 \sspeq \widehat Y\sspp k\sspp X \sspgr 2\,c_3\sspeq \widehat Y \sspm k$, meaning that $2\,k \sspp X\sspgr 0$, which is true because $X$ and $k$ are positive. \hskip 7cm $\square$
\psn
From $t^2\sspgr 0$ the upper bound becomes $k\sspkl \sqrt{\frac{a\sspm 1}{2}}$.\psn
\psn
{\bf Lemma 22:} {\it Case b) does not qualify for a final solution because $c_2\sspgr c_3$}. 
\psn
{Proof:} For this case $2\,q\sspeq \widehat Y\sspp k \sspm t$, with $t\sspgr 0$. $2\, c_3\sspeq  2\,(q\sspm k) \sspeq \widehat Y \sspm k\sspm t$. Hence $2\,c_2\sspeq \widehat Y\sspp k \sspp X \sspgr 2\,c_3$ because $2\,k \sspp X \sspp t\sspgr 0$, which is satisfied because each term is positive. \hskip 16.3cm $\square$ 
\psn
Therefore only case a) has to be treated.
\psn
Concerning the representable values $a\sspgeq 2$ of the two quadratic forms (the case $a\sspeq 1$ is not relevant for nonvanishing $k$; see also the remark in {\it Lemma 18 a)}) the results from {\sl Proposition 8} for case ii) apply. For $X^2\sspm 2\,\widehat Y^2\sspeq -a$ there is {\it Lemma 18}, \Eq{55}, and for {$t^2\sspp 2\,k^2\sspeq a$ see {\it Lemma 19}, \Eq{59}. The combined solutions come also in the three cases $a),\,b), $ and $c)$ of {\sl Proposition 8}, eqs. (60), (61), and (62), respectively. \pn
The uniqueness of the solution $(X^{\prime},\, \widehat Y^{\prime})$ with $a\sspgeq 2$ and $0\sspkl X^{\prime}\sspkl \widehat Y^{\prime} \sspkl \sqrt{a}$ for one of the members of family I (with a certain index $i\sspeq i^{\prime}$) of a conjugate pairs of solutions is also guaranteed, like in {\it Lemma 11 b)}, \Eq {46}, but with a change of variables $s\sspto \sqrt{a}$, $Y\sspto \widehat Y$, $\widehat Y\sspto \widehat Y^{\prime}$ and $\widehat X\sspto X^{\prime}$. 
\psn
For final solutions one has to check whether $c_2\sspkl c_3$ and the primitivity of the curvature triple holds for those solutions corresponding to the unique one in each family $I$ which satisfies $0\sspkl X^{\prime} \sspkl \widehat Y^{\prime}\sspkl Y$. 
\pn
But because now  $Y\sspeq \widehat Y\sspp k$, in contrast to case ii) with $q\sspeq c_3 \sspm k$ and  $Y\sspeq \widehat Y\sspm k$, further candidates may qualify.  It could happen in both families that also solutions with $\widehat Y_i\sspkl X_i$ and $X_i\sspgr 0$ (after a possible common sign flip) satisfy the necessary condition $Y_i\sspeq \widehat Y_i\sspp k\sspgr X_i $ provided $k$ is large enough. Fortunately, this will be proved not to happen, except for one member in each family $II$ for which a proof is missing. This member is $(X_{i^{\prime\prime}\sspp 1}), \widehat Y_{i^{\prime\prime}\sspp 1}$. For these instances the conjecture is that they cannot lead to a final solution because the additional condition  $c_3\sspgr c_2$ fails. This is the content of the next {\sl Lemma}.\psn    
\psn
{\bf Lemma 23: Excluding certain candidates for final solutions.} 
\psn
{\it There is no qualifying solution $(X,\,\widehat Y)$ with $X\sspgr 0$ and $X \sspgr \widehat Y$, with $Y\sspeq \widehat Y\sspp k\sspgr X$. This means that  $X\sspm \widehat Y\ \sspgeq k$, except for possibly one solution of each family $II$. \pn
If such a solution of family $II$, \viz\, $(X_i,\, \widehat Y_i)$, where $i\sspeq i^{\prime\prime} \sspp 1$, (\ie the first positive solution) satisfies $0\sspkl X_i\sspkl Y_i$ then the {\bf conjecture} is that it cannot be  a final solution because the necessary condition $c_3\sspgr c_2$ is not satisfied, \ie $X_i\sspgr 0$ and $X_i\sspgeq t\sspm 2\,k$,      
}\psn
{\bf Proof:}\pn
a) Family $I$ (see {\sl Definition 3} and {\sl Lemma 13} with the above mentioned adaption of $s,\,Y,\, \widehat Y$ and $\widehat X$)\pn
Case i) For $i\sspeq i^{\prime}$ there is no such solution because $0\sspkl X^{\prime}\sspkl \widehat Y^{\prime}$.
\psn
Case ii) With $i\sspgr i^{\prime}$, one has solutions $X_i\sspgr \widehat Y_i \sspgr 0$. An inductive proof, using the recurrence \Eq{40} and the upper bound $k\sspkl \sqrt{a/2}$, shows that a qualifying solutions with $0\sspkl X_i\sspkl Y_i\sspeq \widehat Y_i\sspp k$, \ie $X_i\sspm \widehat Y_i \sspkl  k$, does not exist. \pn 
Start with $i\sspeq i^{\prime}\sspp 1$, and the recurrence $X_{i^\prime\sspp 1} = 4\,X^\prime  \sspp 3\,\widehat Y^{\prime}$ and  $\widehat Y_{i^\prime \sspp 1} = 4\,X^\prime  \sspp 3\,\widehat Y^{\prime}$, with $0\sspkl X^{\prime}\sspkl \widehat Y^{\prime}$. Hence $X_i\sspm \widehat Y_i\sspeq X^\prime\sspp \widehat Y^\prime$.  But $X^\prime\sspp \widehat Y^\prime \sspkl k$ implies, with the upper bound of $k$, that this is also $\sspkl \sqrt{ {\widehat Y{^\prime}}^2 \sspm {X^\prime}^2/2}$. Squaring leads to $3\, X^{\prime}\sspp 4\,\widehat Y^{\prime}\sspkl 0$, because $X^{\prime}\sspgr 0$; a contradiction, because also $\widehat Y^{\prime}\sspgr 0$.\pn
The same contradiction appears for any solution $(X_i,\,\widehat Y_i)$ for $i \sspgeq i^{\prime}\sspp 2$ with $X_i\sspgr \widehat Y_i \sspgr 0$. Assuming that $ X_i\sspm \widehat Y_i \sspeq X_{i\sspm 1}\sspp \widehat Y_{i\sspm 1} \sspkl  k$. then the upper bound of $k$ yields $3\,X_{i\sspm 1} \sspp 4\,\widehat Y_{i\sspm 1}\sspkl 0$; a contradiction because both terms are positive.
\psn
Case iii) We use a common sign flip in $X_i$ and $\widehat Y_i$ (a proper equivalence transformation) for $i \sspkl i^{\prime}$ of the adapted {\sl Lemma 13 a) iii)} such that $Y_i\sspkl 0$ and $0\sspkl \abs{\widehat Y_i}\sspkl X_i$. In this case no recurrence is needed. The proof that $0\sspkl X_i \sspkl Y_i\sspeq \widehat Y_i\sspp k$ is not possible runs as follows. Assume that $X_i \sspm \widehat Y_i \sspeq X_i\sspp \abs{\widehat Y_i} \sspkl k$. But using $k\sspkl \sqrt{{\widehat Y_i}^2 \sspm X_i^2/2}$ leads, after squaring and using $X_i\sspgr 0$, to the contradiction $3\,X_i \sspp 4\,\abs{\widehat Y_i} \sspkl 0$. \pn
Note that this direct proof, not using the recurrence equation, does not work for the above case ii). It would lead for $i\sspgr i^{\prime}$ only to $3\,X_i\sspm 4\,\widehat Y_i$, which cannot be proved to be $\sspgeq 0$; it only follows that $\sspgr -X_i$ or $\sspgr -\widehat Y_i$. 
\psn
b) Family $II$ \psn
Case i) After adaption of {\sl Lemma 13 b) i)} with $i\sspeq i^{\prime\prime}$ and a common sign flip $0\sspkl X_i \sspkl \abs{Y_i}$. Assume $0 \sspkl X_i\sspkl \widehat Y_i\sspeq \widehat Y_i$, \ie $X_i\sspm \widehat Y_i\sspeq X_i\sspp \abs{\widehat Y_i}\sspkl k$. Then  
the upper bound on $k$ leads, after squaring and division by $X_i\sspgr 0$, to the contradiction $X_i \sspp \abs{\widehat Y_i}\sspkl 0$.
\psn
Case ii) for $i\sspgr i^{\prime\prime}$ has $0\sspkl \widehat Y_i\sspkl X_i$. One has to use the recurrence \Eq{40} like above in case a) ii). Assume $0\sspkl X_i\sspkl Y_i\sspeq \widehat Y_i\sspp k$, then by the recurrence $X_i\sspm \widehat Y_i \sspeq  X_{i-1}\sspp \widehat Y_{i-1}\sspkl k$. For all $i\sspgeq i^{\prime\prime} + 2$ the upper bound on $k$ leads, after squaring and dividing by $X_{i-1}$, to $3\,X_{i-1}\sspp 4\, \widehat Y_{i-1}\sspkl 0$, a contradiction because both terms are positive, because $i-1 \sspgr i^{\prime\prime}$.
\pn
However, here the mentioned gap in the proof appears for $i\sspeq i^{\prime\prime}\sspp 1$. \pn
With $0\sspkl \widehat Y_i\sspkl X_i$ assume that $X_i\sspm \widehat Y_i \sspkl k$. With Recurrence \Eq{40} this becomes $X_{i^{\prime\prime}}\sspp \widehat Y_{i^{\prime\prime}}\sspeq -X^{\prime}\sspp \widehat Y^{\prime}\sspkl k$. Using the upper bound for $k$ like above will not lead to a contradiction because it leads to $3\,X^{\prime} \sspm 4\widehat Y^{\prime}\sspkl 0$. \Eg $a\sspeq 89$. In family $I$ $(X^\prime,\,\widehat Y^\prime)\sspeq (3,\,7)$ and $(t, k)\sspeq (9,\,2)$. The upper bound for $k$ is $6$, hence $3\cdot 3\sspm 7\sspeq 4\sspkl 6$, and the assumption seems to be valid but with the proper $k$ value it leads to a contradiction. The used upper bound is too large.
\psn
{\bf Conjecture:}
\pn
{\it For the proper solution $(X_i,\, \widehat Y_i)$ with $i\sspeq i^{\prime\prime}\sspp 1$ in family $II$ one checks first if the proper positive $(t,\,k)$ solution of $t^2 \sspm 2\,k^2 \sspeq a$ satisfies the $c_3\sspgr c_2$ condition  $ X_i\sspkl t\sspm 2\,k$. If it fails no final solution is possible, even if $0\sspkl X_i\sspkl Y_i$ may hold. Otherwise one has $X_i\sspm \widehat Y_i \sspgeq k$. Using the recurrence this means that $-X^{\prime} \sspp \widehat Y^{\prime} \sspgeq k$, where $(X^{\prime},\, \widehat Y^{\prime})$ is the unique positive fundamental solution of $X^2\sspm 2\,\widehat Y^2\sspeq -a$ in the corresponding family $I$, satisfying $0\sspkl X^{\prime}\sspkl \widehat Y^{\prime}$. 
\pn
Using the correct $k$ value the conjecture means that $X^{\prime}\,(3\,X^{\prime} \sspm 4\, \widehat Y^{\prime}) \sspp t^2\sspgeq 0$, where $(t,\, k)$ is the unique solution of $t^2\sspp 2\,k^2\sspeq a$ with positive $k$ and $t$, where the $a$ value is from \seqnum{A192453} without members $1$ and $2$. }
\psn
Case iii) For $i\sspkl i^{\prime\prime}$ we use,as above in case a), a common sign flip, to get   $Y_i\sspkl 0$ and $0\sspkl \abs{\widehat Y_i}\sspkl X_i$. Like in case a) above assume that $0\sspkl X_i\sspkl Y_i\sspeq \widehat Y_i\sspp k$, \ie $X_i\sspp \abs{\widehat Y_i} \sspkl k$.
Using the upper bound of $k$, squaring and dividing by $X_i\sspgr 0$, leads to the contradiction $3\,X_i\sspp 4\,\abs{Y_i} \sspkl 0$. \hskip 4cm $\square$
\pbn
Thus, {\sl Proposition 8} can be taken over as {\bf Proposition 9}, just replacing there $k\sspto -k$ in the remarks at the end of each of the three cases. \psn
Examples for the three types follow.
\psn
{\bf Example 14:}
\psn
Necessary conditions for the final solution are $0\sspkl X\sspkl Y$ and $c_2\sspkl c_3$, \ie $X\sspkl t\sspm 2\,k$.\pn 
In the following examples the prime superscript on $X$ and $\widehat Y$ for the unique family $I$ solution with $O\sspkl X^{\prime}\sspkl \widehat Y^{\prime}$ will be omitted. 
\psn
{\bf Type a) Both binary quadratic forms are solved properly.} 
\psn
{\bf i)} $a\sspeq 2\cdot17^2\sspeq 578\sspeq$\seqnum{A192453}$(30)$. $(X,\,\widehat Y)\sspeq (12,\,19)$, $(t,\,k)\sspeq (24,\,1)$ (see {\it Example 11, a)i)}). $Y\sspeq \widehat Y\sspp k\sspeq 20$. Hence $0\sspkl X\sspkl Y$. $c_3\sspeq (\widehat Y \sspp t\sspm k)/2\sspeq 21$. $q\sspeq (\widehat Y \sspp t\sspp k)/2\sspeq 22$. $c_1\sspeq x\sspeq (Y\sspm X)/2\sspeq 4$ and $c_2\sspeq x\sspeq (Y\sspp X)/2\sspeq 16$. Thus, there is the primitive $DS$-triple $[4,\,16,\,21]$ with $c_{4,\pm}\sspeq c_1\sspp c_2\sspp c_3 \,\pm\, 2\,q$, hence $c_{4,-}\sspeq -3$ and $c_{4,+}\sspeq 85$. See {\it Table 5}.
\pn
The conjecture in {\sl Lemma 23} for the family $II$ element $(X_{i^{\prime\prime}\sspp 1},\, \widehat Y_{i^{\prime\prime}\sspp 1})\sspeq (40,\,33)$ is true because already the necessary condition $X_{i^{\prime\prime}\sspp 1}\sspkl t\sspm 2\,k$ fails, viz $40 \sspgr 24\sspm 1$. Indeed, $c_3\sspeq 23 \sspkl 26\sspeq c_2$. In this case also $40 \sspm 33\sspgeq 1$.
\psn
{\bf ii)}: $a\sspeq 17\cdot 41\sspeq 697\sspeq$\seqnum{A192453}$(47)$. Two solutions $(X,\,\widehat Y)_1\sspeq (5,\, 19)$ and $(X,\,\widehat Y)_2\sspeq (19,\, 23)$. Two solutions  $(t,\,k)_1 \sspeq (25,\,6)$ and  $(t,\,k)_2 \sspeq (7,\,18)$. The last solution has a negative $t\sspm 2\,k$ value, hence no positive $X$ solution qualifies. The other solution has $t\sspm 2\,k\speq 13$, hence the second $(X,\,\widehat Y)$ solution fails. Thus, only $(X,\, \widehat Y)\sspeq (5,\, 19)$ with $(X, \,Y)\sspeq (5,\,25)$ qualifies, and  $c_3\sspeq 19$, $2\,q \sspeq 50$, and  $[c_1,\,c_2]\sspeq [10,\,15]$. Thus, $[c_1,\,c_2,\,c_3]\sspeq [10,\,15,\, 19]$ qualifies, is primitive, and $[c_{4,-},\,c_{4,+}]\sspeq [-6,\, 94]$. This is the only final solution for $a\sspeq 697$, provided the conjecture is true. See {\sl Table 50}.
\pn
The conjecture in {\sl Lemma 23} for the family $II$ element is true because for $(X_{i^{\prime\prime} \sspp 1}, \, \widehat Y_{i^{\prime\prime} \sspp 1})\sspeq (61,\,47)$ the necessary condition $5\sspkl t\sspm 2\,k\speq 13$ fails. Indeed, $c_2\sspeq 57 \sspgr 33 \sspeq c_3$. In this case also $47\sspp 6 \sspkl 61$.
\psn
{\bf iii)}: $a\sspeq 17^3\sspeq 4913\sspeq$\seqnum{A192453}$(268)$.
\pn 
The positive proper solution $(t,\, k)$ representing $17^3$ is $(45,\,38)$, with $t\sspm 2\,k$ negative. Therefore the family $I$ qualifying solution $(X,\, \widehat Y)\sspeq (55,\,63)$, with $(X,\, Y)\sspeq (55,\,101)$ will not lead to a final solution.
\pn  
The conjecture for the first\ positive family $II$ element $(X_{i^{\prime\prime}\sspp  1},\,\widehat Y_{i^{\prime\prime}\sspp 1})\sspeq (87,\,79)$ is true because $t\sspm 2\,k$ is negative. This is an example where one has for $(X_i,\, \widehat Y_i)$, for $i\sspeq i^{\prime\prime}\sspp 1$, the promising looking condition $0\sspkl X_i\sspkl Y_i\sspeq \widehat Y_i\sspp k$, because $0\sspkl 87\sspkl 79\sspp 38$.
\psn
{\bf iv)}: $a\sspeq 17^2\cdot 41\sspeq 11849\sspeq$\seqnum{A192453}$(598)$.
\pn 
The two positive proper $(t,\,k)$ solutions representing $11849$ are $(93,\,40)$ and $(99,\,32)$ with $t\sspm 2\,k$ values $13$ and $35$. The two positive proper fundamental solutions $(X,\,\widehat Y)$ from family $I$ representing $-17^2\cdot 41$ are $(3,\,77)$ and $(101,\,105)$. Only the first solution survives the test $X\sspkl t\sspm 2\,k$ for both $(t,\, k)$ solutions. The two $(X,\, Y)$ solutions are then $(3,\, 117)$ and $(3,\,109)$. The corresponding $DS$-triples are $[57,\,60,\,65]$ and $[53,\, 56,\,72]$ with $[c_{4,-},\,c_{4,+}]$ values $[-28,\,392]$ and $[-27,\,389]$.  
\pn
The conjecture for the first positive family $II$ element $(X_i,\,\widehat Y_i)$, for $i\sspeq i^{\prime\prime}\sspp 1$, is true, because the two solutions are $(299,\,225)$ and $(117,\,113)$, and in both cases  $X_i\sspgeq t\sspm 2\,k$, for each of the two $t\sspm 2\,k$ values $13$ and $35$. In this case both  $(X_i,\,\widehat Y_i)$ solutions look promising for both $k$ values $40$ and $32$, but fail as final solutions.
\psn
\vfill\eject\noindent
{\bf Type b) $\bf (X,\,\widehat Y)$ improper and  $\bf (t,\,k)$ proper}
\psn
{\bf i)} $a\sspeq 3^2\cdot 73\sspeq 657$. $3\sspeq$\seqnum{A225771}$(2)$ and  $73\sspeq$\seqnum{A192453}$(6)$. There is one qualifying solution from family $I$ \viz\, $(X,\,\,\widehat Y)\sspeq 3\cdot(5,\, 7)\sspeq (15,\,21)$, and two positive solutions $(t,\,k)_1\sspeq (25,\,4)$ and $(t,\,k)_2\sspeq (23,\,8)$. The second solution fails because $X \sspgr t\sspm 2\,k$, \ie $15\sspgr 7$. The surviving solution has $Y\sspeq \widehat Y\sspp k_1 \sspeq 21\sspp 4\sspeq 25$, $c_3\sspeq (\widehat Y\sspp t_1\sspm k_1)/2\sspeq 21$ and $q\sspeq 25$. $c_1\sspeq (25\sspm 15)/2\sspeq 5$ and $c_2\sspeq (25\sspp 15)/2\sspeq 20$. Thus, the final primitive curvature solution is $(5,\,20,\,21)$ with $[c_{4,-},\,c_{4,+}]\sspeq [-4,\,96]$. See {\sl Table 5}. 
\pn
The first instance of type $b)$ appears for $a\sspeq 3^2\cdot 41\sspeq 369$ with the $DS$-triple $[4,\,13,\,16]$ with $[c_{4,-},\,c_{4,+}]\sspeq [-3,\,69]$. See {\sl Table 5}. 
\pn
The conjecture in {\sl Lemma 23} for the family $II$ proper representation for $-\tilde a\sspeq -73$, \viz \, $({\tilde X}_i,\, \tilde{\widehat Y}_i)\sspeq (13,\, 11)$, for $i\sspeq i^{\prime\prime} \sspp  1$, is true because $(t,\, k)\sspeq (1, \,6)$, hence $t\sspm 2\,k$ is negative. Here also $X_i \sspeq 3\cdot 13 \sspgr 17 \sspeq t_1\sspm 2\,k_1$, and also $\sspgr 7 \sspeq t_2\sspm 2\,k_2$, confirming the conjecture.  ${\tilde X}_i\sspm\tilde{\widehat Y}_i\sspeq 2\sspkl 6$ looks promising, but it fails as final solution representing  $-73$. In this example also $X_i\sspm \widehat Y_i \sspeq 39\sspm 33\sspeq 6\sspgr 4\sspeq k_1$ fails, but looks promising for $k_2\sspeq 8$. 
\psn
{\bf ii)} $a\sspeq (3^2\cdot 17^2)\,17\sspeq 44217$. $A\sspeq3\cdot 17\sspeq 51\sspeq $\seqnum{A225771}$(11)$ and $\tilde a \sspeq 17\sspeq$\seqnum{A192453}$(3)$. The two positive proper solutions representing $a\sspeq 3^2\,17^3$ are $(t,\,k)_1\sspeq (107,\,128)$ and $(t,\,k)_2\sspeq (197,\,52)$. The first solutions fails because $t\sspm 2\,k$ is negative. The positive fundamental proper representation of $-17$ from family $I$ is $(\tilde X,\,\widetilde{\widehat Y})\sspeq (1,\,3)$ This corresponds to the improper solution $(X,\,\widehat Y)\sspeq (51,\,153)$, leading to the proper solution $(X,\,Y)\sspeq (51,\, 153\sspp 52)\sspeq (51,\,205)$.
$c_3 \sspeq (153\sspp 197\sspm 52)/2\sspeq 149$ and $q\sspeq 201$. $c_1\sspeq 77$
 and $c_2\sspeq 128$. The primitive $DS$-trip;e is $[77,\,128,\,149]$, with $[c_{4,-},\,c_{4,+}]\sspeq [-48,\, 756]$.
\pn
The conjecture for the first positive family $II$ element representing $-17$, \viz\, $(\tilde X_i,\, \tilde{\widehat Y}_i)\sspeq (9,\,7)$, for $i\sspeq i^{\prime\prime}\sspp 1$, is true,because $(t,\,k) \sspeq (3,\,2)$ representing $17$ has negative value $t\sspm 2\,k$. Here with $\tilde X_i\sspm\tilde{\widehat Y}\sspeq 9\sspm 7\sspkl 2\sspeq k$ looks promising, but it fails as final solution The improper solution $(X_i,\,\widehat Y_i)\sspeq (459,\, 357)$ fails because $t_2\sspm 2\,k_2\sspeq 93 \sspkl X_i$. Also  $X_i\sspm \widehat Y_i\sspeq 102\sspgr 52\sspeq k_2$, hence  fails, but it looks promising for $k_1 \sspeq 104$, but fails as final solution.
\pn
Therefore the only final $DS$-triple for this decomposition of $a\sspeq 44217$ is the found one, \viz\, $[77,\,128,\,149]$.
\pn
It will be shown next that the other possible decomposition will not survive.
\psn
{\bf iii)} $a\sspeq 3^2\cdot 17^3\sspeq 44217$. $A\sspeq 3\sspeq$\seqnum{A225771}$(1)$, $\tilde a\sspeq 17^3\sspeq 4913\sspeq $\seqnum{A192453}$(269)$. From the two positive proper solutions $(t,\, k)$ given in the previous section {\bf ii)}, again only $(197,\,52)$ survives. The unique positive proper solution representing $-17^3$, have been given above in type a)iii) $(\tilde X,\,\tilde{\widehat Y})\sspeq (55,\,63)$. Now $(X,\,Y)\sspeq (3\cdot 55,\,3\cdot 63\sspp 52)\sspeq (165,\,241)$ looks promising, but $165\sspgr 83\sspeq t\sspm 2\,k$, thus fails as final solution.
\pn 
The conjecture for the first positive family $II$ solution $(\tilde X_i,\, \tilde{\widehat Y}_i)\sspeq (87,\,79)$ for $-\tilde a\sspeq -17^3$, from above type {\bf a), iii)}, is true, because  $t\sspm 2\,k\sspeq 45\sspm 2\cdot 38$ is negative. Here $\tilde X_i\sspm \tilde{\widehat Y}_i\sspeq 8 \sspkl 38$ looks promising. The improper solution $(X_i,\,\widehat Y_i)\sspeq (261,\, 237)$ fails also because $X_i\sspgr 93\sspeq t_2\sspm 2\,k_2$. Here $X_i\sspm \widehat Y_i\sspeq 24$ looks promising for $k_2\sspeq 52$ as well as for $k_1\sspeq 128$. Therefore the conjecture is true for this improper as well as proper case.
\pn
Thus, with this splitting of $a$ no final solution exists.
\psn
Conclusion: different splittings of $a$ with improper solutions for $(X,\,\widehat Y)$ may lead to different answers. For $a\sspeq 44217 $ only the final solution of part {\bf ii)} is possible.
\psn
{\bf iv)} A similar situation appears for $a\sspeq 3^2\cdot 17^2\cdot 41\sspeq 106641$.
\psn 
$\boldsymbol \alpha)$ The splitting $a\sspeq A^2\,\tilde a$ with $A\sspeq 3\cdot 17\sspeq$\seqnum{A225771}$(11)$ and $\tilde a\sspeq 41\sspeq$\seqnum{A192453}$(5)$ does not lead to a final solution, because all four positive proper solutions $(t,\, k)$ of $a$ have negative $t\sspm 2\,k$ values. These solutions are $(29,\,230)$, $(67, 226)$, $(227,\, 166)$ and $(253,\,146)$. Therefore no final solution exists. 
\pn
The only qualifying improper positive fundamental solutions from family $I$ and for this splitting is $(X,\, \widehat Y)\sspeq (153,\, 255)$, but it will not lead to a final result.
\pn
The conjecture for the first positive family $II$ proper solution  representing $-41$, \viz\, $(\tilde X_i,\,\tilde {\widehat Y_i})\sspeq (11,\,9)$,  where $i\sspeq i^{\prime\prime}\sspp 1$, representing $-41$, is true because $(t,\,k)\sspeq (3,\,4)$. with a negative $t\sspm 2\,k$ value.Here $(X_i,\,Y_i)\sspeq (11,\, 9\sspp 4)$ looks promising but does not lead to a final solution for $-41$ because $c_2\sspgr c_3$. \pn
Also the improper solution $X_i,\,\widehat Y_i\sspeq (561,\,459)$ representing $-106641$ does not lead to a final solution, because all four $t\spm 2\,k$ values are negative,  even though it looks promising for all four $k$ values given above. 
\psn
$\boldsymbol \beta)$ The splitting with $A\sspeq 3\sspeq$\seqnum{A225771}$(2)$ and $\tilde a\sspeq 17^2\cdot 41\sspeq 11849\sspeq$\seqnum{A192453}$(598)$ has the four positive proper $(t,\,k) $ solutions given in part $\boldsymbol \alpha)$, and $t\sspm 2,\,k$ is negative for all four solutions. Therefore the two improper positive $(X,\widehat Y)$ solutions $(9,\,231)$ and $(303,\,315)$  will not lead to a final solution.\pn
Together with $(153,\, 255)$ from the $\bf \boldsymbol \alpha)$ part they give all three positive improper solutions of $a\sspeq 3^2\cdot 17^2\cdot 41$.
\pn  
For the proper solutions representing $-17^2\cdot41$ see type {\bf a) iv)}. 
\pn 
For the truth of the the conjecture for the two family $II$ proper solutions representing $-17^2\cdot41$, \viz\, $(\tilde X_i,\,\tilde{\widehat Y_i})$, with $i\sspeq i^{\prime\prime}\sspp 1$, see type {\bf a)iv)} above. 
\pn
The corresponding two first positive improper solutions representing $a$, \viz\, $(X_i,\,\widehat Y_i)$, for  $i\sspeq i^{\prime\prime}\sspp 1$, are $(897,\, 675$ and $(351,\, 339)$. They satisfy the conjecture because all four above given $t\sspm 2\,k)$ values are negative.
\psn   
Therefore there is no $DS$-triple for $a\sspeq 3^2\cdot 17^2\cdot 41\sspeq 106641 $ of this type $\bf b)$.
\pbn
{\bf Type c) $(\bf X,\,\widehat Y)$ proper and  $\bf (t,\,k)$ improper}. 
\psn
{\bf i)} $a\sspeq 7^2\cdot 2\cdot 17 \sspeq 1666$. $A\sspeq 7\sspeq$\seqnum{A058529}$(2)$, $\tilde a\sspeq 2\cdot 17\sspeq$\seqnum{A192453}$(4)$. $(\tilde t,\,\tilde k)\sspeq (4,\, 3)$,  $(t,\,k)\sspeq (28,\, 21)$.
Because $t\sspm 2\,k$ is negative there will be no final solution, whatever the two proper fundamental family $I$ solutions $(X,\, \widehat Y)$ representing  $a\sspeq 1666$ are (they are $(X,\, \widehat Y)_1\sspeq (4,\, 29)$ and $(X,\, \widehat Y)_2\sspeq (16,\, 31)$ with $(X,\,Y)_1\sspeq (4,\,50)$ and $(X,\,Y)_2\sspeq (16,\,52)$, but $c_2 \sspgr c_3$ in both cases).
\pn
The conjecture for the two first proper positive family $II$ solutions $(X_i,\, \widehat Y_i)$, for $i\sspeq i^{\prime\prime}\sspp 1$, representing $1666$ (they are $(104,\, 79)$ and $(76,\,61)$) is true because $t\sspm 2\, k$ is negative. Only the second solution satisfies $0\sspkl X_i\sspkl Y_i$, with $k\sspeq 21$, but it fails the $c_2\sspkl c_3$ test.
\psn
{\bf ii)} $a\sspeq 7^2\cdot 89\sspeq 4361$. $A\sspeq 7$, $\tilde a\sspeq 89\sspeq$\seqnum{A192453}$(8)$. $(\tilde t,\,\tilde k)\sspeq (9,\, 2)$,  $(t,\,k)\sspeq (63,\, 14)$. The two family $I$ positive proper fundamental solutions are $(X,\, \widehat Y)_1\sspeq (29,\,51)$ and $(X,\, \widehat Y)_2\sspeq (51,\,59)$. Because $t\sspm 2\,k\sspeq 35$ only the first solution has to be tested. We skip the subscript $1$: $(X,\,Y)\sspeq (29,\,65)$, $c_1\sspeq (65\sspm 29)/2\sspeq 18$, $c_2\sspeq 47$, $c_3\sspeq (\widehat Y \sspp t\sspm k)/2 \sspeq 50$, and  $2\,q\sspeq Y\sspp t\sspeq 128$. Therefore a final $DS$-triple for $a\sspeq 4361$ is 
$[c_1,\,c_2,\,c_3]\sspeq [18,\, 47,\, 50)]$, with $(c_{4.-},\,c_{4.+}) \sspeq (-13,\,243)$. 
\pn
The conjecture for the two proper family $II$ members $(X_i,\,\widehat Y_i)$, for $i \sspeq i^{\prime\prime}\sspp 1$, representing $4361$, \viz\,  $(117,\,95)$ and $(83,\,75)$, is true because with $k\sspeq 14$ only the second solution looks promising, but because $83 \sspgr 35$ it fails the $c_2\sspkl c_3$ test.
\pn
Therefore the found $DS$-triple is the only one for $a\sspeq 4361$.
\pn
In fact $a \sspeq 7^2\cdot 89\sspeq 4361$ is the smallest instance of type $c)$ with a $DS$-triple. This is because $7^2\cdot 17$,\, $7^2\cdot 34$,\, $7^2\cdot 41$ and  $7^2\cdot 82$ all do not pass the $X\sspkl t\sspm 2\,k$ test. The instance $17^2\cdot 17\sspeq 4913$ (a candidate to be tested) is already larger. The type $\bf a)$ solution for $a\sspeq17^3$ does also not pass this test because $t\sspm 2\,k \sspeq 45\sspm 2\cdot 38 \sspkl 0$. 
\psn
{\bf iii)} $a\sspeq 7^2\cdot 17^2\cdot 2\sspeq 28322$. 
\pn
{\boldmath$\alpha$)} $A\sspeq 7\cdot 17\sspeq 119\sspeq$\seqnum{A058529}$(16)$, $\tilde a\sspeq 2$. There is no final solution because $(\tilde t,\,\tilde k)\sspeq (0,\,1)$, hence $t\sspeq 0$. See {\sl Lemma 17}.
\psn
{\boldmath $\beta$)} $A\sspeq 7$, $\tilde a\sspeq 17^2\cdot 2\sspeq 578\sspeq$\seqnum{A192453}$(39)$. $(\tilde t,\,\tilde k)\sspeq (24,\,1)$, $(t,\,k)\sspeq (168,\,7)$, with $t\sspm 2\,k\sspeq 154$. $(X,\, \widehat Y)_1\sspeq (44,\, 123)$ with $(X,\,Y)_1\sspeq (44,\,130)$. $c_{3,1}\sspeq (123\sspp 161)/2\sspeq 142$, $2\,q_1\sspeq 130\sspp 168 \sspeq 298$. $c_{1,1}\sspeq (130\sspm 44)/2\sspeq 43$, $c_{2,1}\sspeq 87$, Hence $(c_1,\,c_2,\,c_3)_1\sspeq (43,\, 87,\,142)$ is a $DS$-triple for $a\sspeq 28322$ with $(c_{4.-},\,c_{4.+})_1 \sspeq (-26,\,570)$.
\pn
$(X,\, \widehat Y)_2\sspeq (96,\, 137)$ with $(X,\,Y)_2\sspeq (96,\,144)$. $c_{3,2}\sspeq (137\sspp 161)/2\sspeq 149$, $2\,q_2\sspeq 144\sspp 168\sspeq 312$.      
$c_{1,2}\sspeq (144\sspm 96)/2\sspeq 24$, $c_{2,2}\sspeq 120$, Hence $(c_1,\,c_2,\,c_3)_2\sspeq (24,\, 120,\,149)$ is a second $DS$-triple for $a\sspeq 28322$ with $(c_{4.-},\,c_{4.+})_2 \sspeq (-19,\,605)$.
\psn
There is no other splitting $A\sspeq 17$ and $\tilde a\sspeq 7^2\cdot2$ because $\modan{7}{8}\neq 1$.
\pn
The critical first proper positive family $II$ members for the two solutions $(X_i,\,\widehat Y_i)$, corresponding to the two family $I$ solutions, are $(624,\, 457)$ and $(220,\,219)$. Because $k\sspeq 7$ only the second solution looks promising because $0\sspkl 220\sspkl 226$, but $X_i\sspeq 220\sspgeq 154$, failing the $c_2\sspkl c_3$ test. Also the first. already discarded, solution fails this test. Therefore the conjecture is true for these two examples. 

\pbn
\hrule
\pbn
Keywords: Descartes circle Theorem, Pythagorean triangle, binary quadratic forms: indefinite and positive definite, sums of two, three and four squares.  \psn
AMS MSC number: 11D09, 11E16, 15B36    \psn
OEIS A-numbers: \seqnum{A000010}, \seqnum{A000037}, \seqnum{A000217}, \seqnum{A001109}, \seqnum{A001132}, \seqnum{A001541}, \seqnum{A001653}, \seqnum{A007519}, \seqnum{A007520}, \seqnum{A007522}, \seqnum{A020882}, \seqnum{A023022}, \seqnum{A033200}, \seqnum{A033203}, \seqnum{A035251}, \seqnum{A038873}, \seqnum{A054521}, \seqnum{A055034}, \seqnum{A057126}, \seqnum{A058529}, \seqnum{A192453}, \seqnum{A222946}, \seqnum{A225771}, \seqnum{A249866}, \seqnum{A343238},\seqnum{A343240}.
\vfill
\eject
\noindent
\begin{center}
{\large {\bf Table 1: Primitive Cartesian triples of type $\bf [c, c, d]$  and $\bf [c,d,d]$}} \end {center}
\begin{center}  
\begin{tabular}{|c|c|c|c||c|c|c|c|}\hline
&& && && &\\
\hskip .2cm$\bf [n, m]$ & $\bf [c ,\,c ,\, d]$  &$\bf q $ & $\bf [c_{4,-},\,c_{4,+}]$& $\bf  [n, m]$ & $\bf [c,\, d,\,d]$  &$\bf q $ & $\bf [c_{4,-},c_{4,+}]$  \\
&& && && &\\ \hline\hline
$\bf [2, 1]$ & $[1, 1, 4]$ & $3$ & $[0, 12] $ &                    $\bf  [3, 2] $& $[5, 8, 8] $ & $12 $ & $[-3, 45] $\\  
$\bf [2, 1]$ & $[2, 2, 3]$ & $4$ & $[-1, 15] $ &                   $\bf  [4, 1] $& $[8, 9, 9] $ & $15 $ & $[-4, 56] $\\  
$\bf [3, 2]$ & $[1, 1, 12 ]$ & $5 $ & $[4, 24] $ &                $\bf [4, 3] $& $[7, 18, 18]$ & $ 24$ & $ [-5, 91] $\\ 
$\bf [4, 1]$ & $[2, 2, 15] $ & $8$ & $ [3, 35] $ &                $\bf [5, 4] $ & $[9, 32, 32]$ & $40 $ & $[-7, 153 ]$  \\
$\bf [4, 3] $ & $[1, 1, 24]$ & $ 7 $ & $[12, 40]  $ &            $\bf [6, 1] $& $[12, 25, 25 ]$ & $35  $ & $[-8, 132]$ \\
$\bf  [5, 2]$ & $[9, 9, 20]$ & $21 $ & $[-4, 80] $ &              $\bf  [6, 5] $ &  $[11, 50, 50] $ & $60 $ & $[-9, 231] $\\
$\bf  [5, 2]$ & $[8, 8, 21]$ & $ 20$ & $[-3, 77] $ &              $\bf  [7, 6] $ &  $[13, 72, 72] $&  $ 84 $ & $[-11, 325] $\\
$\bf [5, 4] $& $[1, 1, 40]$ & $9 $ & $[24, 60] $ &                 $\bf [8, 1]$ & $[16, 49, 49]$ & $63 $ & $[-12, 240] $ \\
$\bf [6, 1]  $& $[2, 2, 35]$ & $12 $ & $[15, 63] $ &              $\bf [8, 5] $ & $[39, 50, 50]$ & $ 80$ & $[-21, 299] $ \\
$\bf  [6, 5] $& $[1, 1, 60]$ & $11 $ & $[40, 84]$ &                $\bf [8, 7]$ & $[15, 98, 98]$ & $112$ & $[-13, 435] $ \\
$\bf [7, 2] $& $[25, 25, 28]$ & $45 $ & $[-12, 168 ]$ &         $\bf [9, 2] $ & $[36, 49, 49]  $ & $77 $ & $[-20, 288]$ \\
$\bf [7, 2] $& $[8, 8, 45]$ & $28 $ & $[5, 117 ]$ &                 $\bf [9, 8]  $ & $[17, 128, 128]$ & $ 144$ & $[-15, 561]$ \\
$\bf [7, 4] $& $[9, 9, 56]$ & $33 $ & $[8, 140] $ &                 $\bf [10, 1] $ & $[20, 81, 81] $ & $ 99 $ & $[-16, 380 ] $ \\
$\bf [7, 4] $& $[32, 32, 33]$ & $56 $ & $[-15, 209] $ &          $\bf [10, 7] $ & $ [51, 98, 98]$ & $140 $ & $[-33, 527] $ \\
$\bf [7, 6] $& $[1, 1, 84]$ & $13$ & $[60, 112] $ &                $\bf [10, 9]   $ & $[19, 162, 162]$ & $180 $ & $[-17, 703] $ \\
$\bf [8, 1]  $& $[2, 2, 63]$ & $16$ & $[35, 99]$ &                 $\bf [11, 2]$ & $[44, 81, 81]$ & $117$ & $[-28, 440] $ \\
$\bf [8, 3]  $& $[25, 25, 48]$ & $ 55 $ & $[-12, 208]$ &         $\bf [11, 8]$ & $[57, 128, 128]$ & $176$ & $[-39, 665] $ \\
$\bf [8, 3]  $& $[18, 18, 55]$ & $ 48 $ & $[-5, 187]$ &           $\bf [11, 10]   $ & $[21, 200, 200]$ & $220$ & $[-19, 861] $ \\
$\bf [8, 5] $& $[9, 9, 80]$ & $39$ & $[20, 176] $&                  $\bf [12, 1]$ & $[24, 121, 121]$ & $143$ & $[-20, 552] $ \\
$\bf [8, 7] $& $[1, 1, 112]$ & $15$ & $[84, 144] $ &               $\bf [12, 7]$ & $[95, 98, 98]$ & $168$ & $[-45, 627] $ \\
$\bf [9, 2]$& $[8, 8, 77]$ & $36 $ & $[21, 165] $ &                 $\bf [12, 11]$ & $[23, 242, 242]$ & $264$ & $[-21, 1035] $ \\
$\bf [9, 4]$& $[25, 25, 72]$ & $65$ & $[-8, 252] $ &                $\bf [13, 2] $ & $[52, 121, 121]$ & $165 $ & $[-36, 624 ] $ \\
$\bf [9, 4]$& $[32, 32, 65]$ & $72 $ & $[-15, 273] $ &             $\bf [13, 8] $ & $[105, 128, 128]$ & $208$ & $[-55, 777] $ \\
$\bf [9, 8]$& $[1, 1, 144]$ & $17$ & $[112, 180] $ &                $\bf [13, 10] $ & $[69, 200, 200]$ & $260$ & $[-51, 989] $ \\
$\bf [10, 1] $& $[2, 2, 99]$ & $20 $ &  $[63, 143] $ &               $\bf [13, 12]$ & $[25, 288, 288]$ & $312$ & $[-23, 1225] $ \\
$\bf [10, 3]$& $[49, 49, 60]$ & $91$ & $[-24, 340] $ &             $\bf  [14, 1] $ & $[28, 169, 169]$ & $195 $ & $[-24, 756] $ \\
$\bf [10, 3]$& $[18, 18, 91]$ & $ 60$ & $[7, 247] $ &                $\bf [14, 3]$ & $[84, 121, 121]$ & $ 187$ & $[-48, 700] $ \\
$\bf [10, 7] $& $[9, 9, 140]$ & $ 51$ & $[56, 260] $ &               $\bf [14, 9]$ & $[115, 162, 162 ]$ & $  252$ & $[-65, 943] $ \\
$\bf [10, 9] $& $[1, 1, 180]$ & $19$ & $[144, 220] $ &              $\bf [14, 11]$ & $[75, 242, 242]$ & $308 $ & $[-57, 1175] $ \\
$\bf [11, 2] $& $[8, 8, 117]$ & $ 44 $ & $[45, 221] $ &              $\bf [14, 13]$ & $[27, 338, 338]$ & $364$ & $[-25, 1431] $ \\
$\bf [11, 4] $& $[49, 49, 88]$ & $ 105$ & $[-24, 396] $ &          $\bf [15, 2]$ & $[60, 169, 169] $ & $221 $ & $[-44, 840]  $ \\
$\bf [11, 4] $& $[32, 32, 105]$ & $88  $ & $[-7, 345] $ &           $\bf [15, 4]$ & $[120, 121, 121] $ & $209 $ & $[-56, 780 ]  $ \\
$\bf [11, 6] $& $[25, 25, 132]$ & $ 85 $ & $[12, 352] $ &           $\bf [15, 14]$ & $[29, 392, 392]$ & $420$ & $[-27, 1653] $ \\
$\bf [11, 6] $& $[72, 72, 85]$ & $132  $ & $[-35, 493] $ &           $\bf [16, 1]$ & $[32, 225, 225]$ & $ 255$ & $[-28, 992] $ \\
$\bf [11, 8]    $& $[9, 9, 176]$ & $57 $ & $[80, 308] $ &              $\bf [16, 3]$ & $[96, 169, 169]$ & $ 247$ & $[-60, 928] $ \\
$\bf [11, 10] $& $[1, 1, 220]$ & $21$ & $[180, 264] $ &              $\bf [16, 11] $ & $[135, 242, 242]$ & $ 352$ & $[-85, 1323] $ \\
$\bf [12, 1]    $& $[2, 2, 143]$ & $24$ & $[99, 195] $ &               $\bf [16, 13] $ & $[87, 338, 338]$ & $416$ & $[-69, 1595] $ \\
$\bf  [12, 5]   $& $[49, 49, 120]$ & $119 $ & $[-20, 456] $ &        $\bf [16, 15] $ & $[31, 450, 450]$ & $ 480$ & $[-29, 1891] $ \\
$\bf  [12, 5]   $& $[50, 50, 119]$ & $ 120  $ & $[-21, 459] $ &       $\bf [17, 2] $ & $[68, 225, 225]$ & $ 285$ & $[-52, 1088] $ \\
$\bf  [12, 7]   $& $[25, 25, 168]$ & $ 95 $ & $[28, 408 ] $ &          $\bf [17, 4]$ & $[136, 169, 169]$ & $273$ & $[-72, 1020] $ \\
$\bf [12, 11] $& $[1, 1, 264]$ & $23 $ & $[220, 312] $ &                $\bf [17, 10]$ & $[189, 200, 200] ]$ & $ 340$ & $[-91, 1269] $ \\
$\bf [13, 2]    $& $[8, 8, 165]$ & $ 52$ & $[77, 285] $ &                $\bf  [17, 12]$ & $[145, 288, 288]$ & $ 408$ & $[-95, 1537] $ \\
$\bf [13, 4]    $& $[81, 81, 104]$ & $153$ & $[-40, 572] $ &           $\bf [[17, 14]$ & $[93, 392, 392]$ & $ 476$ & $[-75, 1829] $ \\
$\bf [13, 4]  $& $[32, 32, 153]$ & $104 $ & $[9, 425 ] $ &             $\bf  [17, 16]$ & $[33, 512, 512]$ & $ 544$ & $[-31, 2145] $ \\
$\bf  cont'd  $& $ ...$ & $... $ & $... $ &                                          $\bf  ...$ & $...$ & $...$ & $ ... $ \\
\hline
\end{tabular}
\end{center}
\vfill
\eject
\noindent
\begin{center}
{\large {\bf Table 1 continued: Primitive Cartesian triples of type $\bf [c, c, d]$,  $n\sspleq 17$}} 
\end{center}
\begin{center}  
\begin{tabular}{|c|c|c|c||c|c|c|c|}\hline
&& && && &\\
\hskip .2cm$\bf [n, m]$ & $\bf [c ,\,c ,\, d]$  &$\bf q $ & $\bf [c_{4,-},\,c_{4,+}]$& $\bf  [n, m]$ & $\bf [c,\, d,\,d]$  &$\bf q $ & $\bf [c_{4,-},c_{4,+}]$  \\
&& && && &\\ \hline\hline
$\bf [13, 6]$ & $[49, 49, 156]$ & $133$ & $[-12, 520] $ &                    $\bf ...   $& $... $ & $... $ & $... $\\  
$\bf [13, 6]$ & $[72, 72, 133]$ & $156$ & $[-35, 589] $ &                   $\bf   $& $ $ & $$ & $$\\  
$\bf [13, 8]$ & $[25, 25, 208]$ & $105 $ & $[48, 468] $ &                $\bf  $& $$ & $ $ & $  $\\ 
$\bf [13, 10]$ & $[9, 9, 260] $ & $69$ & $ [140, 416] $ &                $\bf  $ & $$ & $$ & $$  \\
$\bf [13, 12] $ & $[1, 1, 312]$ & $ 25$ & $[264, 364]  $ &             $\bf $& $$ & $  $ & $$ \\
$\bf  [14, 1]$ & $[2, 2, 195]$ & $28 $ & $[143, 255] $ &                $\bf  $ &  $ $ & $$ & $$\\
$\bf  [14, 3]$ & $[18, 18, 187]$ & $ 84$ & $[55, 391] $ &               $\bf  $ &  $ $&  $  $ & $ $\\
$\bf [14, 5] $& $[81, 81, 140]$ & $171 $ & $[-40, 644] $ &              $\bf $ & $$ & $ $ & $$ \\
$\bf [14, 5]  $& $[50, 50, 171]$ & $140$ & $[-9, 551] $ &               $\bf $ & $$ & $ $ & $$ \\
$\bf  [14, 9] $& $[25, 25, 252]$ & $115 $ & $[72, 532]$ &                $\bf $ & $  $ & $$ & $ $ \\
$\bf [14, 11] $& $[9, 9, 308]$ & $ 75 $ & $[176, 476]$ &                 $\bf  $ & $   $ & $$ & $ $ \\
$\bf [14, 13] $& $[1, 1, 364]$ & $27 $ & $[312, 420] $ &                  $\bf  $ & $   $ & $$ & $ $ \\
$\bf [15, 2] $& $[8, 8, 221]$ & $60$ & $[117, 357] $ &                  $\bf  $ & $   $ & $$ & $ $ \\
$\bf [15, 4] $& $[32, 32, 209]$ & $120,6 $ & $ [33, 513] $ &           $\bf  $ & $   $ & $$ & $ $ \\
$\bf [15, 8] $& $[49, 49, 240]$ & $161$ & $[16, 660] $ &                 $\bf  $ & $   $ & $$ & $ $ \\
$\bf [15, 8 ]  $& $[128, 128, 161]$ & $240$ & $[-63, 897]$ &           $\bf  $ & $   $ & $$ & $ $ \\      
$\bf [15, 14]  $& $[1, 1, 420]$ & $ 29 $ & $[364, 480]$ &          $\bf  $ & $   $ & $$ & $ $ \\
$\bf [16, 1]  $& $[2, 2, 255]$ & $ 32 $ & $[195, 323]$ &            $\bf  $ & $   $ & $$ & $ $ \\
$\bf [16, 3] $& $[18, 18, 247]$ & $96$ & $[91, 475] $&                   $\bf  $ & $   $ & $$ & $ $ \\
$\bf [16, 5] $& $[121, 121, 160]$ & $231$ & $[-60, 864] $ &         $\bf  $ & $   $ & $$ & $ $ \\
$\bf [16, 5 ]$& $[50, 50, 231]$ & $160$ & $[11, 651] $ &              $\bf  $ & $   $ & $$ & $ $ \\
$\bf [16, 7]$& $[81, 81, 224]$ & $ 207$ & $[-28, 800] $ &               $\bf  $ & $   $ & $$ & $ $ \\
$\bf [16, 7]$& $[98, 98, 207]$ & $224 $ & $[-45, 851] $ &          $\bf  $ & $   $ & $$ & $ $ \\
$\bf [16, 9]$& $[49, 49, 288]$ & $175$ & $[36, 736] $ &               $\bf  $ & $   $ & $$ & $ $ \\
$\bf [16, 9] $& $[162, 162, 175]$ & $ 288 $ &  $[-77, 1075] $ &        $\bf  $ & $   $ & $$ & $ $ \\    
$\bf [16, 11]$& $[25, 25, 352]$ & $135$ & $[132, 672] $ &              $\bf  $ & $   $ & $$ & $ $ \\    
$\bf [16, 13]$& $[9, 9, 416]$ & $ 87$ & $[260, 608] $ &                 $\bf  $ & $   $ & $$ & $ $ \\    
$\bf [16, 15] $& $[1, 1, 480]$ & $ 31$ & $[420, 544] $ &                $\bf  $ & $   $ & $$ & $ $ \\    
$\bf [17, 2] $& $[8, 8, 285]$ & $68$ & $[165, 437] $ &             $\bf  $ & $   $ & $$ & $ $ \\    
$\bf [17, 4] $& $[32, 32, 273]$ & $ 136 $ & $[65, 609] $ &              $\bf  $ & $   $ & $$ & $ $ \\    
$\bf [17, 6] $& $[121, 121, 204]$ & $ 253$ & $[-60, 952] $ &         $\bf  $ & $   $ & $$ & $ $ \\    
$\bf [17, 6] $& $[[72, 72, 253]$ & $204  $ & $[-11, 805] $ &            $\bf  $ & $   $ & $$ & $ $ \\    
$\bf [17, 8] $& $[81, 81, 272]$ & $225, $ & $[-16, 884] $ &            $\bf  $ & $   $ & $$ & $ $ \\    
$\bf [17, 8] $& $[128, 128, 225]$ & $272  $ & $[-63, 1025] $ &        $\bf  $ & $   $ & $$ & $ $ \\    
$\bf [17, 10] $& $[49, 49, 340]$ & $189, $ & $[60, 816] $ &              $\bf  $ & $   $ & $$ & $ $ \\    
$\bf [17, 12] $& $[25, 25, 408]$ & $145,$ & $[168, 748] $ &           $\bf  $ & $   $ & $$ & $ $ \\    
$\bf [17, 14] $& $[9, 9, 476]$ & $ 93$ & $[308, 680] $ &              $\bf  $ & $   $ & $$ & $ $ \\    
$\bf  [17, 16]   $& $[1, 1, 544]$ & $ 33$ & $[480, 612] $ &             $\bf  $ & $   $ & $$ & $ $ \\    
$\bf  ...  $& $ ...$ & $... $ & $... $ &                                          $\bf   $ & $ $ & $ $ & $  $ \\
\hline
\end{tabular}
\end{center}
\vfill
\eject
\noindent
\begin{center}
{\large {\bf Table 2: Distinct primitive Cartesian triples with curvatures $\bf c_3 \sspeq 1.. 38$}} 
\end {center}
\begin{center}  
\begin{tabular}{|l|c|c|c||l|c|c|c|}\hline
&& && && &\\
\hskip .2cm$\bf c_3$ & $\bf [c_1,\,c_2,\,c_3]$  &$\bf q $ & $\bf [c_{4,-},\,c_{4,+}]$& $\bf c_3$ & $\bf [c_1,c_2,c_3]$  &$\bf q $ & $\bf [c_{4,-},c_{4,+}]$  \\
&& && && &\\ \hline\hline
$\bf 6$ & $[2, 3, 6]$ & $6$ & $[-1, 23] $ &                      $\bf     $& $[4, 9, 28]$ & $ 20$ & $[1, 81] $\\  
$\bf 7$ & $[3, 6, 7]$ & $9$ & $[-2, 34] $ &                      $\bf     $& $[4, 21, 28]$ & $ 28$ & $[-3, 109] $\\ 
$\bf 8$ & $ \emptyset$ & $\emptyset$ & $\emptyset $ &  $\bf     $& $[13, 21, 28]$ & $ 35$ & $[-8, 132]$  \\
$\bf 9$ & $[1, 4, 9]$ & $7$ & $[0, 28]  $ &                       $\bf     $& $[21, 24, 28]$ & $ 42$ & $[-11, 157]$ \\
$\bf 10$ & $[3, 7, 10]$ & $11$ & $[-2, 42] $ &                 $\bf 29$ & $[5, 24, 29]$ & $31$ & $[-4, 120] $\\
$\bf 11 $& $[2, 6, 11]$ & $10$ & $[-1, 39] $ &                 $\bf     $ & $[8, 12, 29]$ & $26$ & $[-3, 101] $ \\
$\bf 12 $& $[1, 4, 12]$ & $8$ & $[1, 33] $ &                     $\bf 30$ & $[6, 11, 30]$ & $24$ & $[-1, 95] $ \\
$\bf      $& $[5, 8, 12]$ & $14$ & $[-3, 53]$ &                   $\bf    $ & $[7, 22, 30]$ & $32$ & $[-5, 123] $ \\
$\bf 13 $& $[4, 12, 13]$ & $16$ & $[-3, 61]$ &                $\bf 31$ & $[3, 22, 31]$ & $29$ & $[-2, 114]$ \\
$\bf 14 $& $[3, 6, 14]$ & $12$ & $[-1, 47] $ &                 $\bf     $ & $[6, 30, 31]$ & $36$ & $[-5, 139] $ \\
$\bf 15 $& $[2, 3, 15]$ & $9$ & $[2, 38] $ &                     $\bf    $ & $[10, 19, 31]$ & $33$ & $[-6, 126] $ \\
$\bf     $& $[3, 10, 15]$ & $15$ & $[-2, 58]$ &                 $\bf 32$ & $[5, 20, 32]$ & $30$ & $[-3, 117] $ \\
$\bf     $& $[11, 14, 15]$ & $23$ & $[-6, 86]$ &                $\bf    $ & $[21, 29, 32]$ & $47$ & $[-12, 176] $ \\
$\bf 16 $& $[1, 9, 16]$ & $13$ & $[0, 52] $&                    $\bf 33$ & $[1, 4, 33]$ & $13$ & $[12, 64] $ \\
$\bf     $& $[4, 13, 16]$ & $18$ & $[-3, 69] $ &                 $\bf     $ & $[1, 12, 33]$ & $21$ & $[4, 88] $ \\
$\bf 17$& $[8, 9, 17]$ & $19$ & $[-4, 72] $ &                   $\bf      $ & $[4, 12, 33]$ & $24$ & $[1, 97] $ \\
$\bf 18$& $[2, 11, 18]$ & $16$ & $[-1, 63] $ &                 $\bf     $ & $[8, 17, 33]$ & $31$ & $[-4, 120] $ \\
$\bf 19$& $[6, 7, 19]$ & $17$ & $[-2, 66] $ &                  $\bf     $ & $[12, 20, 33]$ & $36$ & $[-7, 137] $ \\
$\bf    $& $[10, 15, 19]$ & $25$ & $[-6, 94] $ &                $\bf 34$ & $[3, 6, 34]$ & $18$ & $[7, 79] $ \\
$\bf 20$& $[5, 12, 20]$ & $20$ & $[-3, 77] $ &                 $\bf      $ & $[3, 7, 34]$ & $19$ & $[6, 82] $ \\
$\bf    $& $[12, 17, 20]$ & $28$ & $[-7, 105] $ &              $\bf      $ & $[6, 7, 34]$ & $22$ & $[3, 91] $ \\
$\bf 21 $& $[4, 16, 21]$ & $22$ & $[-3, 85] $ &                 $\bf      $ & $[6, 31, 34]$ & $38$ & $[-5, 147] $ \\
$\bf      $& $[5, 20, 21]$ & $25$ & $[-4, 96] $ &                $\bf     $ & $[18, 19, 34]$ & $40$ & $[-9, 151] $ \\
$\bf 22 $& $[3, 15, 22]$ & $21$ & $[-2, 82 ] $ &               $\bf     $ & $[27, 31, 34]$ & $53$ & $[-14, 198] $ \\
$\bf     $& $[7, 18, 22]$ & $26$ & $[-5, 99] $ &                 $\bf 35$ & $[2, 15, 35]$ & $25$ & $[2, 102] $ \\
$\bf     $& $[18, 19, 22]$ & $34$ & $[-9, 127] $ &              $\bf     $ & $[6, 14, 35]$ & $28$ & $[-1, 111] $ \\
$\bf 23 $& $[2, 3, 23]$ & $11$ & $[6, 50] $ &                    $\bf     $ & $[14, 15, 35]$ & $35$ & $[-6, 134 ] $ \\
$\bf     $& $[2, 6, 23]$ & $14$ & $[3, 59] $ &                    $\bf     $ & $[14, 26, 35]$ & $42$ & $[-9, 159] $ \\
$\bf     $& $[3, 6, 23]$ & $15$ & $[2, 62] $ &                     $\bf   $ & $[18, 23, 35]$ & $40$ & $[-10, 162] $ \\
$\bf     $& $[11, 14, 23]$ & $27$ & $[-6, 102] $ &             $\bf 36$ & $[1, 25, 36]$ & $31$ & $[0, 124] $ \\
$\bf 24 $& $[1, 12, 24]$ & $18$ & $[1, 73] $ &                  $\bf     $ & $[5, 29, 36]$ & $37$ & $[-4, 144] $ \\
$\bf     $& $[5, 21, 24]$ & $27$ & $[-4, 104] $ &               $\bf     $ & $[9, 17, 36]$ & $33$ & $[-4, 12] $ \\
$\bf     $& $[12, 17, 24]$ & $30$ & $[-7, 113] $ &             $\bf     $ & $[9, 32, 36]$ & $42$ & $[-7, 161] $ \\
$\bf     $& $[13, 21, 24]$ & $33$ & $[-8, 124] $ &             $\bf 37$ & $[4, 28, 37]$ & $36$ & $[-3, 141] $ \\
$\bf 25 $& $[1, 16, 25]$ & $21$ & $[0, 84] $ &                  $\bf       $ & $[13, 24, 37]$ & $41$ & $[-8, 156] $ \\
$\bf     $& $[4, 9, 25]$ & $19$ & $[0, 76] $ &                     $\bf       $ & $[16, 36, 37]$ & $50$ & $[-11, 189] $ \\
$\bf 26 $& $[3, 14, 26]$ & $22$ & $[-1, 87] $ &                  $\bf      $ & $[21, 28, 37]$ & $49$ & $[-12, 184] $ \\
$\bf     $& $[11, 15, 26]$ & $29$ & $[-6, 110] $ &               $\bf 38$ & $[2, 3, 38]$ & $14$ & $[15, 71] $ \\
$\bf 27 $& $[2, 18, 27]$ & $24$ & $[-1, 95] $ &                 $\bf       $ & $[2, 15, 38]$ & $26$ & $[3, 107] $ \\
$\bf     $& $[7, 10, 27]$ & $23$ & $[-2, 90] $ &                  $\bf       $ & $[2, 27, 38]$ & $34$ & $[-1, 135] $ \\
$\bf     $& $[14, 26, 27]$ & $ 38$ & $[-9, 143] $ &             $\bf       $ & $[3, 15, 38]$ & $27$ & $[2, 110] $ \\
$\bf     $& $[18, 23, 27]$ & $ 39$ & $[-10, 146] $ &            $\bf       $ & $[14, 27, 38]$ & $44$ & $[-9, 167] $ \\
$\bf 28 $& $[1, 4, 28]$ & $12$ & $[9, 57]$ &                      $\bf       $ & $[23, 30, 38]$ & $52$ & $[-13, 195] $ \\          
$\bf     $& $[1, 9, 28]$ & $ 17$ & $[4, 72] $&                       $\bf ...    $& $...$ & $... $ & $...$\\                                           \hline
\end{tabular}
\end{center}
 \vfill
\eject
\noindent
\begin{center}
{\large {\bf Table 3: Case i) $\bf q\sspeq c_3$. Proper solutions $\bf [x,\, y]\sspeq [c_1,\,c_2]$ of \Eq{35}, \\
and $\bf [X,\, Y]\sspeq [y\sspm x,y\sspp x]$ of \Eq{38}, $\bf q$,\\
and rpapfs $\bf t-tuples$, for $\bf s\sspin [1,313]$}} \\

\end {center}
\begin{center}  
\begin{tabular}{|c|c|c|c|l|}\hline
&& && \\
\hskip .2cm
$\bf s$ & $\bf [x,\,y]$  & $\bf [X,\,Y] $ &$\bf q\sspeq c_3$ & $\bf t- tuples$\\
&& &&\\ 
\hline\hline
$\bf 7$ & $[2,\,3]$ & $[1,\,5]$ &  $6$& $ (-1,4,1),\, [(-5,1)] $\\
\hline
$\bf 17$ & $[3,\,10]$ & $[7,\,13]$ &  $15$& $ (-1,6,2,1),\,[(-6,2,2)] $\\
\hline
$\bf 23$ & $[5,\,12]$ & $[7,\,17]$ &  $20$& $ (-1,3,2,3,1),\,[(-3,3,2,2)] $\\
\hline
$\bf 31$ & $[4,\,21]$ & $[17,\,25]$ &  $28$& $ (-1,8,2,2,1),\,[(-8,3,2)] $\\
\hline
$\bf 41$ & $[14,\,15]$ & $[1,\,29]$ &  $35$& $ (-1,28,1),\,[(-29,1)] $\\
\hline
$\bf 47$ & $[7,\,30]$ & $[23,\,37]$ &  $42$& $ (-4,2,3,2,1),\,[(-1,4,3,3,2)] $\\
\hline
$\bf \underline{49}$ & $[5,\,36]$ & $[31,\,41]$ &  $45$& $ (-1,10,2_3,1),\,[(-10,4,2)] $\\
\hline
$\bf 71$ & $[6,\,55]$ & $[49,\,61]$ &  $66$& $ (-1,12,2_4,1),\,[(-12,5,2)] $\\
\hline
$\bf 73$ & $[18,\,35]$ & $[17,\,53]$ &  $63$& $ (-1,2_8,4,1),\,[(-2,9,2_3)]  $\\
\hline
$\bf 79$ & $[9,\,56]$ & $[47,\,65]$ &  $72$& $ (-1,3_3,2,2,1),\,[(-3,2,3,4,2)] $\\
\hline
$\bf 89$ & $[21,\,44]$ & $[23,\,65]$ &  $77$& $ (-1,3,6,3,1),\,[(-3,2_4,3,2,2)]  $\\
\hline
$\bf 97$ & $[7,\,78]$ & $[71,\,85]$ &  $91$& $ (-1,14,2_5,1),\,[(-14,6,2)]  $\\
\hline
$\bf 103$ & $[33,\,40]$ & $[7,\,73]$ &  $88$& $ [(-1,3,2,11,1),\,[(-3,3,2_{10})] $\\
\hline
$\bf 113$ & $[22,\,63]$ & $[41,\,85]$ &  $99$& $ (-2,2_{13},3,1),\,[(-1,2,15,2,2)]  $\\
\hline
$\bf \underline{119}$ & $[24,\,65]$ & $[41,\,89]$ &  $104$& $ (-6,2_5,3,1),\,[(-1,6,7,2,2)]  $\\
& $[11,\,90]$ & $[79,\,101]$ &  $110$& $ (-4,3,3,2_3,1),\,[(-1,4,2,3,5,2)]   $\\
\hline
$\bf 127$ & $[8,\,105]$ & $[97,\,113]$ &  $120$& $ (-1,16,2_6,1),\,[(-16,7,2)]$\\
\hline
$\bf 137$ & $[45,\,52]$ & $[7,\,97]$ &  $117$& $(1,6,14,1),\,[(-6,2_{14})] $\\
\hline
$\bf 151$ & $[39,\,70]$ & $[31,\,109]$ &  $130$& $ (-15,2,4,1),\,[(-1,15,3,2_3)] $\\
\hline
$\bf \underline{161}$ & $[26,\,99]$ & $[73,\,125]$ &  $143$& $ (-1,10,2,4,2,1),\,[(-10,3,2,3,2)] $\\
&$[9,\,136]$ & $[127,\,145]$ &  $153$& $  (-1,18,2_7,1),\,[(-18,8,2)]  $\\
\hline
$\bf 167$ & $[13,\,132]$ & $[119,\,145]$ &  $156$& $ (-1,3,4,3,2_4,1),\,[(-3,2,2,3,6,2)] $\\
\hline
$\bf 191$ & $[30,\,119]$ & $[89,\,149]$ &  $170$& $ (-1,2_{14},4,2,1),\,[(-2,15,2,3,2)]  $\\
\hline
$\bf 193$ & $[60,\,77]$ & $[17,\,137]$ &  $165$& $ (-17,8,1),\,[(-1,17,2_8)] $\\
\hline
$\bf 199$ & $[10,\,171]$ & $[161,\,181]$ &  $190$& $ (-1,20,2_8,1),\,[(-20,9,2)]  $\\
\hline
$\bf \underline{217}$ & $[55,\,102]$ & $[47,\,157]$ &  $187$& $ (-15,2,2,4,1),\,[(-1,15,4,2_3)]  $\\
& $[30,\,143]$ & $[113,\,173]$ &  $195$& $ (-1,3,3,2_6,3,2,1),\,[(-3,2,9,3,2)]  $\\
\hline
$\bf 223$ & $[15,\,182]$ & $[167,\,197]$ &  $210$& $ (-4,4,3,2_5,1),\,[(-1,4,2,2,3,7,2 )] $\\
\hline
$\bf 233$ & $[33,\,152]$ & $[119,\,185]$ &  $209$& $ (-1,13, 2_3,3,2,1),\,[(-13,5,3,2)] $\\
\hline
$\bf 239$ & $[84,\,85]$ & $[1,\,169]$ &  $204$& $ (-1,168,1),\,[(-169,1)] $\\
\hline
$\bf 241$ & $[11,\,210]$ & $[199,\,221]$ &  $231$& $ (-1,22,2_9,1),\,[(-22,10,2)] $\\
\hline
$\bf 257$ & $[68,\,117]$ & $[49,\,185]$ &  $221$& $ (-1,5,2,5,4,1),\,[(-5,3,2,2,3,2_3)]  $\\
\hline
$\bf 263$ & $[60,\,133]$ & $[73,\,193]$ &  $228$& $ (-5,5,3,3,1),\,[(-1,5,2_3,3,3,2,2)]  $\\
\hline
$\bf 271$ & $[51,\,154]$ & $[103,\,205]$ &  $238$& $ (-1,102,2,1),\,[(-102,2,2)] $\\
\hline
$\bf 281$ & $[34,\,195]$ & $[161,\,229]$ &  $255$& $ (-3,2,4,4,2,2,1),\,[(-1,3,3,2,3,2,4,2)]  $\\
\hline
$\bf \underline{287}$ & $[17,\,240]$ & $[223,\,257]$ &  $272$& $ (-1,3,5,3,2_6, 1),\, [(-3,2,2,2,3,8,2)] $\\
&$[12,\,253]$ & $[241,\,265]$ &  $276$& $ (-1,24,2_{10},1),\,[(-24,11,2)]$\\
\hline
$\bf \underline{289}$ & $[91,\,114]$ & $[23,\,205]$ &  $247$& $ (-2,11,9,1),\,[(-1,2_{10},3,2_{8})]  $\\
\hline
$\bf 311$ & $[95,\,126]$ & $[31,\,221]$ &  $266$& $ (-3,2_7,8,1),\,[(-1,3,9,2_7)]  $\\
\hline
$\bf 313$ & $[65,\,168]$ & $[103,\,223]$ &  $273$& $  (-2,3,2_3,3,2,2,3,1),\,[(-1,2,2,6,5,2,2)] $\\
\hline
$...$ &$...$&$...$& $...$& $...$\\
\hline
\end{tabular}
\end{center}
\vfill
\eject
\begin{center}
{\large {\bf Table 4: Case ii) $\bf q\sspkl c_3$}}.\\ 
\end {center}
\begin{center}
\begin{tabular}{|c|l|l|c|c|c|l|l|}\hline
&& && && &\\
\hskip .2cm
$\bf c_3$ & $\bf [X,\,\widehat Y]$ & $\bf[c_1,c_2]$ & $\bf q$ &$\bf k$ & $\bf t$& $\bf a$ & $\bf [c_{4,-},\,c_{4,+}]$\\
&& && && &\\ 
\hline\hline
$\bf 9$ & $[3,\,7]$ & $[1,\,4]$ & $7$ & $2$ & $9$ & $89$ & $[0,\,28]$ \\
$\bf 11$ & $[4,\,9]$ & $[2,\,6]$ & $10$ & $1$ & $12$ & $146\sspeq 2\cdot 73$ & $[-1,39]$\\
$\bf 12$ & $[3,\,9]$& $[1,\,4]$& $8$ & $4$ & $11$ & $153\sspeq 3^2\cdot 17$ & $[1,\,33]$\\
$\bf 14$ & $[3,\,11]$ & $[3,\,6]$    & $12$ & $2$ & $15$ & $233$ & $[-1,\,47]$\\
$\bf 15$ & $[1,\,11]$ & $[2,\,3]$    & $9$ & $6$ & $13$ & $241$ & $[2,\,38]$\\
$\bf 16$ & $[8,\,13]$ & $[1,\,9]$    & $13$ & $3$ & $16$ & $274\sspeq 2\cdot 137$ & $[0,\,52]$\\
$\bf 18$ & $[9,\,15]$ & $[2,\,11]$   & $16$ & $2$ & $19$ & $369\sspeq 3^2\cdot 41$& $[-1,\,63]$ \\
$\bf 19$ & $[1,\,15]$ & $[6,\,7]$    & $17$ & $2$ & $21$ & $449$ & $[-2,\,66]$\\
$\bf 22$ & $[12,\,19]$ & $[3,\,15]$  & $21$ & $1$ & $24$ & $578\sspeq 2\cdot 17^2$& $[-2,\,82]$ \\
$\bf 23$ & $[1,\,17]$ & $[2,\,3]$   & $11$ & $12$ & $17$ & $577$ & $[6,\,50]$\\
         & $[4,\,17]$ & $[2,\,6]$   & $14$ & $9$ & $20$ & $562\sspeq 2\cdot 281$ & $[3,\,59]$\\
         & $[3,\,17]$ & $[3,\,6]$   & $15$ & $8$ & $21$ & $569$ & $[2,\,62]$\\
$\bf 24$ & $[11,\,19]$& $[1,\,12]$   & $18$ & $6$ & $23$ & $601$ & $[1,\,73]$\\ 
$\bf 25$ & $[15,\,21]$& $[1,\,16]$   & $21$ & $4$ & $25$ & $657\sspeq 3^2\cdot 73$ & $[0,\,84]$\\
         & $[5,\,19]$ & $[4,\,9]$   & $19$ & $6$ & $25$ & $697\sspeq 17\cdot 41$ & $[0,\,76]$\\
$\bf 26$ & $[11,\,21]$ & $[3,\,14]$  & $22$ & $4$ & $27$ & $761$ & $[-1,\,87]$\\
$\bf 27$ & $[16,\,23]$& $[2,\,18]$   & $24$ & $3$ & $28$ & $802\sspeq 2\cdot 401$ & $[-1,\,95]$\\
         & $[3,\,21]$ & $[7,\,10]$   & $23$ & $4$ & $29$ & $873\sspeq 3^2\cdot 97$ & $[-2,\,90]$ \\
$\bf 28$ & $[3,\,21]$ & $[1,\,4]$   & $12$ & $16$ & $19$ & $873\sspeq 3^2\cdot 97$& $[9,\,57]$ \\   
         & $[8,\,21]$ & $[1,\,9]$   & $17$ & $11$ & $24$ & $818\sspeq 2\cdot 409$& $[4,\,72]$ \\  
         & $[5,\,21]$ & $[4,\,9]$   & $20$ & $8$ & $27$ & $857$ & $[1,\,81]$\\
$\bf 29$ & $[4,\,23]$ & $[8,\,12]$   & $26$ & $3$ & $32$ & $1042\sspeq 2\cdot 521$& $[-3,\,101]$ \\
$\bf 30$ & $[5,\,23]$ & $[6,\,11]$   & $24$ & $6$ & $31$ & $1033$ & $[-1,\,95]$\\
$\bf 31$ & $[19,\,27]$ & $[3,\,22]$  & $29$ & $2$ & $33$ & $1097$ & $[-2,\,114]$\\
$\bf 32$ & $[15,\,27]$ & $[5,\,20]$ & $30$ & $2$ & $35$ & $1233\sspeq 3^\cdot 137$ & $[-3,\,117]$\\
$\bf 33$ & $[3,\,25]$  & $[1,\,4]$  & $13$ & $20$ & $21$ & $1241\sspeq 17\cdot 73$ & $[12,\,64]$\\
         & $[11,\,25]$ & $[1,\,12]$  & $21$ & $12$ & $29$ & $1129$ & $[4,\,88]$\\  
         & $[8,\,25]$  & $[4,\,12]$   & $24$ & $9$ & $32$ & $1186\sspeq 2\cdot 593$ & $[1,\,97]$\\  
         & $[9,\,27]$ & $[8,\,17]$   & $31$ & $2$ & $37$ & $1377\sspeq 3^4\cdot 17$ & $[-4,\,120]$\\ 
$\bf 34$ & $[3,\,25]$ & $[3,\,6]$   & $18$ & $16$ & $27$ & $1241\sspeq 17\cdot 73$& $[7,\,79]$ \\  
         & $[4,\,25]$ & $[3,\,7]$   & $19$ & $15$ & $28$ & $1234\sspeq 2\cdot 617$   & $[6,\,82]$\\
         & $[1,\,25]$  & $[6,\,7]$  & $22$ & $12$ & $31$ & $1249$ & $[3,\,91]$\\
$\bf 35$ & $[13,\,27]$ & $[2,\,15]$   & $25$ & $10$ & $33$ & $1289$ & $[2,\,102]$\\ 
         & $[8,\,27]$  & $[6,\,14]$   & $28$ & $7$ & $36$ & $1394\sspeq 2\cdot 17\cdot 41$ & $[-1,\,111]$\\ 
$\bf 36$ & $[24,\,31]$ & $[1,\,25]$   & $31$ & $5$ & $36$ & $1346 \sspeq 2\cdot 673 $ & $[0,\,124]$\\ 
         & $[8,\,29]$ & $[9,\,17]$   & $33$ & $3$ & $40$ & $1618 \sspeq 2\cdot 809$ & $[-4,\,128]$\\
$\bf 37$ & $[24,\,33]$ & $[4,\,28]$   & $36$ & $1$ & $40$ & $1602 \sspeq 2\cdot 3^2\cdot 89$ & $[-3,\,141]$\\        
$...$ & $...$ & $...$ & $...$ & $...$ & $...$ & $...$& $...$\\ 
\hline
\end{tabular}
\end{center}
\pbn
{\bf Legend for Table 4:\ $\bf q\sspeq \sqrt{c_1\, c_2 \sspp (c_1\sspp c_2)\,c_3}$, and case $\bf q\sspeq c_3 \sspm k$, $\bf k\sspgeq 1$.\pn
Solutions $\bf [X,\, \widehat Y]$ of $\bf X^2 \sspm 2\,\widehat Y^2\sspeq -(t^2 \sspp 2\,k^2)\sspfed -a$, with $\bf 0\sspkl X\sspkl \widehat Y\sspm k$. }\psn
\dstyle{\bf  [c_1,\,c_2]\sspeq \left[\frac{1}{2}(\widehat Y\sspm X\sspm k),\, \frac{1}{2}(\widehat Y\sspp X\sspm k)\right]},\, $\bf c_1^2 \sspp c_2^2 \sspp 6\,c_1\,c_2 \sspp 4\,k\,(c_1\sspp c_2)\sspeq t^2$.\pn
$\bf c_3 \sspeq (\widehat Y\sspp k\sspp t)/2$,\, $\bf q \sspeq (\widehat Y\sspm k\sspp t)/2$.
\vfill
\eject
\begin{center}
{\large {\bf Table 5 Case ii) $\bf q\sspgr c_3$}}.\\ 
\end {center}
\begin{center}
\begin{tabular}{|c|l|l|c|c|c|l|l|}\hline
&& && && &\\
\hskip .2cm
$\bf c_3$ & $\bf [X,\,\widehat Y]$ & $\bf[c_1,c_2]$ & $\bf q$ &$\bf k$ & $\bf t$& $\bf a$ & $\bf [c_{4,-},\,c_{4,+}]$\\
&& && &&\\ 
\hline\hline
$\bf 7$ & $[3,\,7]$ & $[3,\,6]$ & $9$ & $2$ & $9$ & $89$ & $[-2,\,34]$ \\
$\bf 10$ & $[4,\,9]$ & $[3,\,7]$ & $11$ & $1$ & $12$ & $146\sspeq 2\cdot 73$ & $[-2,\,42]$  \\
$\bf 12$ & $[3,\,11]$& $[5,\,8]$& $14$ & $2$ & $15$ & $233$ & $ [-3,\,53]$ \\
$\bf 13$ & $[8,\,13]$ & $[4,\,12]$ & $16$ & $3$ & $16$ & $274\sspeq 2\cdot 137$ & $[-3,\,61]$\\
$\bf 15$ & $[3,\,17]$& $[11,\,14]$    & $23$ & $8$ & $21$ & $569$ & $[-6,\,86]$\\
$\bf 16$ & $[9,\,15]$ & $[4,\,13]$   & $18$ & $2$ & $19$ & $369\sspeq 3^2\cdot 41$ & $[-3,\,69]$\\
$\bf 17$ & $[1,\,15]$ & $[8,\,9]$   & $19$ & $2$ & $21$ & $449$ &$[-4,\,72]$ \\
$\bf 19$ & $[5,\,19]$ & $[10,\,15]$   & $25$ & $6$ & $25$ & $697\sspeq 17\cdot 41$ & $[-6,\,94]$\\
$\bf 20$ & $[5,\,21]$ & $[12,\,17]$  & $28$ & $8$ & $27$ & $857$& $[-7,\,105]$ \\
$\bf 21$ & $[12,\,19]$ & $[4,\,16]$   & $22$ & $1$ & $24$ & $578\sspeq 2\cdot 17^2$ & $[-3,\,85]$\\
         & $[15,\,21]$ & $[5,\,20]$   & $25$ & $4$ & $25$ & $657\sspeq 3^2\cdot 73$ & $[-4,\,96]$\\
$\bf 22$ & $[11,\,21]$& $[7,\,18]$   & $26$ & $4$ & $27$ & $761$ & $[-5,\,99]$\\ 
           & $[1,\,25]$ & $[18,\,19]$ & $34$ & $12$ & $31$ & $1249$ & $[-9,\,127]$\\
$\bf 23$ & $[3,\,21]$& $[11,\,14]$   & $27$ & $4$ & $29$ & $873\sspeq 3^2\cdot 97$ & $[-6,\,102]$ \\
$\bf 24$ & $[16,\,23]$ & $[5,\,21]$  & $27$ & $3$ & $28$ & $802\sspeq 2\cdot 401$ & $[-4,\,104]$\\
         & $[5,\,23]$ & $[12,\,17]$ & $30$ & $6$ & $31$ & $1033$ & $[-7,\,113]$\\ 
         & $[8,\,25]$ & $[13,\,21]$ & $33$ & $9$ & $32$ & $1186\sspeq 2\cdot 593$ & $[-8,\,124]$\\
$\bf 26$ & $[4,\,23]$& $[11,\,15]$   & $29$ & $3$ & $32$ & $1042\sspeq 2\cdot 521$ & $[-6,\,110]$\\
$\bf 27$ & $[12,\,29]$ & $[14,\,26]$ & $38$ & $11$ & $36$ & $1538\sspeq 2\cdot 769$& $[-9,\,143]$ \\   
         & $[5,\,29]$ & $[18,\,23]$ & $39$ & $12$ & $37$ & $1657$ & $[-10,\,146]$\\  
$\bf 28$ & $[8,\,27]$ & $[13,\,21]$ & $35$ & $7$ & $36$ & $1394\sspeq 2\cdot 17\cdot 41$ & $[-8,\,132]$\\ 
         & $[3,\,31]$ & $[21,\,24]$ & $42$ & $14$ & $39$ & $1913$ & $[-11,\,157]$\\
$\bf 29$ & $[19,\,27]$ & $[5,\,24]$ & $31$ & $2$ & $33$ & $1097$ & $[-4,\,120]$\\
$\bf 30$ & $[15,\,27]$ & $[7,\,22]$ & $32$ & $2$ & $35$ & $1233\sspeq 3^2\cdot 137$ & $[-5,\,123]$\\
$\bf 31$ & $[24,\,31]$ & $[6,\,30]$ & $36$ & $5$ & $36$ & $1346\sspeq 2\cdot 673$ & $[-5,\,139]$\\
         & $[9,\,27]$ & $[10,\,19]$ & $33$ & $2$ & $37$ & $1377\sspeq 3^4\cdot 17$ & $[-6,\,126]$\\
$\bf 32$ & $[8,\,35]$  & $[21,\,29]$  & $47$ & $15$ & $44$ & $2386\sspeq 2\cdot 1193$& $[-12,\,176]$ \\
$\bf 33$ & $[8,\,29]$ & $[12,\,20]$   & $36$ & $3$ & $40$ & $1618\sspeq 2\cdot 809$ &$[-2,\,137]$\\  
         & $[13,\,31]$ & $[12,\,25]$   & $39$ & $6$ & $41$ & $1753$ & $[-8,\,148]$\\  
$\bf 34$ & $[25,\,33]$ & $[6,\,31]$   & $38$ & $4$ & $39$ & $1533$ & $[-5,\,147]$\\ 
         & $[1,\,31]$  & $[18,\,19]$   & $40$ & $6$ & $43$ & $1921\sspeq 17\cdot 113$ & $[-9,\,151]$\\ 
         & $[4,\,39]$  & $[27,\,31]$   & $53$ & $19$ & $48$ & $3026\sspeq 2\cdot 17\cdot 89$ & $[-14,\,198]$\\
$\bf 35$ & $[12,\,33]$ & $[14,\,26]$   & $42$ & $7$ & $44$ & $2034 \sspeq 2\cdot 3^2\cdot 113 $ & $[-9,\,159]$\\ 
         & $[5,\,33]$ & $[18,\,23]$   & $43$ & $8$ & $45$ & $2153$ & $[-10,\,162]$\\
$...$ & $...$ & $...$ & $...$ & $...$ & $...$ & $...$& $...$\\ 
\hline
\end{tabular}
\end{center}
\pbn
{\bf Legend for Table 5:\ $\bf q\sspeq \sqrt{c_1\, c_2 \sspp (c_1\sspp c_2)\,c_3}$, and case $\bf q\sspeq c_3 \sspp k$, $\bf k\sspgeq 1$.\pn
Solutions $\bf [X,\, \widehat Y]$ of $\bf X^2 \sspm 2\,\widehat Y^2\sspeq -(t^2 \sspp 2\,k^2)\sspfed -a$, with $\bf 0\sspkl X\sspkl \widehat Y\sspp k$. }\psn
\dstyle{\bf  [c_1,\,c_2]\sspeq \left[\frac{1}{2}(\widehat Y\sspm X\sspp k),\, \frac{1}{2}(\widehat Y\sspp X\sspp k)\right]},\, $\bf c_1^2 \sspp c_2^2 \sspp 6\,c_1\,c_2 \sspm 4\,k\,(c_1\sspp c_2)\sspeq t^2$.\pn
$\bf c_3 \sspeq (\widehat Y\sspm k\sspp t)/2$,\, $\bf q \sspeq (\widehat Y\sspp k\sspp t)/2$.
\vfill
\eject
\begin{center}
{\large {\bf Table 6: Primitive Cartesian triples with $\bf c_{4,-} \sspeq 0$ and $\bf c_3\sspeq n^2, n\sspeq 2 .. 24$}} 
\end {center}
\begin{center}  
\begin{tabular}{|l|c|c|c||l|c|c|c|}\hline
&& && && &\\
\hskip .2cm$\bf n$ & $\bf [c_1,\, c_2,\, c_3]$  &$\bf q $ & $\bf [c_{4,-},\,c_{4,+}]$& $\bf n$ & $\bf [c_1,c_2,c_3]$  &$\bf q $ & $\bf [c_{4,-},c_{4,+}]$  \\
&& && && &\\ \hline\hline
$\bf 2$ & $[1, 1, 4]$ & $3$ & $[0, 12] $ &                       $\bf 17cont'd$ & $[9, 196, 289]$ & $247$ & $ [0, 988] $ \\  
$\bf 3$ & $[1, 4, 9]$ & $7$ & $[0, 28] $ &                       $\bf  $& $[4, 225, 289 ]$ & $  259$ & $[0, 1036] $\\ 
$\bf 4$ & $  [1, 9, 16] $ & $13$ & $[0, 52] $ &                $\bf       $& $[1, 256, 289]$ & 273$ $ & $[0, 1092] $\\ 
$\bf 5$ & $[4, 9, 25]$ & $19$ & $[0,76]  $ &                   $\bf 18 $& $[49, 121, 324]$ & $247 $ & $[0, 988]$  \\
$\bf  $ & $[1, 16, 25]$ & $21 $ & $[0, 84] $ &                 $\bf      $& $[25, 169, 324]$ & $ 259$ & $[0, 1036]$ \\
$\bf 6$& $[1, 25, 36 ]$ & $31$ & $[0, 124] $ &                $\bf      $ & $[1, 289, 324]$ & $307$ & $[0, 1228] $\\
$\bf 7 $& $[9, 16, 49]$ & $37$ & $[0, 148] $ &                $\bf 19 $ & $[81, 100, 361 ]$ & 271$$ & $[0, 1084$ \\
$\bf   $& $[4, 25, 49]$ & $39$ & $[0, 156]$ &                  $\bf      $ & $[64, 121, 361]$ & $273$ & $[0,  1092] $ \\
$\bf   $& $[1, 36, 49]$ & $43$ & $[0, 172]$ &                 $\bf    $ & $[49, 144, 361]$ & $ 277$ & $[0, 1108] $ \\
$\bf 8 $& $[9, 25, 64]$ & $49$ & $[0, 196] $ &                $\bf     $ & $[36, 169, 361 ]$ & $283 $ & $[0, 1132]$ \\
$\bf   $& $[1, 49,  64]$ & $57 $ & $[0, 228] $ &               $\bf     $ & $[25, 196, 361]$ & $291$ & $[0, 1164] $ \\
$\bf 9$& $[16, 25, 81]$ & $61$ & $[0, 244]$ &                 $\bf     $ & $[16, 225, 361]$ & $ 301$ & $[0, 1204] $ \\
$\bf     $& $[4, 49, 81]$ & $67$ & $[0, 268]$ &                $\bf    $ & $[9, 256, 361]$ & $313$ & $[0, 1252] $ \\
$\bf     $& $[1, 64, 81]$ & $73$ & $[0, 292] $&                 $\bf    $ & $[4, 289, 361]$ & $327$ & $[0, 1308] $ \\
$\bf 10 $& $[9, 49, 100 ]$ & $79 $ & $[0, 316] $ &            $\bf   $ & $[1, 324, 361]$ & $343$ & $[0, 1372] $ \\
$\bf   $& $[1, 81, 100]$ & $91$ & $[0,364] $ &                  $\bf  20   $ & $[81, 121, 400]$ & $301$ & $[0, 1204] $ \\  
$\bf 11$& $[25, 36, 121]$ & $ 91$ & $[0,364] $ &              $\bf      $ & $[49, 169, 400]$ & $309$ & $[0, 1236] $ \\
$\bf    $& $[16, 49, 121 ]$ & $93 $ & $[0, 372] $ &             $\bf     $ & $[9, 289, 400]$ & $349$ & $[0, 1396] $ \\
$\bf    $& $[9, 64, 121]$ & $ 97$ & $[0, 388] $ &                $\bf     $ & $[1, 361, 400]$ & $381$ & $[0, 1524] $ \\
$\bf    $& $ [4, 81, 121]$ & $ 103 $ & $[0, 412] $ &             $\bf 21$ & $[100, 121, 441]$ & $331$ & $[0, 1324] $ \\
$\bf    $& $[1, 100, 121]$ & $111$ & $[0, 444] $ &               $\bf      $ & $[64, 169, 441]$ & $337$ & $[0, 1348] $ \\
$\bf 12 $& $[25, 49, 144]$ & $109$ & $[0, 436] $ &              $\bf      $ & $[25, 256, 441]$ & $ 361$ & $[0, 1444] $ \\   
$\bf      $& $[1, 121, 144]$ & $133$ & $[0, 532] $ &              $\bf      $ & $[16, 289, 441]$ & $ 373 $ & $[0, 1492] $ \\  
$\bf 13  $& $[36, 49, 169]$ & $127 $ & $[0, 508] $ &            $\bf     $ & $[4, 361, 441]$ & $ 403$ & $[0, 1612] $ \\   
$\bf    $& $[25, 64, 169]$ & $129$ & $[0, 516] $ &               $\bf     $ & $[1, 400, 441]$ & $421 $ & $[0, 1684] $ \\  
$\bf     $& $[16, 81, 169]$ & $133$ & $[0, 532] $ &               $\bf 22 $ & $[81, 169, 484]$ & $367$ & $[0, 1468] $ \\
$\bf    $& $[9, 100, 169]$ & $139 $ & $[0, 556] $ &               $\bf     $ & $[49, 225, 484]$ & $379$ & $[0, 1516] $ \\
$\bf     $& $[4, 121, 169 ]$ & $147$ & $[0, 588] $ &             $ \bf     $ & $[25, 289, 484]$ & $399$ & $[0, 1596] $ \\
$\bf     $& $[1, 144, 169 ]$ & $157$ & $[0, 628] $ &             $\bf     $ & $[9, 361, 484]$ & $ 427$ & $[0, 1708] $ \\
$\bf 14 $& $[25, 81, 196]$ & $151$ & $[0, 604] $ &             $\bf   $ & $[1, 441, 484]$ & $ 463$ & $[0, 1852] $ \\
$\bf     $& $[9, 121, 196]$ & $163$ & $[0, 652] $ &              $\bf 23$ & $[121, 144, 529]$ & $397$ & $[0, 1588] $ \\ 
$\bf     $& $[1, 169, 196]$ & $183$ & $[0, 732] $ &              $\bf     $ & $[100, 169, 529]$ & $399$ & $[0, 1596] $ \\ 
$\bf 15 $& $[49, 64, 225]$ & $169$ & $[0, 676] $ &             $\bf     $ & $[81, 196, 529]$ & $ 403$ & $[0, 1612] $ \\
$\bf     $& $[16, 121, 225]$ & $181$ & $[0, 724] $ &             $\bf     $ & $[64, 225, 529]$ & $ 409$ & $[0, 1636] $ \\
$\bf     $& $[4, 169, 225]$ & $ 199$ & $[0, 796] $ &                $\bf    $ & $[49, 256, 529 ]$ & $417$ & $[0, 1668] $ \\ 
$\bf     $& $[1, 196, 225]$ & $211$ & $[0, 844] $ &                 $\bf       $ & $[36, 289, 529 ]$ & $427$ & $[0, 1708] $ \\     
$\bf 16 $& $[49, 81, 256]$ & $193$ & $[0, 772] $ &               $\bf       $ & $[25, 324, 529]$ & $439$ & $[0, 1756] $ \\ 
$\bf     $& $[25, 121, 256]$ & $201$ & $[0, 804] $ &              $\bf      $ & $[16, 361, 529 ]$ & $453$ & $[0,  1812] $ \\
$\bf     $& $[9, 169, 256]$ & $217$ & $[0,  868] $ &               $\bf  $ & $[9, 400, 529]$ & $469$ & $[0,  1876] $ \\
$\bf     $& $[1, 225, 256]$ & $241$ & $ [0, 964] $ &               $\bf       $ & $[ 4, 441, 529]$ & $487$ & $[0, 1948] $ \\
$\bf 17 $& $[64, 81, 289]$ & $217 $ & $ [0, 868] $ &               $\bf  $ & $[1, 484, 529]$ & $ 507$ & $[0, 2028] $ \\
$\bf     $& $[49, 100, 289]$ & $219 $ & $ [0, 876] $ &              $\bf 24 $ & $[121, 169, 576]$ & $433$  & $[0, 1732] $ \\
$\bf     $& $[36, 121, 289 ]$ & $223 $ & $ [0,  892] $ &           $\bf  $ & $[49, 289, 576]$ & $457$  & $[0, 1828] $ \\   
$\bf     $& $[25, 144, 289]$ & $ 229$ & $ [0, 916] $ &              $\bf  $ & $[25, 361, 576]$ & $481$  & $[0, 1924 ] $ \\  
$\bf     $& $[16, 169, 289]$ & $237 $ & $ [0, 948] $ &             $\bf  $ & $[1, 529, 576]$ & $553$  & $[0, 2212] $ \\   
\hline
\end{tabular}
\end{center}
\vfill
\eject

\end{document}